\def\no{\if01}
\def\iftwelvept{\no}

\def\ifusepdf{\no}
\def\ifpsfont{\no}

\iftwelvept
\documentclass[leqno,12pt]{amsart}
\else
\documentclass[leqno,10pt]{amsart}
\fi
\usepackage{amssymb}
\usepackage{amscd}
\usepackage{latexsym}
\usepackage{verbatim}
\usepackage[all]{xy}

\ifusepdf
\usepackage{hyperref}
\else\fi
\ifpsfont
\usepackage[T1]{fontenc}
\usepackage{times}
\else\fi

\iftwelvept
\else\fi

\theoremstyle{plain}
\newtheorem{Theorem}{Theorem}[section]
\newtheorem{Warning}{Warning}[Theorem]

\newtheorem{Proposition}[Theorem]{Proposition}
\newtheorem{Lemma}[Theorem]{Lemma}
\newtheorem{Corollary}[Theorem]{Corollary}

\theoremstyle{definition}

\newtheorem{Definition}[Theorem]{Definition}
\newtheorem{Remark}[Theorem]{Remark}
\newtheorem{Construction}[Theorem]{Construction}
\newtheorem{Example}[Theorem]{Example}

\newtheorem{Notation}[Theorem]{Notation}

\newcommand{\ZZ}{\mathbf{Z}}

\newcommand{\DD}{\mathbb{D}}

\newcommand{\NNNN}{\operatorname{N}}

\newcommand{\HH}{\operatorname{\mathcal{HH}}}

\newcommand{\EE}{\mathcal{E}}

\newcommand{\AAA}{\mathcal{A}}
\newcommand{\BBB}{\mathcal{B}}

\newcommand{\DDD}{\mathcal{D}}
\newcommand{\uni}{\mathbf{1}}
\newcommand{\CCC}{\mathcal{C}}

\newcommand{\PR}{\textup{Pr}}
\newcommand{\PRT}{\textup{Pr}^{\textup{t}}}
\newcommand{\PRTT}{\textup{Pr}^{\textup{t}+}}
\newcommand{\PRTTT}{\textup{Pr}^{\textup{t}\pm}}

\newcommand{\OO}{{\mathcal{O}}}

\newcommand{\LLL}{\mathcal{L}}

\newcommand{\MMM}{\mathcal{M}}

\newcommand{\PPP}{\mathcal{P}}

\newcommand{\Hom}{\operatorname{Hom}}

\newcommand{\Spec}{\operatorname{Spec}}

\newcommand{\SP}{\operatorname{Sp}}

\newcommand{\SPC}{\operatorname{Sp}^{\textup{conn}}}
\newcommand{\Mod}{\operatorname{Mod}}

\newcommand{\SSS}{\mathcal{S}}
\newcommand{\wSSS}{\widehat{\mathcal{S}}}

\newcommand{\colim}{\operatorname{colim}}

\newcommand{\Cat}{\textup{Cat}_{\infty}}

\newcommand{\Map}{\operatorname{Map}}

\newcommand{\Fun}{\operatorname{Fun}}
\newcommand{\Alg}{\operatorname{Alg}}

\newcommand{\End}{\operatorname{End}}

\newcommand{\wCat}{\widehat{\textup{Cat}}_{\infty}}

\newcommand{\CAlg}{\operatorname{CAlg}}

\newcommand{\vCat}{\widehat{\operatorname{Cat}}_\dagger}

\newcommand{\QC}{\operatorname{QC}}

\newcommand{\LM}{\operatorname{\mathcal{LM}}}

\newcommand{\hhh}{\operatorname{h}}
\newcommand{\Ind}{\operatorname{Ind}}

\newcommand{\XX}{\mathcal{X}}

\newcommand{\PRST}{\PR_{\textup{St}}}

\newcommand{\assoc}{\operatorname{As}}

\newcommand{\eone}{\mathbf{E}_1}
\newcommand{\etwo}{\mathbf{E}_2}
\newcommand{\eenu}{\mathbf{E}_n}
\newcommand{\einf}{\mathbf{E}_\infty}

\newcommand{\LMod}{\operatorname{LMod}}
\newcommand{\RMod}{\operatorname{RMod}}
\newcommand{\BMod}{\operatorname{BMod}}

\newcommand{\HL}{\widehat{L}}
\newcommand{\UUUU}{\mathbb{U}}
\newcommand{\VVVV}{\mathbb{V}}

\newcommand{\Sym}{\operatorname{Sym}}
\newcommand{\Proof}{{\sl Proof.}\quad}
\newcommand{\QED}{{\unskip\nobreak\hfil\penalty50\quad\null\nobreak\hfil
{$\Box$}\parfillskip0pt\finalhyphendemerits0\par\medskip}}

\begin{document}

\title{Categorified Koszul duality and Deformation theory of higher categories}

\author{Isamu Iwanari}


\maketitle

\setcounter{tocdepth}{1}
\tableofcontents

\section{Introduction}

In this paper, we develop a deformation theory of higher categories, which provides a generalization
of the deformation theory in the context of stable $(\infty,1)$-categories to the setting of stable $(\infty,n)$-categories for $n>1$.
Informally, 
an $(\infty,n)$-category is a higher category which possesses 
$i$-morphisms between $(i-1)$-morphisms for every $i\ge1$
($0$-morphisms are objects)
such that $i$-morphisms are invertible for $i> n$. The $(\infty,n)$-category theory is extensively studied from various respects; see e.g. \cite{BSP}, \cite{GH}, \cite{JFS}, \cite{Lu}, \cite{Si}.
To explain the purpose of this paper, let us briefly review the deformation theory of stable $(\infty,1)$-categories (that is, stable $\infty$-categories) with an emphasis on its difficulties.
Let $\DDD$ be a $k$-linear stable $\infty$-category over a field $k$ (or a $k$-linear dg category), and let $R$ be an Artin algebra equipped with a homomorphism $R\to k$. 
The deformation theory of $\DDD$ studies deformations (or liftings) $(\DDD_R,\DDD_R\otimes_Rk\simeq \DDD)$ of $\DDD$ to $R$, where $\DDD_R$ is an $R$-linear stable $\infty$-category.
These are analyzed
by means of algebraic invariants such as Hochschild cohomology $\HH^\bullet(\DDD/k)$.
It is well known that unlike the deformation theory of associative algebras \cite{Ger1}, \cite{Ger2} and more general $\mathbf{E}_{n}$-algebras (algebras over little $n$-disk operads) \cite{FR}, the deformation theory of plain $k$-linear stable $\infty$-categories is not well-behaved.
Indeed, the deformation problems of plain $k$-linear stable $\infty$-categories are in general not formal moduli problems, and
they often cannot be described in terms of Hochschild cohomology: a typical example is the derived $\infty$-category of coherent sheaves on $\Spec k[x]/(x^2)$.
Consequently, a good deformation theory essentially requires refined data.
In \cite{IH}, we developed a deformation theory of $k$-linear stable $\infty$-categories equipped with left complete $t$-structures
and showed that this deformation theory works well: the deformation problems are formal moduli problems, and they are described by means of Hochschild
cohomology of the underlying stable $\infty$-categories.
We turn to the setting of higher categories, that is, stable $(\infty,n)$-categories.
The above difficulties persist in this context. Even worse, it was shown in \cite{Ch}, using an argument of \cite[X]{DAG},
that the deformation theory of plain stable $(\infty,n)$-categories is in general only an $(n+1)$-proximated formal moduli problem.
Roughly speaking, this means that the larger $n$ is, the further the deformation problem is from being a formal moduli problem, and hence the larger the a priori discrepancy between the deformation problem and the Hochschild cohomology.
We are therefore led to the following question:
\begin{center}
What is an appropriate framework in which the deformation problem for $(\infty,n)$-categories becomes well-behaved?
\end{center}
To address this question,
we introduce the notion of a $t$-structured stable presentable $(\infty,n)$-category. We then develop a deformation theory of complete $t$-structured $k$-linear $(\infty,n)$-categories by establishing a categorified Koszul duality for them.

{\it Categorified Koszul duality and $t$-structured stable presentable higher categories.}
We now consider the categorified Koszul duality from the viewpoint of deformation theory.
If $A$ is an ordinary associative $k$-algebra (more generally, a dg algebra satisfying suitable conditions),
the deformation theory of $A$-modules can be understood via a classical Koszul duality (cf. \cite[X, Section 5.2]{DAG}):
for an Artin local $k$-algebra $R$ with $R\to k$ and its Koszul dual $\DD_1(R)=\End_{R}(k)$, Koszul duality gives
an equivalence between right $R$-modules and left $\DD_1(R)$-modules in a suitable formulation.
The Koszul duality for $k$-linear stable $\infty$-categories generalizes this equivalence to
the following one
\[
\xymatrix{
\textup{right}\ \LMod_{R}^\otimes\textup{-module objects in}\ \PRTTT_k \ar[r]^(0.5){\simeq} & \textup{left}\ \LMod_{\DD_2(R)}^\otimes\textup{-module objects in}\ \PRTTT_k  \ar[l]
}
\]
where $\PRTTT_k$ is the $\infty$-category of $k$-linear stable presentable $\infty$-categories equipped with both left and right complete $t$-structures,
and $\LMod_R^\otimes$ (resp. $\LMod_{\DD_2(R)}^\otimes$) is the monoidal $\infty$-category of left $R$-modules (resp. left $\DD_2(R)$-modules)
with the right complete $t$-structure. Here $\DD_2(R)$ is the $\etwo$-Koszul dual of $R$.
The $t$-structures play the role of filtrations (or topologies) in classical Koszul duality, and $t$-structures are essential for realizing
the categorified Koszul duality (see \cite{IH} for details).
We turn to $k$-linear stable $(\infty,n)$-categories. We use the notion of a presentable $(\infty,n)$-category introduced in \cite{Ste}, 
see also \cite{Ao} for a useful reformulation. Let us focus on $(\infty,2)$-categories.
Roughly speaking, a $k$-linear (stable) presentable $(\infty,2)$-category is an $\infty$-category $\DDD$ admitting small colimits and equipped with
a $\PR_k$-module structure $\PR_k\times \DDD\to \DDD$ satisfying a certain finiteness condition, where $\PR_k$ is the symmetric monoidal $\infty$-category of $k$-linear presentable $\infty$-categories.
A complete $t$-structured $k$-linear (stable) presentable $(\infty,2)$-category is an $\infty$-category $\DDD$ admitting small colimits and equipped with a $\PRTTT_k$-module structure $\PRTTT_k\times \DDD\to \DDD$ satisfying a certain finiteness condition.
The notion of a (complete) $t$-structured $k$-linear (stable) presentable $(\infty,n)$-category for arbitrary $n$ is defined along the same lines.
We prove
a higher categorified Koszul duality (i.e., Koszul duality for $(\infty,n)$-categories), stated informally as follows (see Theorem~\ref{HCKD}):

\begin{Theorem}
\label{intromain1}
Let $n\ge1$ and let $R$ be an almost Artin $\mathbf{E}_{n+1}$-algebra. Let $\DD_{n+1}(R)$ denote the $\mathbf{E}_{n+1}$-Koszul dual.
Then there exists an equivalence
\[
\xymatrix{
\textup{right}\ (n-1)\PRTTT_{R}\textup{-module objects in}\ n\PRTTT_k \ar[r]^(0.48){\simeq} & \textup{left}\ (n-1)\PRTTT_{\DD_{n+1}(R)}\textup{-module objects in}\ n\PRTTT_k  \ar[l]
}
\]
\end{Theorem}

Roughly, for a complete $t$-structured
stable presentable $(\infty,n)$-category $\CCC$, a (right) $(n-1)\PRTTT_{R}$-module structure on $\CCC$ provides a (right) $R$-linear structure on $\CCC$.
Similarly, a (left) $(n-1)\PRTTT_{\DD_{n+1}(R)}$-module structure gives a (left) $\DD_{n+1}(R)$-linear structure.
Thus, Theorem~\ref{intromain1} gives an equivalence between right $R$-linear complete $t$-structured presentable $(\infty,n)$-categories and left $\DD_{n+1}(R)$-linear complete $t$-structured presentable $(\infty,n)$-categories.
The symbol $n\PRTTT_k$ denotes the $\infty$-category of complete $t$-structured $k$-linear presentable $(\infty,n)$-categories,
and $(n-1)\PRTTT_{R}$ and $(n-1)\PRTTT_{\DD_{n+1}(R)}$ are higher module categories: 
for example, $0\PRTTT_{R}$ is the $\LMod_R$ with $t$-structure determined by the full subcategory of connective objects,
$1\PRTTT_{R}=\LMod_{\LMod_R}(\PRTTT_k)$, $2\PRTTT_R=\LMod_{1\PRTTT_{R}}(2\PRTTT_k)$ and so on.
Here $\LMod_{X}(P)$ denotes the $\infty$-category of left $X$-module object in $P$, and $2\PRTTT_k$ is the symmetric monoidal category of complete $t$-structured $k$-linear presentable $(\infty,2)$-categories (see Section~\ref{PRA}). When $n=1$, the result is the one above.

\vspace{3mm}

{\it Deformation theory.}
We will give an application of this categorified Koszul duality to the deformation theory of $(\infty,n)$-categorical context.
Let $\DDD$ be a complete $t$-structured $k$-linear presentable $(\infty,n)$-category, that is, an object of $n\PRTTT_k$.
For an Artin $\mathbf{E}_{n+1}$-algebra $R$, a deformation of $\DDD$ to $R$
is a right $(n-1)\PRTTT_R$-module $\DDD_R$ in $n\PRTTT_k$, together with an equivalence $\phi:\DDD_R\otimes_{(n-1)\PRTTT_R}(n-1)\PRTTT_k\simeq \DDD$,
where the tensor product is taken along the functor $(n-1)\PRTTT_R\to(n-1)\PRTTT_k$ induced by $R\to k$.
We prove the following (see Theorem~\ref{deformmain}, Theorem~\ref{FMPmain} for details):

\begin{Theorem}
\label{intromain2}
Suppose that $\DDD$ has Hochschild cohomology $\HH^\bullet(\DDD/k)$ (see Section~\ref{hochsec}).
It is then naturally an $\mathbf{E}_{n+1}$-algebra (cf. the higher Deligne conjecture).
Let $\textup{Deform}_{\DDD}(R)$ be the space of deformations of $\DDD$ to $R$. 
The functor
$\textup{Deform}_{\DDD}$ from the $\infty$-category of Artin (dg) $\mathbf{E}_{n+1}$-algebras to that of spaces is a formal $\mathbf{E}_{n+1}$-moduli problem determined by
Hochschild cohomology $\HH^\bullet(\DDD/k)$ (see Theorem~\ref{FMPmain} for the formulation).

For any Artin $\mathbf{E}_{n+1}$-algebra $R$, there exists a canonical equivalence of spaces
\[
\textup{Deform}_{\DDD}(R)\simeq \Map_{\Alg_{n+1}(\Mod_k)}(\DD_{n+1}(R),\HH^\bullet(\DDD/k))
\]
where the right-hand side is the mapping space of $\mathbf{E}_{n+1}$-algebras.
\end{Theorem}

In characteristic zero, if we restrict the deformation problem to commutative Artin (dg) algebras, it becomes a formal moduli problem
determined by the underlying dg Lie algebra $\HH^\bullet(\DDD/k)[n]$.
We also deduce the obstruction theory in this context (see Theorem~\ref{obstructionmain}):

\begin{Theorem}
\label{intromain3}
Let $i\ge1$.
Let $p:R'\to R$ be an elementary extension of $R$ by an $(i-1)$-shifted $k$-vector space $V[i-1]$ (see Section~\ref{obssec}).
Let $(\DDD_R,\phi)$ be a deformation of $\DDD$ to $R$.
Let $\HH^{i+n+1}(\DDD/k)$ denote $H^{i+n+1}(\HH^\bullet(\DDD/k))$.
There exists a canonically defined element $\textup{Obs}(\DDD_R,p)\in \HH^{i+n+1}(\DDD/k)\otimes_kV$ 
such that $\textup{Obs}(\DDD_R,p)$ vanishes if and only if there exists $\DDD_{R'}\in \pi_0(\textup{Deform}_{\DDD}(R'))$ lying over $\DDD_R$ along  the map $p$.
\end{Theorem}

{\it Hochschild cohomology.}
The Hochschild cohomology $\HH^\bullet(\DDD/k)$ appearing in Theorems~\ref{intromain2} and~\ref{intromain3} is defined from the $t$-structured presentable
$(\infty,n)$-category $\DDD$, and it is defined as a higher version of Hochschild cohomology.  
It depends a priori on the $t$-structure on $\DDD$. We also note that, as the first assumption of Theorem~\ref{intromain2} suggests, the existence of the Hochschild cohomology of a presentable $(\infty,n)$-category is a nontrivial problem even without $t$-structures. We prove that the Hochschild cohomology of $t$-structured $k$-linear presentable $(\infty,n)$-categories exists in several important cases. Moreover, in these cases it depends only on the underlying (non-$t$-structured) $k$-linear presentable $(\infty,n)$-category.
We discuss Hochschild cohomology in detail in Section~\ref{hochsec}.

\vspace{2mm}

{\it Application to sheaves of categories.}
We give an application of our deformation theory to sheaves of categories on stacks and their higher generalization.
For a scheme $\XX$ over $k$ (or, more generally, for an algebraic stack satisfying suitable conditions),
the notion of a sheaf of categories on $\XX$ is obtained by categorifying the notion of a quasi-coherent sheaf. Informally, a sheaf of categories $\MMM$ assigns an $R$-linear presentable $\infty$-category $\Gamma(U,\MMM)$
to each $U=\Spec R\to \XX$ (see \cite{Gai1aff}, \cite{TV}). Sheaves of categories on $\XX$ constitute an $(\infty,2)$-category, and by
$1$-affineness of $\XX$ \cite{Gai1aff} (under a suitable condition) this $(\infty,2)$-category is given by the module category $\Mod_{\QC(\XX)}(\PR_k)$. 
Here $\QC(\XX)$ denotes the monoidal $\infty$-category of quasi-coherent sheaves (complexes) on $\XX$. It has a canonical enhancement to a complete $t$-structured presentable $(\infty,2)$-category
$\Mod_{\QC(\XX)}(\PRTTT_k)$ (see Section~\ref{catsheaves}). This can be generalized to the context of sheaves of $(\infty,n)$-categories, and we define
a complete $t$-structured presentable $(\infty,n+1)$-category $n\Mod_{\QC(\mathcal{X})}[\PRTTT_k]$ accordingly (see Section~\ref{highercatsheaves}).
This is a higher categorical generalization of $\Mod_{\QC(\XX)}(\PRTTT_k)$ in the sense that $1\Mod_{\QC(\mathcal{X})}[\PRTTT_k]=\Mod_{\QC(\XX)}(\PRTTT_k)$.
We prove the following (see Proposition~\ref{catsheafhoch}, Proposition~\ref{highercatprop}, Theorem~\ref{higherfinal} for details),
which provides a deformation-theoretic interpretation of the higher Hochschild cohomology of $\QC(\XX)$.
Note that $n\Mod_{\QC(\mathcal{X})}[\PRTTT_k]$ is an $(\infty,n+1)$-category, so Theorems~\ref{intromain2} and~\ref{intromain3} apply with $n$ replaced by $n+1$:

\begin{Corollary}
\label{intromain4}
Let $\XX$ be a perfect stack over $k$. 
Consider the deformation functor $\textup{Deform}_{n\Mod_{\QC(\mathcal{X})}[\PRTTT_k]}$.
Then $\textup{Deform}_{n\Mod_{\QC(\mathcal{X})}[\PRTTT_k]}$ is an $\mathbf{E}_{n+2}$-formal moduli problem
determined by Hochschild cohomology $\HH^\bullet(n\Mod_{\QC(\mathcal{X})}[\PRTTT_k]/k)$. Moreover, we have
\begin{enumerate}

\item $\HH^\bullet(n\Mod_{\QC(\mathcal{X})}[\PRTTT_k]/k)$ is equivalent to the Hom complex $\mathcal{H}om_{\QC(\Map(S^n,\mathcal{X}))}(\OO_{\Delta}, \OO_{\Delta})$
in $\Mod_k$. Here $\Map(S^n,\mathcal{X})$ is the iterated free loop space, that is, the mapping stack from the sphere $S^n$ to $\XX$ in the framework of derived geometry, and
$\OO_{\Delta}$ is the pushforward of $\OO_{\XX}$ along the morphism $\XX\to\Map(S^n,\XX)$ given by constant loops.

\item 
We further assume that $\XX$ is a scheme separated of finite type over $k$ of characteristic zero. 
Then there is a Koszul duality type isomorphism
\[
\HH^\bullet(n\Mod_{\QC(\mathcal{X})}[\PRTTT_k]/k)\simeq \Gamma(\XX,\Sym_{\OO_{\XX}}^{\bullet}(\mathbb{T}_{\XX/k}[-n-1]))
\]
where the right-hand side denotes the derived global sections.

\item 
The isomorphism classes of deformations of $n\Mod_{\QC(\mathcal{X})}[\PRTTT_k]$ to the trivial square zero extension $k\oplus k[j]$ ($j\ge0$) (i.e. $\pi_0(\textup{Deform}_{n\Mod_{\QC(\mathcal{X})}[\PRTTT_k]})(k\oplus k[j])$)
is isomorphic to
\[
H^{j+n+2}(\mathcal{H}om_{\QC(\Map(S^n,\mathcal{X}))}(\OO_{\Delta}, \OO_{\Delta})).
\]

\item if $p:R'\to R$ is an elementary extension of $R$ by a shifted $k$-vector space $V[i-1]$
and $(n\Mod_{\QC(\mathcal{X})}[\PRTTT_k]_R,\phi)$ is a deformation of $n\Mod_{\QC(\mathcal{X})}[\PRTTT_k]$ to $R$,
then there exists an obstruction class $\textup{Obs}(n\Mod_{\QC(\mathcal{X})}[\PRTTT_k]_R,p)$
in $H^{i+n+2}(\mathcal{H}om_{\QC(\Map(S^n,\mathcal{X}))}(\OO_{\Delta}, \OO_{\Delta}))\otimes_kV$ for a deformation to $R'$.

\end{enumerate}

\end{Corollary}

We conclude the introduction with some directions for future research.
First, our framework applies to many $t$-structured presentable $(\infty,n)$-categories beyond the examples of this paper (i.e., sheaves of categories).
For example, 
one can consider the category of  (a suitable class of) $k$-linear $\infty$-categories with $t$-structures
whose $\Hom(\CCC_1,\CCC_2)$ from $\CCC_1$ to $\CCC_2$ is the $\infty$-category of 
right $t$-exact $k$-linear functors, which can be defined as $k$-linear $\infty$-categories with a $t$-structure.
After an appropriate procedure 
(right completions and left completions of $t$-structures on Hom), one obtain an $\PRTTT_k$-enriched category, which gives rise to
a $t$-structured $k$-linear presentable $(\infty,2)$-category.
One level higher, we can also consider the Morita category whose objects are (a suitable class of) 
$k$-linear monoidal $\infty$-categories with $t$-structures, and whose mapping object from $\mathcal{C}_1^\otimes$ to $\mathcal{C}_2^\otimes$ is the $\infty$-category of $\mathcal{C}_1^\otimes$-$\mathcal{C}_2^\otimes$-bimodule objects in $\PRTTT_k$. Such a Morita category should give rise to a complete $t$-structured $k$-linear presentable $(\infty,3)$-category, provided that the construction of Morita categories of algebras is carried out in the enriched ($k$-linear) setting. We expect that such a construction will become available in the near future, and that it will provide a large supply of examples to which our deformation theory applies.
Second, in Corollary~\ref{intromain4}, the Hochschild cohomology $\HH^\bullet(n\Mod_{\QC(\XX)}[\PRTTT_k]/k)$ is naturally an $\mathbf{E}_{n+2}$-algebra,
and is equivalent, as a $k$-module spectrum, to $H:=\Hom_{\QC(\Map(S^n,\XX))}(\OO_{\Delta},\OO_{\Delta})$.
On the other hand, To\"en \cite{Toen} constructed an $\mathbf{E}_{n+1}$-monoidal structure on $\QC(\Map(S^n,\XX))$
by a cobordism-theoretic construction using cospans, which endows $H$ with another $\mathbf{E}_{n+2}$-algebra structure.
We expect that these two structures coincide (see Remark~\ref{cobor}).
This would relate the deformation theory of categorified sheaves to topological field theory
and to deformation quantization of shifted Poisson structures.
Third, it would also be interesting to compare the deformation theory of an $\mathbf{E}_{n+1}$-algebra $A$ (see Francis \cite{FR})
with that of the $t$-structured $(\infty,n+1)$-category $n\PRTTT_A$: a deformation of $A$ induces that of $n\PRTTT_A$.

\vspace{2mm}

This paper is organized as follows.
In Section~\ref{Pre} we fix some notation and terminology
and review notions and some results.
In Section~\ref{TPHCsec} we introduce $t$-structured stable presentable $\infty$-categories
after reviewing generalities on $t$-structures and presentable $(\infty,n)$-categories.
In Section~\ref{HCKDsec} and Section~\ref{n1sec},
we prove the categorified Koszul duality equivalence, that is, Theorem~\ref{intromain1}.
In Section~\ref{hochsec}, we study Hochschild cohomology of $t$-structured presentable $(\infty,n)$-categories. 
Section~\ref{deformsec} discusses the deformation theory. We prove Theorem~\ref{intromain2} and Theorem~\ref{intromain3}.

\vspace{3mm}

Acknowledgements. The author is supported by JSPS KAKENHI grant.

\section{Preliminaries}

\label{Pre}

Let $k$ denote a fixed field, which serves as the base field.

\subsection{Notation}

\label{notation}

{\it Notation and Terminology.}
We use 
the language of $(\infty,1)$-categories.
We use the theory of {\it quasi-categories} as a model of $(\infty,1)$-categories.
Our main references are \cite{HTT}
 and \cite{HA}.
Following \cite{HTT}, we shall refer to quasi-categories
as {\it $\infty$-categories}.
To an ordinary category, we can assign an $\infty$-category by taking
its nerve, and therefore
when we treat ordinary categories, we often omit the nerve $\NNNN(-)$
and directly regard them as $\infty$-categories.
We use presentable higher categories, i.e., presentable $(\infty,n)$-categories in the sense of Stefanich \cite{Ste}, which we briefly review in Section~\ref{PHCsec}.
This is formulated within the theory of $(\infty,1)$-categories, so the reader is not required to have a background of models of
$(\infty,n)$-categories for $n>1$.

Here is a list of some of the notation and conventions we frequently use:

\begin{itemize}

\item $\simeq$: categorical equivalence or homotopy equivalence.  The symbol $\simeq$ is used to indicate isomorphisms and equivalences in categories and
$\infty$-categories.

\item $\CCC^{op}$: the opposite $\infty$-category of an $\infty$-category.

\item $\CCC^{\simeq}$: the largest Kan subcomplex of an $\infty$-category of $\CCC$.

\item $\Map_{\mathcal{C}}(C,C')$: the mapping space from an object $C\in\mathcal{C}$ to $C'\in \mathcal{C}$ where $\mathcal{C}$ is an $\infty$-category.

\item $\Fun(A,B)$: the function complex for simplicial sets $A$ and $B$. If $A$ and $B$ are $\infty$-categories, we regard $\Fun(A,B)$ as the functor category.

\item $\Cat$: the $\infty$-category of small $\infty$-categories. The word ``small'' means that it belongs to a fixed universe $\mathbb{U}$.
We write $\wCat$ for the $\infty$-category of $\infty$-categories in an enlarged universe $\mathbb{V}$ ($\mathbb{U}\in \mathbb{V}$). 
We write $\vCat$ for the (nonfull) subcategory of $\wCat$ which consists of $\infty$-categories having small colimits,
and whose $1$-morphisms are spanned by those functors which preserve small colimits.

\item $\SSS$: the $\infty$-category of $\mathbb{U}$-small spaces. We write $\wSSS$ for the $\mathbb{V}$-small spaces.

\item $\Ind(\CCC)$: the Ind-category of $\CCC$.

\end{itemize}

We will use the theory of higher algebras.
We refer the reader to \cite{HA}.
Here is a list of (some) of the notation about $\infty$-operads and algebras over them that we will use:

\begin{itemize}
\label{notation2}


\item $\Alg_n(\MMM)$:
For a symmetric monoidal $\infty$-category $\MMM^\otimes$,
we write $\Alg_{n}(\MMM)$ for
the $\infty$-category of algebra objects over the $\infty$-operad of little $n$-cubes $\mathbf{E}^\otimes_n$ in $\MMM^\otimes$  (cf. \cite[5.1]{HA}).
We refer to an object of $\Alg_{n}(\MMM)$ as an $\eenu$-algebra in $\MMM$.
If we denote by $\assoc^\otimes$ the associative operad (\cite[4.1.1]{HA}),
there is the standard equivalence $\assoc^\otimes\simeq \eone^\otimes$
of $\infty$-operads. We usually identify $\Alg_{1}(\MMM)$
with the $\infty$-category $\Alg_{\assoc}(\MMM)$, that is, the 
$\infty$-category of associative algebras in $\MMM$.

\item $\Alg_n^+(\MMM)=\Alg_n(\MMM)_{/\mathsf{1}}$: the $\infty$-category of augmented algebras. Here $\mathsf{1}$ is a unit algebra.

\item $\CAlg(\mathcal{M})$: the $\infty$-category of commutative
algebra objects in a symmetric
monoidal $\infty$-category $\mathcal{M}^\otimes$.
Instead of $\Alg_{\einf}(\MMM)$ or $\Alg_{\textup{comm}}(\MMM^\otimes)$
we write $\CAlg(\mathcal{M})$ for the $\infty$-category of algebra objects over the $\einf$-operad or the commutative operad in $\MMM^\otimes$.

\item $\Mod_R(\mathcal{M})$: 
the $\infty$-category of
$R$-module objects,
where $\mathcal{M}^\otimes$
is a symmetric monoidal $\infty$-category and $R\in \CAlg(\MMM)$.

\item $\Mod_R$: When $\MMM^\otimes$ is the symmetric monoidal $\infty$-category $\SP^\otimes$
of spectra, 
we write $\Mod_R^\otimes$
for
$\Mod_R^\otimes(\SP)$.
We denote by $\Mod_R$ the underlying category, that is, $\Mod_R(\SP)$.
When $R$ is the Eilenberg-Maclane ring spectrum $HA$, we write $\Mod_A$ for $\Mod_{HA}$.
For example, $\Mod_k$ is equivalent to the derived $\infty$-category of $k$-modules.

\item $\CAlg_k$: $\CAlg_k=\CAlg(\Mod_k)$.

\item 
$\LMod(\MMM)$:
the $\infty$-category $\Alg_{\LM}(\MMM)$ of $\LM^\otimes$-algebras
over $\LM^\otimes$: the $\infty$-operad defined in \cite[4.2.1.7]{HA}.
We may consider an object of  $\LMod(\MMM)$ to be a pair 
of $(A,M)$ such that $A$ is an associative algebra object in $\MMM$
and a left $A$-module $M$.
There is a Cartesian fibration $\LMod(\MMM)\to \Alg_{\assoc}(\MMM)\simeq \Alg_1(\MMM)$
which sends $(A,M)$ to $A$.
Another canonical morphism is the forgetful functor $\LMod(\MMM)\to \MMM$
which is informally given by $(A,M)\mapsto M$.

\item $\LMod_A(\MMM)$, $\LMod_A$: the $\infty$-category of left $A$-module objects in $\MMM$.
For $A\in \Alg_1(\MMM)$, we define $\LMod_A(\MMM)$
to be the fiber of $\LMod(\MMM)\to \Alg_1(\MMM)$ over $A$ in $\Cat$.
If $\MMM$ is the $\infty$-category of spectra (or $\Mod_k$),
we write $\LMod_A$ for $\LMod_A(\MMM)$.

\item $\RMod(\MMM)$: the right module version of $\LMod(\MMM)$

\item $\RMod_B(\MMM)$, $\RMod_A$: the $\infty$-category of right $B$-module objects in $\MMM$.
$\RMod_A$ is the right module version of $\LMod_A$.

\item $\PR$: $\infty$-category of presentable $\infty$-categories such that the mapping space consists of functors which preserve small colimits. 

\item $\PRST$, $\PR_k$: $\PRST$ is the $\infty$-category of presentable stable $\infty$-categories. A morphism $\CCC\to \DDD$ in $\PRST$ is a functor that preserves small colimits. It has a closed symmetric monoidal structure whose unit is the stable $\infty$-category $\SP$ of spectra  (see \cite[Definition 5.5.3.1]{HTT}, \cite[4.8]{HA}).
We write $\PR_k$ for $\Mod_{\Mod_k}(\PRST)$ and write $\otimes_k$ for the tensor products.

\end{itemize}

\subsection{Hom objects}
\label{homendobjectsec}
Let $\mathcal{M}$ be an $\infty$-category.
Let $\mathcal{A}^\otimes$ be a monoidal $\infty$-category.
Suppose that $\MMM$ is a left $\mathcal{A}^\otimes$-module, that is,
it has a left $\mathcal{A}^\otimes$-module action in $\wCat$.
To be precise, 
we let $\wCat^\times\times_{\Gamma}\mathcal{LM}^\otimes\to \mathcal{LM}^\otimes$
be the base change of the Cartesian symmetric monoidal $\infty$-category
$\wCat^\times\to \Gamma$, and we consider 
the $\infty$-category $\LMod(\wCat)$ of $\mathcal{LM}^\otimes$-algebra objects
(see \cite[4.2.1]{HA}).
This $\infty$-category $\LMod(\wCat)$ classifies pairs
$(\mathcal{A}^\otimes,\MMM)$ and admits a functor 
$\LMod(\wCat)\to \Alg_1(\wCat)$ which carries $(\mathcal{A}^\otimes,\MMM)$ to $\mathcal{A}^\otimes$.
We define $\LMod_{\AAA^\otimes}(\wCat)$ to be $\LMod(\wCat)\times_{\Alg_1(\wCat)}\{\AAA^\otimes\}$
and refer to it as the $\infty$-category of left $\mathcal{A}^\otimes$-modules.
In this context,
the Hom object from $M$ to $M'$ (defined) in $\mathcal{A}$ is an object $\mathcal{H}om(M,M')$ of $\mathcal{A}$
together with a morphism $f:\mathcal{H}om(M,M')\otimes M\to M'$, 
which has the universal property described as follows:
For any $E\in \mathcal{A}$, the composition with $f$ induces an equivalence of mapping spaces
\[
\Map_{\mathcal{A}}(E,\mathcal{H}om(M,M'))\to \Map_{\mathcal{M}}(E\otimes M, M')
\]
(in \cite{HA}, Hom objects are called morphism objects).

The endomorphism algebra object of $M\in \mathcal{M}$ (defined) in $\mathcal{A}$ is an associative
algebra $\End(M)$ in $\mathcal{A}$ (i.e., an object of $\Alg_1(\mathcal{A})$)
together with a left module action of $\End(M)$ on $M$, 
which has the universal property described informally as follows:
For any associative algebra $R$ in $\mathcal{A}$ and a left module action $R\otimes M\to M$
on $M$, there exists an essentially unique morphism $f:R\to \End(M)$
such that the left module action of $R$ is induced by the restriction of the action of $\End(M)$ along $f$.
To be precise, $\End(M)$ is defined to be an object of $\Alg_1(\mathcal{A}^\otimes)$ which represents
the functor $\Alg_1(\mathcal{A}^\otimes)^{op}\to \SSS$ corresponding to
the right fibration 
\[
\LMod(\MMM)\times_{\MMM}\{M\}\to \Alg_1(\mathcal{A}^\otimes).
\]
Here $\LMod(\MMM)\to \MMM$ is the forgetful functor.

The more specialized case is that $\MMM^\otimes$
is a monoidal $\infty$-category and the underlying $\infty$-category $\MMM$ is regarded as a left $\MMM^\otimes$-module via the tensor product functor.
In this case, Hom objects are referred to as internal Hom objects.

\subsection{Koszul duality of algebras}
\label{koszulsec}
Let $\MMM^\otimes$ be a symmetric monoidal $\infty$-category which has small colimits and limits such that the tensor product preserves sifted colimits. We write $\uni$ for a unit object.
We review the Koszul duality of algebras in $\MMM^\otimes$ in the formalism of pairings (see \cite{HA}, \cite[X]{DAG} for more details).
Let $\Alg^+_{1}(\MMM)=\Alg_1(\MMM)_{/\uni}$ be the $\infty$-category of augmented $\eone$-algebra objects (associative algebra objects) in $\MMM^\otimes$.
For a pair $(\epsilon_A:A\to \uni,\epsilon_B:B \to \uni)\in \Alg^+_{1}(\MMM)\times \Alg^+_{1}(\MMM)$ we put
\[
\textup{Pair}(\epsilon_A,\epsilon_B):= \Map_{\Alg_{1}(\MMM)}(A\otimes B,\uni)\times_{\Map(A,\uni)\times\Map(B,\uni)}\{(\epsilon_A,\epsilon_B)\}.
\]
Fix $\epsilon_B:B \to \uni$ and consider the functor
$u_{\epsilon_B}:\Alg_1^+(\MMM)^{op}\to \SSS$ given informally by $[\epsilon_A:A\to \uni]\mapsto \textup{Pair}(\epsilon_A,\epsilon_B)$. 
The Koszul dual $\DD_{\MMM^\otimes}(B)\in \Alg_1^+(\MMM)$ of $B$ (we abuse notation since $\DD_{\MMM^\otimes}(B)$ depends on the augmentation $\epsilon_B:B \to \uni$)
is defined by the following universal property that $\DD_{\MMM^\otimes}(B)$ represents the functor $u_{\epsilon_B}$.
Under a mild condition on $\MMM^\otimes$, 
there exist $\DD_{\MMM^\otimes}(B)$ and functorial equivalences
\[
\Map_{\Alg_1^+(\MMM)}(A,\DD_{\MMM^\otimes}(B))\simeq \textup{Pair}(\epsilon_A,\epsilon_B)
\]
for $\epsilon_A:A\to \uni$. We shall refer to the pairing $B\otimes\DD_{\MMM^\otimes}(B)\to \uni$ 
corresponding to the identity of $\DD_{\MMM^\otimes}(B)$ as the universal pairing for $B\to k$.
In that case, there is a functor 
\[
\DD_{\MMM^\otimes}:\Alg_1(\MMM) \longrightarrow \Alg_1(\MMM)^{op}
\]
which carries $A\to \uni$ to $\DD_{\MMM^\otimes}(A)\to \uni$.

The Koszul duals can be described in terms of bar constructions and endomorphism algebras.
Let $\textup{Bar}:\Alg_{1}^+(\MMM)\to \MMM$ be the functor
given by $[B\to \uni] \mapsto \uni\otimes_B\uni$, which we refer to as the bar construction of augmented algebras.
Then $\uni\otimes_B\uni$ admits a structure of a coaugmented coalgebra in a suitable way, and $\textup{Bar}$
is promoted to $\textup{Bar}:\Alg_{1}^+(\MMM)\to \Alg^+_{1}((\MMM)^{op})^{op}$
(notice that the target is the $\infty$-category of coaugmented coalgebras).
According to \cite[X, Example 4.4.20]{DAG}, for $B\in \Alg^+_{1}(\MMM)$, the Koszul dual $\DD_{\MMM^\otimes}(B)$
is equivalent to
the dual $\textup{Bar}(B)^\vee$ in $\MMM$, that is, the endomorphism algebra $\End_{\LMod_B(\MMM)}(\uni)\in \Alg_1(\MMM)$.

We recall the Koszul duality between $\eenu$-algebras. Set $\MMM^\otimes =\Mod_k$: in this case
Koszul duals exist (cf. \cite[X, 3.1.5, 4.3.8, 4.4.8]{DAG}).
Let $n\ge1$ be a natural number and  let $\Alg_{n}(\Mod_k)$ be the $\infty$-category of $\eenu$-algebras in $\Mod_k$.
We write $\Alg^+_{n}(\Mod_k)=\Alg_n(\Mod_k)_{/k}$ for the $\infty$-category
$\Alg_{n}(\Mod_k)_{/k}$ of augmented $\eenu$-algebras.
We consider the Koszul duals of augmented $\eenu$-algebras.
We take $\MMM^\otimes$ to be $\Alg_{n-1}^\otimes(\Mod_k)$. There exists the Koszul duality functor
\[
\DD_n:\Alg_{n}^+(\Mod_k)\simeq \Alg_1^+(\Alg_{n-1}(\Mod_k))\longrightarrow\Alg_1^+(\Alg_{n-1}(\Mod_k))^{op}\simeq  \Alg_{n}^+(\Mod_k)^{op}
\]
(see e.g. \cite[X, 4.4.8]{DAG}).
Let $\epsilon_{B}:B\to k$ be an augmented $\eenu$-algebra.
As mentioned above,
$\DD_n(B)$ is determined by the following universal property.
For an augmented $\eenu$-algebra $\epsilon_C:C\to k$, we write
\[
\textup{Pair}(\epsilon_B,\epsilon_C)= \Map_{\Alg_{n}(\Mod_k)}(B\otimes_k C,k)\times_{\Map(B,k)\times\Map(C,k)}\{(\epsilon_B,\epsilon_C)\}.
\]
Then the functor $\Alg^+_{n}(\Mod_k)^{op}\to \SSS$ given by $[\epsilon_C:C\to k]\mapsto \textup{Pair}(\epsilon_B,\epsilon_C)$
is representable by $\DD_{n}(B)$.
There is a universal (Koszul dual) pairing $B\otimes_k\DD_{n}(B)\to k$ corresponding to
the identity map $\textup{id}\in \Map_{\Alg^+_{n}(\Mod_k)}(\DD_{n}(B),\DD_{n}(B))$.
We refer to $\DD_{n}(B)$ as the $\eenu$-Koszul dual of $B$.

\section{$t$-structured presentable higher categories}

\label{TPHCsec}

\subsection{Generalities on $t$-structures}

\label{tsec}

Unless otherwise stated, we adopt the cohomological index.
Let us recall the notion of $t$-structures on stable $\infty$-categories.
Let $\CCC$ be a stable $\infty$-category.
Recall that the homotopy category $\hhh(\CCC)$ has the structure of a triangulated category.
A $t$-structure on $\CCC$ is a pair of full subcategories $(\CCC^{\le0},\CCC^{\ge0})$ of $\CCC$
such that $(\hhh(\CCC^{\le0}), \hhh(\CCC^{\ge0}))$ is a $t$-structure on the triangulated category $\hhh(\CCC)$ (see e.g. \cite{HA}, \cite[Section 2.1]{IH}).
We shall write $\CCC^{\le n}$ and $\CCC^{\ge n}$ for $\CCC^{\le 0}[-n]$ and $\CCC^{\ge 0}[-n]$, respectively.

Suppose that $(\CCC^{\le0},\CCC^{\ge0})$ is a $t$-structure.
The inclusion $\CCC^{\le n}\hookrightarrow\CCC$ has a right adjoint 
$\tau^{\le n}:\CCC\longrightarrow \CCC^{\le n}$
and the inclusion $\CCC^{\ge n}\hookrightarrow\CCC$
has a left adjoint $\tau^{\ge n}:\CCC\longrightarrow \CCC^{\ge n}$.
There exists a canonical equivalence $\tau^{\ge 0}\circ \tau^{\le0}\simeq \tau^{\le 0}\circ \tau^{\ge0}$.
We often use the notation $(\CCC^{\le0}, \CCC^{\ge0})$ to denote a $t$-structure on $\CCC$.
Note that $\CCC^{\ge0}$ is uniquely determined by $\CCC^{\le0}$: the full subcategory $\CCC^{\ge0}$ is spanned by those objects $C'$ such that $\Map_{\CCC}(C,C')$ is contractible for any $C\in \CCC^{\le -1}$. Dually, $\CCC^{\le0}$ is uniquely determined by $\CCC^{\ge0}$.
By $(\CCC,\CCC^{\le0})$ we mean a stable $\infty$-category endowed with the $t$-structure determined by a full subcategory $\CCC^{\le0}$.

\begin{Definition}
Let $(\CCC^{\le 0}, \CCC^{\ge 0})$ be a $t$-structure on a stable $\infty$-category $\CCC$.
\begin{enumerate}
\item For any $n \in \ZZ$, the left adjoint functor $\tau^{\ge n}: \CCC \to \CCC^{\ge n}$ to the inclusion
$\CCC^{\ge n} \subset \CCC$ induces a morphism $C \to \tau^{\ge n}(C)$, determined by the unit map. For all $n \in \ZZ$, this yields the sequence
\[
C \to \cdots \to \tau^{\ge n-1}(C) \to \tau^{\ge n}(C) \to \tau^{\ge n+1}(C) \to \cdots
\]
in $\CCC$.
The $t$-structure is left complete if, for any $C \in \CCC$, the induced morphism $C \to \varprojlim_{n \in \ZZ} \tau^{\ge n}(C)$ is an equivalence.
\item The right adjoint functor $\tau^{\le n}: \CCC \to \CCC^{\le n}$ to the inclusion
$\CCC^{\le n} \subset \CCC$ induces a morphism $\tau^{\le n}(C) \to C$, determined by the counit map. For all $n \in \ZZ$, this yields the sequence
\[
\cdots \to \tau^{\le n-1}(C) \to \tau^{\le n}(C) \to \tau^{\le n+1}(C) \to \cdots \to C
\]
in $\CCC$.
The $t$-structure is right complete if, for any $C \in \CCC$, the induced morphism $\varinjlim_{n \in \ZZ} \tau^{\le n}(C) \to C$ is an equivalence.

\item Suppose that $\CCC$ is presentable. The $t$-structure $(\CCC^{\le 0}, \CCC^{\ge 0})$ is accessible if the full subcategory $\CCC^{\le 0}$ is presentable.
\item Suppose that $\CCC$ has a monoidal structure. If the tensor product functor $\CCC \times \CCC \to \CCC$ sends $\CCC^{\le 0} \times \CCC^{\le 0}$ to $\CCC^{\le 0}$ and the unit object lies in $\CCC^{\le 0}$, we say that the $t$-structure is compatible with the monoidal structure.
\end{enumerate}
\end{Definition}

\begin{Example}
\label{connectivetstr}
Let $A$ be an associative connective ring spectrum, and let $\LMod_A$ be the stable $\infty$-category of left $A$-module spectra.
Define $\LMod_A^{\le 0}$ as the full subcategory of $\LMod_A$ spanned by objects $M$ such that $\pi_n(M) = 0$ for all $n < 0$, and $\LMod_A^{\ge 0}$ as the full subcategory spanned by objects $M$ such that $\pi_n(M) = 0$ for all $n > 0$.
Then $(\LMod_A^{\le 0}, \LMod_A^{\ge 0})$ defines an accessible, left complete, and right complete $t$-structure (see \cite[7.1.1.13]{HA}).
The full subcategory $\LMod_A^{\le 0}$ is also characterized as the smallest full subcategory which contains $A$ and is closed under small colimits.
Moreover, the $t$-structure is compatible with filtered colimits, as $\LMod_A^{\ge 0} \subset \LMod_A$ is closed under filtered colimits.
This $t$-structure is compatible with the monoidal structure when $A$ is an $\etwo$-algebra.
The heart $\LMod_A^{\heartsuit}$ is equivalent to the abelian category of left $\pi_0(A)$-modules.
\end{Example}

\begin{Remark}
\label{unifiedtstr}
Example~\ref{connectivetstr} is generalized to the nonconnective case as follows.
For $A\in \Alg_1(\SP)$ we define $\LMod_A^{\le0}$ to be the smallest full subcategory
which contains $A$ and is closed under colimits and extensions. By \cite[1.4.4.11]{HA},
$\LMod_A^{\le0}$ determines an accessible $t$-structure on $\LMod_A$ such that $\LMod_A^{\ge0}$
consists of those objects $M$ such that $\pi_0(\Map_{\LMod_A}(A[n],M))=\pi_n(M)=H^{-n}(M)=0$ for $n>0$.
The full subcategory $\LMod_A^{\ge0}$ is closed under filtered colimits, and the $t$-structure is right complete.
We shall refer to this $t$-structure as the standard $t$-structure.
When $A$ is an $\etwo$-algebra, the $t$-structure is compatible with the monoidal structure on $\LMod_A$.
If $A$ is not connective, the resulting $t$-structure is generally not left complete.
\end{Remark}

\begin{Example}
\label{Indtstr}
Let $\CCC$ be a small stable $\infty$-category.
Let $(\CCC^{\le0}, \CCC^{\ge0})$ be a $t$-structure on $\CCC$.
Then $(\Ind(\CCC^{\le0}),\Ind(\CCC^{\ge0}))$ defines a $t$-structure on $\Ind(\CCC)$.
By the construction, the $t$-structure is compatible with filtered colimits.
If the $t$-structure on $\CCC$ is right bounded, then 
$(\Ind(\CCC^{\le0}),\Ind(\CCC^{\ge0}))$ is right complete (see \cite[VIII, Lemma 5.4.1]{DAG}).

For instance, suppose that $R$ is a noetherian connective $\eone$-algebra over a field $k$, 
that is, an object of $\Alg_1(\Mod_k^{\le0})$.
We let $\textup{LCoh}(R)$ denote the full subcategory of $\LMod_R$ spanned by those objects $M$
such that $H^n(M)=0$ for $|n|>>0$, and each $H^n(M)$ is a finitely generated $H^0(R)$-module.
This small stable $\infty$-category $\textup{LCoh}(R)$ inherits a bounded $t$-structure from that of $\LMod_R$ in Example~\ref{connectivetstr}.
Then $(\Ind(\textup{LCoh}(R)^{\le0}),\Ind(\textup{LCoh}(R)^{\ge0}))$ 
is a right complete accessible $t$-structure. This is generally not left complete.
\end{Example}

Let $\PRT$ be the $\infty$-category of presentable stable $\infty$-categories equipped with accessible $t$-structures.
An object of $\PRT$ is a pair $(\CCC,\CCC^{\le0})$ where $\CCC$ is a presentable stable $\infty$-category, and $\CCC^{\le0}$ is a  full subcategory
inducing a $t$-structure,
such that $\CCC^{\le0}$ itself is presentable.
A morphism $(\CCC, \CCC^{\le0})\to (\DDD,\DDD^{\le0})$ is a colimit-preserving functor $F:\CCC\to  \DDD$ which is right $t$-exact (i.e., $F(\CCC^{\le0})\subset \DDD^{\le0}$).
Let $\PRTT$ be the full subcategory of $\PRT$, spanned by presentable stable $\infty$-categories equipped with accessible right complete $t$-structures.
Let $\PRTTT$ be the full subcategory of $\PRT$, spanned by presentable stable $\infty$-categories equipped with accessible left complete and right complete $t$-structures.
There exists an adjoint pair
\[
\xymatrix{
I^+: \PRTT \ar@<0.5ex>[r] &  \PRT:\widehat{R}  \ar@<0.5ex>[l]  
}
\]
such that $I^+$ is the inclusion, and the right adjoint $\widehat{R}$ carries $(\CCC,\CCC^{\le0})$
to the right completion $\widehat{R}(\CCC,\CCC^{\le0})$, which is
$(\colim_{n\to \infty} \CCC^{\le n},\CCC^{\le0})$
where $\colim_{n\to \infty} \CCC^{\le n}$ is the colimit of the sequence of inclusions
\[
\CCC^{\le 0}\to \CCC^{\le 1}\to \CCC^{\le 2}\to \cdots\to \CCC^{\le n}\to \cdots
\]
in $\PRST$
(cf. \cite[5.5.3.4, 5.5.3.18]{HTT}). 

The left completion $\HL((\CCC,\CCC^{\le0}))$ of $(\CCC,\CCC^{\le0})$
is defined as
$\lim_{n\to -\infty}\CCC^{\ge n}$ in $\PRST$
of the sequence 
\[
\cdots \to \CCC^{\ge n-1}\to \CCC^{\ge n}\to \CCC^{\ge n+1} \to \cdots
\]
induced by truncations.
We usually write $\HL(\CCC)$ for $\HL((\CCC,\CCC^{\le0}))$.
If $(\CCC,\CCC^{\le0})\in \PRTT$, then the left completion $\HL(\CCC)$ is also accessible and right complete (see \cite[VIII, 4.6.13]{DAG}).
There exists a localization adjoint pair
\[
\xymatrix{
\HL :\PRTT \ar@<0.5ex>[r] &  \PRTTT:I^{\pm}  \ar@<0.5ex>[l]  
}
\]
such that $I^\pm$ is the inclusion, and $\HL$ is the left adjoint given by the left completion.

{\it Symmetric monoidal structures.}
Next, we recall the symmetric monoidal structure on $\PRT$ (see \cite[VIII, 4.6.1]{DAG}).
Let $(\mathcal{B},\mathcal{B}^{\le0})$ and $(\CCC,\CCC^{\le0})$ be presentable $\infty$-categories endowed with accessible $t$-structures. 
The tensor product of $(\mathcal{B},\mathcal{B}^{\le0})$ and $(\CCC,\CCC^{\le0})$
is defined as
$\mathcal{B}\otimes \CCC$ (in $\PRST$) 
together with the smallest full subcategory $(\mathcal{B}\otimes \CCC)^{\le0}$
which contains the essential image of $\mathcal{B}^{\le0}\times \CCC^{\le0}\to \mathcal{B}\otimes \CCC$
and is closed under small colimits and extensions.
By \cite[1.4.5.11]{HA}, $(\mathcal{B}\otimes \CCC)^{\le0}$ defines an accessible $t$-structure.
The universal property of $(\mathcal{B}\otimes \CCC, (\mathcal{B}\otimes \CCC)^{\le0})$ is described as follows: For any $(\DDD, \DDD^{\le0})\in \PRT$
the composition with the canonical functor $\mathcal{B}\times \CCC\to \mathcal{B}\otimes \CCC$ gives an equivalence
\[
\Fun_{\PRT}((\mathcal{B}\otimes \CCC, (\mathcal{B}\otimes \CCC)^{\le0}),(\DDD,\DDD^{\le0}))\to \Fun^{\le}(\mathcal{B}\times \CCC,\DDD)
\]
where $\Fun^{\le}(\mathcal{B}\times \CCC,\DDD)$ is the full subcategory of $\Fun(\mathcal{B}\times \CCC,\DDD)$,
which consists of those functors $F$ that preserve small colimits separately in each variable, such that the essential image of $\mathcal{B}^{\le0}\times \CCC^{\le0}$ 
is contained in $\DDD^{\le0}$.

If $(\mathcal{B},\mathcal{B}^{\le0})$ and $(\CCC,\CCC^{\le0})$ are right complete,
then the tensor product $(\mathcal{B}\otimes \CCC, (\mathcal{B}\otimes \CCC)^{\le0})$ in $\PRT$
is also right complete (see \cite[VIII, 4.6.11]{DAG}).
Thus, $\PRTT$ inherits a symmetric monoidal structure from that of $\PRT$.
The left adjoint $I^+$ is promoted
to a symmetric monoidal functor.
The $\infty$-category $\PRTTT$ also has a symmetric monoidal structure.
Let $(\mathcal{B},\mathcal{B}^{\le0})$ and $(\CCC,\CCC^{\le0})$ be objects of $\PRTTT$.
Explicitly, the tensor product of $(\mathcal{B},\mathcal{B}^{\le0})$ and $(\CCC,\CCC^{\le0})$ in $\PRTTT$
is given by 
\[
\HL((\mathcal{B}\otimes\CCC,  (\mathcal{B}\otimes \CCC)^{\le0})).
\]
We usually write $\HL(\mathcal{B}\otimes\CCC)$ or $\mathcal{B}\widehat{\otimes}\CCC$ for it.
The functor $\HL:\PRTT\to \PRTTT$ is promoted to a symmetric monoidal functor. 
Since $\HL$ is a left adjoint functor, $\HL$ commutes with the relative tensor products.
The inclusion $I^{\pm}$
is not a symmetric monoidal functor but a lax symmetric monoidal functor.

\begin{Remark}
\label{inthomtstr}
The symmetric monoidal $\PRT$ admits internal Hom objects: For $\mathcal{P}, \mathcal{Q}\in \PRT$, an internal Hom object $\Hom_{\PRT}(\mathcal{P},\mathcal{Q})$ from $\mathcal{P}$ to $\mathcal{Q}$ in $\PRT$ is 
given  by the pair $(\Hom_{\PR}(\underline{\mathcal{P}},\underline{\mathcal{Q}}), \Hom_{\PR}(\underline{\mathcal{P}},\underline{\mathcal{Q}})^{\le0})$
such that $\Hom_{\PR}(\underline{\mathcal{P}},\underline{\mathcal{Q}})\in \PR$ 
is the internal Hom object from the underlying object $\underline{\mathcal{P}}$ to $\underline{\mathcal{Q}}$, and 
 $\Hom_{\PR}(\underline{\mathcal{P}},\underline{\mathcal{Q}})^{\le0}$ is the full subcategory, spanned by right $t$-exact functors.
Thus, $\PRTT$ also admits internal Hom objects obtained by taking right completions of internal Hom objects in $\PRT$.
\end{Remark}

\begin{Notation}
Let $\PRTT_k$ denote $\Mod_{\Mod_k}(\PRTT)$ where we regard $\Mod_k$ as an object of $\CAlg(\PRTT)$ by taking $\Mod_k^{\le0}$
to be the full subcategory spanned by connective $k$-module spectra. 
We write $\PRTTT_k$ for $\Mod_{\Mod_k}(\PRTTT)$.
\end{Notation}

\subsection{Presentable higher categories}

\label{PHCsec}

Let $\vCat$ denote the subcategory of the $\infty$-category of $\mathbb{V}$-small $\infty$-categories 
such that objects are spanned by those $\infty$-categories which admit $\UUUU$-small colimits,
and $1$-morphisms are spanned by functors which preserve $\UUUU$-small colimits.
According to \cite[4.8.4.2]{HA}, $\vCat$ is a $\VVVV$-presentable $\infty$-category.

For instance, 
an $\mathbb{U}$-presentable $\infty$-category belongs to $\vCat$.
The $\infty$-category $\PR$ of presentable $\infty$-categories has $\UUUU$-small colimits
so that $\PR$ belongs to $\vCat$. 
In addition,
by \cite[Proposition 5.1.4]{Ste}, it is $\UUUU$-compactly generated $\infty$-category, and
\[
\PR\subset \vCat
\]
is characterized as the full subcategory of $\kappa_0$-compact objects
where $\kappa_0$ is the cardinal of the set of $\UUUU$-small sets.
($\kappa_0$ is the cardinality of $\UUUU$)

The $\infty$-category $\vCat$ has a symmetric monoidal structure (see \cite[4.8.1.3]{HA}.
For $\DDD,\DDD'\in \vCat$, there is $\DDD\otimes \DDD' \in \vCat$ with a functor $u:\DDD\times \DDD'\to \DDD\otimes \DDD'$
which preserves $\UUUU$-small colimits separately in each variable and further satisfies
the universal property that for any $\mathcal{E}\in \vCat$, the composition with $u$ induces
an equivalence $\Map_{\vCat}(\DDD\otimes \DDD',\mathcal{E})\to \Map_{\vCat}^e(\DDD\times \DDD',\mathcal{E})$
of mapping spaces. Here $\Map_{\vCat}^e(\DDD\times \DDD',\mathcal{E})$ is the full subcategory of $\Map_{\vCat}(\DDD\times \DDD',\mathcal{E})$,
which consists of functors preserving small colimits separately in each variable.
Moreover, $\PR$ is closed under tensor products,
and the tensor product functor 
$\otimes:\PR\times \PR\to \PR$ preserves $\UUUU$-small colimits separately in each  variable.
Thus, $\PR$ is a commutative algebra object of $\vCat$.

We consider $\Mod_{\PR}(\vCat)$ and define $2\PR$ to be the full subcategory of $\Mod_{\PR}(\vCat)$,
which consists of $\kappa_0$-compact objects.
The $\infty$-category $2\PR$ itself is closed under $\UUUU$-small colimits 
and is generated by $\PR\otimes \Fun(\Delta^1,\SSS)$ under small colimits.
In particular, $2\PR$ is closed under taking tensor products in $\Mod_{\PR}(\vCat)$
so that $2\PR$ belongs to $\CAlg(\vCat)$.

More generally, this argument is extended to the following:

\begin{Construction}
\label{constructionnc}
For any $\mathcal{\CCC}\in \CAlg(\vCat)$
the full subcategory $\Mod_{\CCC}(\vCat)^{\kappa_0}$ of $\Mod_{\CCC}(\vCat)$,
which consists of $\kappa_0$-compact objects, belongs to $\CAlg(\vCat)$.
This assignment $\CCC\mapsto \Mod_{\CCC}(\vCat)^{\kappa_0}$ determines a functor
\[
\mathsf{M}:\CAlg(\vCat)  \longrightarrow \CAlg(\vCat)
\]
such that any morphism $f:\mathcal{C}\to \mathcal{C}'$ defines the base change functor
$(-)\otimes_{\mathcal{C}}\mathcal{C}':\Mod_{\mathcal{C}}(\vCat)\to \Mod_{\mathcal{C}'}(\vCat)$ whose restriction 
induces $\mathsf{M}(f):\mathsf{M}(\mathcal{C})=\Mod_{\mathcal{C}}(\vCat)^{\kappa_0}\to \mathsf{M}(\mathcal{C}')=\Mod_{\mathcal{C}'}(\vCat)^{\kappa_0}$.
For $n\ge2$ we define $n\CCC$ to be $\mathsf{M}^{n-1}(\CCC)$. By convention, $1\CCC$ is $\CCC$.
\end{Construction}

Using this construction
we define $n\PR$ to be the full subcategory $\Mod_{(n-1)\PR}(\vCat)^{\kappa_0}$ of $\Mod_{(n-1)\PR}(\vCat)$,
which consists of $\kappa_0$-compact objects. That is, $n\PR=\mathsf{M}^{n-1}(\PR)$.
By convention, if $n=1$, $1\PR$ is $\PR$.
We refer to an object of $n\PR$ as a presentable $(\infty,n)$-category.

Let $\PRST$ denote the $\infty$-category of stable presentable $\infty$-categories. 
We often identify it with $\Mod_{\SP}(\PR)$.
For $n\ge2$,
we define $n\PRST$ to be $\mathsf{M}^{n-1}(\PRST)$. 
We put $1\PRST=\PRST$.
We refer to an object of $n\PRST$ as a stable presentable $(\infty,n)$-category.
By convention, we set $0\PRST=\SP$

\subsection{$t$-structured presentable higher categories}

Let $\PRTT$ be the $\infty$-category of stable presentable $\infty$-categories equipped with right complete and accessible $t$-structures.
A $1$-morphism of $\PRTT$ is a functor that preserves small colimits and is right $t$-exact.

\subsubsection{Colimits and limits}
Let
\[
\SP(\CCC)=\CCC\otimes_{\SPC}\SP\simeq \colim_{\Sigma}(\CCC\stackrel{\Sigma}{\to} \CCC  \stackrel{\Sigma}{\to} \CCC \stackrel{\Sigma}{\to}\cdots ) 
\]
be the stabilization of $\CCC$, where $\Sigma C$ is the pushout $\ast\sqcup_{C}\ast$ for $C\in \CCC$, and the right-hand side is a colimit in $\PRST$ (equivalently, in $\PR$).
Consider 
its full subcategory $\SP(\CCC)^{\le0}$ that is the smallest one, which contains the image of $\CCC$ and is closed under extensions and small colimits.
According to \cite[1.4.3.4]{HA}, $\SP(\CCC)^{\le0}$ is presentable, and the pair $(\SP(\CCC),\SP(\CCC)^{\le0})$ constitutes 
an accessible $t$-structure.
As proved below, $\SP(\CCC)^{\le0}$ is a right complete $t$-structure on $\SP(\CCC)$.

\begin{Proposition}
\label{prttadjoint}
Let $\mathfrak{C}:\PRTT\to \Mod_{\SPC}(\PR)$ be the functor given by $(\DDD,\DDD^{\le0})\mapsto \DDD^{\le0}$.
Then there is a left adjoint $\mathfrak{S}:\Mod_{\SPC}(\PR)\to \PRTT$ given by $\CCC \mapsto (\SP(\CCC),\SP(\CCC)^{\le0})$,
which constitutes an adjunction
\[
\xymatrix{
\mathfrak{S}:\Mod_{\SPC}(\PR)  \ar@<0.5ex>[r] &   \PRTT :\mathfrak{C}.  \ar@<0.5ex>[l]   \\
}
\]
Moreover,  $\mathfrak{C}$ is fully faithful.
\end{Proposition}

\begin{Remark}
More precisely, the right adjoint is defined by
\[
\PRTT\simeq \Mod_{(\SP,\SPC)}(\PRTT) \to \Mod_{\SPC}(\PR)
\]
where the left equivalence comes from the fact that $(\SP,\SPC)$ is a unit object of $\PRTT$, and the right functor 
is given by $(\DDD,\DDD^{\le0})\mapsto \DDD^{\le0}$.
\end{Remark}

\Proof
Consider the inclusions $\CCC\to \SP(\CCC)^{\le0}\to \SP(\CCC)$.
For $(\DDD,\DDD^{\le0})\in \PRTT$, the composition induces
\[
\Map_{\PRT}((\SP(\CCC),\SP(\CCC)^{\le0}),(\DDD,\DDD^{\le0})) \stackrel{f}{\to} \Map_{\PR}(\SP(\CCC)^{\le0},\DDD^{\le0}) \stackrel{g}{\to} \Map_{\PR}(\CCC,\DDD^{\le0}).
\]
Using the universal property of $\SP(\CCC)$, that is an equivalence 
$\Map_{\PRST}(\SP(\CCC),\DDD)\simeq \Map_{\PR}(\CCC, \DDD)$ induced by $\CCC\to \SP(\CCC)$,
we see that the composite $g\circ f$ is an equivalence.
Thus, to see that $\mathfrak{S}$ is a left adjoint to $\mathfrak{C}$,
it suffices to prove that $(\SP(\CCC),\SP(\CCC)^{\le0})$ belongs to $\PRTT$.
(This also implies that both $f$ and $g$ are equivalences.) The counit map $\colim_{\Sigma}\DDD^{\le0}\to \DDD$, and the unit map is $\CCC\to \SP(\CCC)^{\le0}$)
We consider the sequence of stabilizations 
\[
\colim_{\Sigma}\CCC=\colim(\CCC\stackrel{\Sigma}{\to}\CCC \stackrel{\Sigma}{\to}\cdots )\to \colim_{\Sigma}\SP(\CCC)^{\le0}=\colim(\SP(\CCC)^{\le0} \stackrel{\Sigma}{\to} \SP(\CCC)^{\le0}\stackrel{\Sigma}{\to}\cdots ) \to \colim_{\Sigma}\SP(\CCC)
\]
induced by $\CCC\to \SP(\CCC)^{\le0}\to \SP(\CCC)$.
This can naturally be identified with $\SP(\CCC)\stackrel{\phi}{\to} \widehat{R}((\SP(\CCC))=\widehat{R}((\SP(\CCC),\SP(\CCC)^{\le0})) \stackrel{\psi}{\to} \SP(\CCC)$.
The composite $\psi\circ\phi$ is the canonical functor $\SP(\CCC)\to \SP(\SP(\CCC))$ so that it is an equivalence.
By construction, the composite $\phi\circ \psi$ induces an identity functor $\widehat{R}(\SP(\CCC))\supset \SP(\CCC)^{\le0}\stackrel{\simeq}{\to} \SP(\CCC)^{\le0}\subset \widehat{R}(\SP(\CCC))$,
so that $\phi\circ \psi$ is an equivalence since $\SP(\CCC)^{\le0}$ determines a right complete $t$-structure on 
$\widehat{R}(\SP(\CCC))$. Thus, $\phi$ and $\psi$ are mutually inverse equivalences $\SP(\CCC)\simeq \widehat{R}(\SP(\CCC))$.
Thus, $\SP(\CCC)^{\le0}$ determines a right complete $t$-structure.

For any $(\DDD,\DDD^{\le0})\in \PRTT$, the canonical map $\colim_{\Sigma}\DDD^{\le0}\to \DDD$ is an equivalence.
It follows that $\mathfrak{S}\circ \mathfrak{C}$ is equivalent to the identity functor. Thus, $\mathfrak{C}$ is fully faithful.
\QED

\begin{Corollary}
\label{tstablecolimits}
The $\infty$-category $\PRTT$ has small colimits. The forgetful functor $\PRTT\to \PRST$ preserves small colimits.
Moreover, $\PRTT$ belongs to $\CAlg(\vCat)$.
\end{Corollary}

\Proof
Since $\PR$ has small colimits, $\Mod_{\SPC}(\PR)$ also has small colimits (cf. \cite[4.2.3.5]{HA}).
If $f:I\to \PRTT$ is a diagram indexed by a small $\infty$-category $I$,
we can take a colimit $\CCC_{\infty}$ of the composite $I\to \PRTT\stackrel{\mathfrak{C}}{\hookrightarrow} \Mod_{\SPC}(\PR)$.
By Proposition~\ref{prttadjoint}, we see that  $\mathfrak{S}(\CCC_{\infty})$ is a colimit of $f$.
Note that if we write $(\CCC_i,\CCC_i^{\le0})=f(i)$, the colimit of $f$ is $(\SP(\colim_{i\in I}\CCC_i^{\le0}), \SP(\colim_{i\in I}\CCC_i^{\le0}))^{\le0})$
where $\colim_{i\in I}\CCC_i^{\le0}$ is a colimit taken in $\Mod_{\SPC}(\PR)$. We note the equivalences
\[
\SP(\colim_{i\in I}\CCC_i^{\le0})\simeq (\colim_{i\in I}\CCC_i^{\le0})\otimes_{\SPC}\SP\simeq \colim_{i\in I}(\CCC_i^{\le0}\otimes_{\SPC}\SP)\simeq \colim_{i\in I}\SP(\CCC_i^{\le0})\simeq \colim_{i\in I}\CCC_i.
\]
Thus, the underlying presentable stable $\infty$-category of the colimit of $f$ is $\colim_{i\in I}\CCC_i$,
and $\SP(\colim_{i\in I}\CCC_i^{\le0}))^{\le0}$ is the smallest full subcategory of $\colim_{i\in I}\CCC_i$ which contains the image of $\colim_{i\in I}\CCC_i^{\le0}$
and is closed under extensions and small colimits.
In particular, the forgetful functor $\PRTT\to \PRST$ preserves small colimits.

For the last assertion, it suffices to show that
the tensor product functor $\otimes:\PRTT\times \PRTT\to \PRTT$ preserves small colimits 
separately in each variable.
The tensor product $\PR$ preserves small colimits separately in each variable.
For $(\DDD,\DDD^{\le0})\in \PRTT$,
we have a commutative diagram in $\PR$:
\[
\xymatrix{
(\colim_{i\in I} \CCC_i^{\le0})\otimes \DDD^{\le0} \ar[r] \ar[d]^{\simeq} &  (\colim_{i\in I} \CCC_i)\otimes \DDD \ar[d]^{\simeq} \\ 
\colim_{i\in I} (\CCC_i^{\le0}\otimes \DDD^{\le0}) \ar[r]  &  \colim_{i\in I} (\CCC_i\otimes \DDD) 
}
\]
Consequently, the tensor product in $\PRTT$ preserves small colimits 
separately in each variable.
\QED

\begin{Corollary}
\label{tstablecolimits2}
The $\infty$-category $\PRTTT$ has small colimits.
Moreover, $\PRTTT$ belongs to $\CAlg(\vCat)$.
\end{Corollary}

\Proof
Use Corollary~\ref{tstablecolimits} and the left completion functor $\widehat{L}:\PRTT\to \PRTTT$.
(Note that $\widehat{L}$ is symmetric monoidal and preserves small colimits.)
\QED

\begin{Proposition}
\label{PRTcolimitlimit}
The following holds.
\begin{enumerate}

\item The $\infty$-category $\PRT$ has small limits. The forgetful functor $\PRT\to \PRST$ preserves small limits.

\item The $\infty$-category $\PRTT$ has small limits.

\end{enumerate}
\end{Proposition}

\Proof
A small limit in $\PR$ can be computed as a small limit in $\wCat$,
so that $\lim_{i\in I}(\CCC_i,\CCC_i^{\le0})$ of a diagram $D:I\to \PRT\subset \Fun(\Delta^1,\PR)$
is $\lim_{i\in I}\CCC_i^{\le0}\to \lim_{i\in I}\CCC_i$
where $\lim_{i\in I}\CCC_i^{\le0}$ and $\lim_{i\in I}\CCC_i$ are limits in $\wCat$ (equivalently, limits in $\PR$).
Thus the induced functor $C:\lim_{i\in I}\CCC_i^{\le0}\to \lim_{i\in I}\CCC_i$
is fully faithful. Since each $\CCC_i\to \CCC_j$ is a right $t$-exact colimit-preserving functor between stable $\infty$-categories,
it follows that the essential image of $C$ is closed under small colimits and extensions.
Thus, a limit of $D\to \Fun(\Delta^1,\PR)$ lies in $\PRT$. Thus, $\PRT$ has small limits,
and 
the forgetful functor $\PRT\to \PRST$ preserves small limits. This proves (1).
The right completion functor $\PRT\to \PRTT$ is right adjoint so that it preserves small limits. This proves (2). 
\QED

\subsubsection{}
We define $t$-structured presentable $(\infty,n)$-categories.

\begin{Definition}
For a positive integer $n\ge1$ and $\PRTT\in \CAlg(\vCat)$ (see Corollary~\ref{tstablecolimits}), we define $n\PRTT\in \CAlg(\vCat)$ (see Construction~\ref{constructionnc}).
By convention, $1\PRTT=\PRTT$ and $0\PRTT=(\SP,\SPC)$.
We refer to an object of $n\PRTT$ as a $t$-structured presentable $(\infty, n)$-category.

For a positive integer $n\ge1$ and $\PRTTT \in \CAlg(\vCat)$ (see Corollary~\ref{tstablecolimits2}), we define $n\PRTTT\in \CAlg(\vCat)$ (see Construction~\ref{constructionnc}).
By convention, $1\PRTTT=\PRTTT$ and $0\PRTTT=(\SP,\SPC)$.
We refer to an object of $n\PRTTT$ as a complete $t$-structured presentable $(\infty, n)$-category.
\end{Definition}

\begin{Remark}
In this paper, we will not consider $\PRT$-modules and their higher variants as the main objects of interest
because they do not appear in categorified Koszul duality.
Thus, we simply call objects of $n\PRTT$ $t$-structured presentable $(\infty, n)$-categories.
\end{Remark}

\subsubsection{Completions of $t$-structured presentable categories}

\begin{Definition}
We consider the symmetric monoidal localization functor $\widehat{L}:\PRTT\to \PRTTT$,
whose right adjoint is the lax symmetric monoidal inclusion $\PRTTT\subset \PRTT$.
The base change $\Mod_{\PRTT}(\vCat)\to \Mod_{\PRTTT}(\vCat)$
induces a symmetric monoidal functor $2\PRTT\to 2\PRTTT$.
Inductively, for any $n\ge2$,
a symmetric monoidal functor $(n-1)\PRTT\to (n-1)\PRTTT$ induces the base change functor
\[
(-)\otimes_{(n-1)\PRTT}(n-1)\PRTTT :n\PRTT\longrightarrow n\PRTTT.
\]
For $\DDD\in n\PRTT$, we define the completion of $\DDD$ as $\DDD\otimes_{(n-1)\PRTT}\PRTTT$.
For $n\ge2$, we write $n\widehat{L}=(-)\otimes_{(n-1)\PRTT}(n-1)\PRTTT$
and refer to it as the completion functor.
\end{Definition}

Let $S$ be the collection of all $1$-morphisms $f$ of $\PRTT$ such  that
$\widehat{L}(f)$ is an equivalence. Note that $S$ is multiplicative 
in the sense that for any $f\in S$ and any object $\CCC\in \PRTT$, $f\otimes \textup{id}_{\CCC}$ also belongs to $S$.
According to \cite[5.2.7.12]{HTT} and \cite[1.3.4.1, 4.1.3.4]{HA}, 
$\PRTTT$ is the symmetric monoidal $\infty$-category obtained from $\PRTT$ by inverting $S$:
For any symmetric monoidal $\infty$-category $\mathcal{D}^\otimes$, 
the functor $\Map(\PRTTT, \DDD^\otimes)\to \Map(\PRTT, \DDD^\otimes)$
induced by the composition with $\widehat{L}$ is
fully faithful (here $\Map$ means the mapping space between symmetric monoidal $\infty$-categories).
The essential image consists of those symmetric monoidal functors which carry $S$ to equivalences in $\DDD$.
The underlying $\infty$-category $\PRTTT$ is obtained from $\PRTT$ by inverting $S$.

\begin{Proposition}
Let $\iota:\Mod_{\PRTTT}(\vCat)\to \Mod_{\PRTT}(\vCat)$ be the restriction functor along $\widehat{L}:\PRTT\to \PRTTT$.
Then $\iota$ is fully faithful.
\end{Proposition}

\Proof
Since $\vCat$ has geometric realizations,
there exists a left adjoint $l:\Mod_{\PRTT}(\vCat)\to \Mod_{\PRTTT}(\vCat)$ of $\iota$ given by the base change along $\widehat{L}:\PRTT\to \PRTTT$.
To prove that $\iota$ is fully faithful, it suffices to show that for any $\mathcal{M}\in \Mod_{\PRTTT}(\vCat)$
the counit map $l(\iota(\mathcal{M}))\to \mathcal{M}$ is an equivalence.
Note that both $l$ and $\iota$ preserve (large) colimits. 
In addition, $\mathcal{M}$ can be described as a colimit of a free $\PRTTT$-modules.
Thus, we may and will assume that $\mathcal{E}=\CCC\otimes \PRTTT$ with some $\CCC\in \vCat$.
The counit map $l(\iota( \CCC\otimes \PRTTT)) \simeq l(\CCC\otimes \PRTT[S^{-1}])\simeq  \CCC\otimes (\PRTT[S^{-1}])[\overline{S}^{-1}]\to  \CCC\otimes \PRTT[S^{-1}]$
is an equivalence since $(\PRTT[S^{-1}])[\overline{S}^{-1}]\to  \PRTT[S^{-1}]$ is an equivalence where $[\overline{S}^{-1}]$ indicates the localization
defined by inverting the image $\overline{S}$ of $S$ in $\PRTT[S^{-1}]$.
\QED

Let $\PRTT_k$ denote $\Mod_{(\Mod_k,\Mod_k^{\le0})}(\PRTT)$.
Let $\PRTTT_k$ denote $\Mod_{(\Mod_k,\Mod_k^{\le0})}(\PRTTT)$.
The symmetric monoidal functor $\widehat{L}:\PRTT\to \PRTTT$ induces $\PRTT_k\to \PRTTT_k$ which we also denote by $\widehat{L}$.

The same argument shows: 

\begin{Proposition}
Let $\iota:\Mod_{\PRTTT_k}(\vCat)\to \Mod_{\PRTT_k}(\vCat)$ denote the restriction functor along $\widehat{L}:\PRTT_k\to \PRTTT_k$.
Then $\iota$ is fully faithful.
\end{Proposition}

In particular, we have the localization adjoint pair
\[
l:\Mod_{\PRTT_k}(\vCat) \rightleftarrows \Mod_{\PRTTT_k}(\vCat):\iota
\]
We also refer to $l$ as the completion functor.

\begin{Remark}
At the time of writing,
we do not know whether or not $\iota$ preserves $\kappa_0$-compact objects.
\end{Remark}

\subsubsection{Underlying stable $(\infty,n)$-categories}
Consider a stable presentable $\infty$-category $(\CCC,\CCC^{\le0})$ equipped with a $t$-structure.
One can consider the underlying stable presentable $\infty$-category $\CCC$
obtained by forgetting the $ t$-structure $\CCC^{\le0}$.
We define underlying stable $(\infty,n)$-categories
obtained from objects of $n\PRTT$.

Recall from Proposition~\ref{prttadjoint} that  
$\PRTT$ can be regarded as a full subcategory
of $\Mod_{\SPC}(\PR)$ through the assignment given by $(\CCC,\CCC^{\le0})\mapsto \CCC^{\le0}$.
We also note that the stabilization $\SP(\CCC^{\le0})$ of $\CCC^{\le0}$ is canonically equivalent to
$\CCC$. 
The tensor product with $\SP$ over $\SPC$ induces 
\[
\SP\otimes_{\SPC}(-):\Mod_{\SPC}(\PR)\longrightarrow \Mod_{\SP}(\PR)
\]
which carries $\PPP$ to the stabilization $\SP(\PPP)\simeq \SP\otimes_{\SPC}\PPP$.
This observation yields a commutative diagram
\[
\xymatrix{
\PRTT \ar[r] \ar[d] &  \PRST \ar[d] \\
\Mod_{\SPC}(\PR)\ar[r]^{\SP\otimes_{\SPC}(-)}  &   \Mod_{\SP}(\PR)
}
\]
where the vertical arrows are fully faithful embeddings.
The upper horizontal arrow is the forgetful functor, which is also regarded as
the functor given informally by $(\CCC,\CCC^{\le0})\mapsto \SP(\CCC^{\le0})$.
The forgetful functor $\PRTT\to \PRST$
plays  the role of the canonical functor $\SPC\to \SP$.
This leads to the following definition:

\begin{Definition}
\label{underlyingdef}
For $\mathcal{E}\in 2\PRTT$, we define the underlying stable presentable $(\infty,2)$-category of $\mathcal{E}$
to be $\PRST\otimes_{\PRTT}\mathcal{E}$.
The assignment $\mathcal{E}\mapsto \PRST\otimes_{\PRTT}\mathcal{E}$ induces
a symmetric monoidal functor $\lambda_2:2\PRTT\to 2\PRST$.

This construction can be iterated:
if a symmetric monoidal functor $\lambda_{n-1}:(n-1)\PRTT\to (n-1)\PRST$
has been defined, we define $\lambda_{n}:n\PRTT\to n\PRST$
by $\lambda_n(\PPP)=(n-1)\PRST\otimes_{(n-1)\PRTT}\mathcal{P}$. 
For $\PPP\in n\PRTT$,
we refer to $\lambda_n(\PPP)$ as the underlying stable presentable $(\infty,n)$-category of $\PPP$.

If we replace $\PRTT\to \PRST$ with the forgetful functor $\PRTTT\to \PRST$,
then the same construction defines $(n-1)\PRTTT\to (n-1)\PRST$. This also induces
\[
(-)\otimes_{(n-1)\PRTTT}(n-1)\PRST :n\PRTTT\to n\PRST.
\]
For $\PPP\in n\PRTTT$,
we refer to $\PPP\otimes_{(n-1)\PRTTT}(n-1)\PRST$ as the underlying stable presentable $(\infty,n)$-category of $\PPP$.
\end{Definition}

\subsubsection{$k$-linear contexts}
Let $k$ be a connective commutative ring spectrum.
In subsequent sections, we work in $k$-linear contexts.
We generalize $t$-structured presentable $(\infty, n)$-categories to the $k$-linear setting.
This generalization is formal. 
We write  $\PRTT_{k}=\Mod_{\Mod_k}(\PRTT)$ and $\PRTTT_{k}=\Mod_{\Mod_k}(\PRTTT)$.
Here we abuse notation by writing $\Mod_k$ for $(\Mod_{k},\Mod_k^{\le0})\in \PRTTT$. 
Using Corollary~\ref{tstablecolimits} and Corollary~\ref{tstablecolimits2} together with \cite[Corollary 4.2.3.5]{HA}
we see that both $\PRTT_k$ and $\PRTTT_k$ belong to $\CAlg(\vCat)$.
For $n\ge0$, we can consider $n\PRTT_k$ and $n\PRTTT_k$. By convention, $1\PRTT_k=\PRTT_k$, $1\PRTTT_k=\PRTTT_k$, $0\PRTT_k=(\Mod_k,\Mod_k^{\le0})$,
 and $0\PRTTT_k=(\Mod_k,\Mod_k^{\le0})$.
We refer to an object of $n\PRTT_k$ (resp $n\PRTTT_k$) as a {\it $t$-structured $k$-linear presentable $(\infty,n)$-category}
(resp. a {\it complete $t$-structured $k$-linear presentable $(\infty,n)$-category}).

There is a completion functor
\[
n\widehat{L}=(-)\otimes_{(n-1)\PRTT_k}(n-1)\PRTTT_k:n\PRTT_k\longrightarrow n\PRTTT_k
\]
which we define in the same way.

The underlying $k$-linear presentable $(\infty,n)$-categories are defined by functors
\[
(-)\otimes_{(n-1)\PRTT_k}(n-1)\PR_k :n\PRTT_k\to n\PR_k,\ \ \ \ (-)\otimes_{(n-1)\PRTTT_k}(n-1)\PR_k :n\PRTTT_k\to n\PR_k
\]
which we define in the same way.
We have the diagram
\[
\xymatrix{
n\PRTT_k \ar[r]^{n\widehat{L}} \ar[rd]_{(-)\otimes_{(n-1)\PRTT_k}(n-1)\PR_k} & n\PRTTT_k  \ar[d]^{(-)\otimes_{(n-1)\PRTTT_k}(n-1)\PR_k} \\
     &   n\PR_k.
}
\]
Note that for $n\ge1$ this diagram does not commute.

\subsection{Higher module categories}

\label{thetasection}

Let $A$ be a commutative ring spectrum, that is, $A\in \CAlg(\SP)$.
Let us consider $\LMod_A(\SP)\simeq \Mod_A(\SP)$. This $\infty$-category naturally has a symmetric monoidal structure.
Moreover, it admits a right complete accessible $t$-structure (cf. Example~\ref{connectivetstr}, Remark~\ref{unifiedtstr}), which exhibits $\LMod_A(\SP)$ as an object of $\CAlg(\PRTT)$.
Next we consider $\LMod_{\LMod_A(\SP)}(\PRTT)$. This category also has a symmetric monoidal structure, which exhibits it as an object of 
$\CAlg(2\PRTT)$. We write $\PRTT_A=\LMod_{\LMod_A(\SP)}(\PRTT)$.
We consider $2\PRTT_A=\LMod_{\PRTT_A}(2\PRTT)$. We can regard it as an object of $\CAlg(3\PRTT)$.
In the next step we have
\[
3\PRTT_A=\LMod_{\LMod_{\LMod_{\LMod_A(\SP)}(\PRTT)}(2\PRTT)}(3\PRTT).
\]
We can repeat this procedure to obtain $m\PRTT_A\in \CAlg((m+1)\PRTT)$.

We will work in a slightly more general setting. 
Let $n$ be a positive integer.
Let $A\in \Alg_{n}(\SP)$.
We will associate $m\PRTT_A$, that is an $\mathbf{E}_{n-m-1}$-monoidal $t$-structured presentable $(\infty,m+1)$-category, i.e., an object of $\Alg_{n-m-1}((m+1)\PRTT)$.
The construction gives a functor
\[
\theta_{\SP,n\mapsto m}^+:\Alg_{n}(\SP)\longrightarrow \Alg_{n-m-1}((m+1)\PRTT)
\]
which carries $A$ to $m\PRTT_A$. The construction is given in Construction~\ref{thetaconst1}.

When $A$ is connective, the stable $\infty$-category $(\LMod_A, \LMod_A^{\le0})$ equipped with the standard $t$-structure is 
 both left and right complete so that $(\LMod_A, \LMod_A^{\le0})$ belongs to $\PRTTT$.
Thus, we can define $\PRTTT_A=\LMod_{\LMod_A(\SP)}(\PRTTT)$, $2\PRTTT=\LMod_{\LMod_{\LMod_A(\SP)}(\PRTTT)}(2\PRTTT)$ and so on.
In the same way as in the setting of $\PRTT$, we have
\[
\theta^{\pm}_{\SP,n\mapsto m}:\Alg_{n}(\SPC)\longrightarrow \Alg_{n-m-1}((m+1)\PRTT)
\]
which carries $A$ to $m\PRTTT_A$ (see Construction~\ref{thetaconst2}).

We can also consider a version of $\PRST$: there is a functor
\[
\theta_{\SP,n\mapsto m}:\Alg_{n}(\SP)\longrightarrow \Alg_{n-m-1}((m+1)\PRST)
\]
which carries $A$ to $m\PR_A$ (see Construction~\ref{thetaconst1}).

In what follows, we construct these functors and discuss relations between them.

\subsubsection{From algebras to higher module categories}

\label{PRA}

\begin{Construction}
\label{cconstruction}
Let $\CCC^\otimes\in \CAlg(\vCat)$ be a symmetric monoidal $\infty$-category which has $\UUUU$-small colimits, and the tensor product
functor preserves $\UUUU$-small colimits separately in each variable.

To $A\in \Alg_1(\CCC)$ we associate the $\infty$-category $\LMod_A(\CCC)$.
Note that $\LMod_A(\CCC)$ has a $\CCC^\otimes$-module structure.
The assignment $A\mapsto \LMod_A(\CCC)$ is functorial: There exists a functor
\[
\theta_{\CCC}:\Alg_1(\CCC)\longrightarrow \Mod_{\CCC}(\vCat)_{\CCC/}
\]
which carries $A$ to $\LMod_A(\CCC)$ equipped with 
an essentially unique $\CCC$-module map $\CCC\to \LMod_A(\CCC)$
sending the unit of $\CCC$ to the left $A$-module $A$ (see \cite[4.8.3]{HA}; right modules are adopted).
A morphism $A\to A'$ in $\Alg_1(\CCC)$ maps to the base change functor
$\LMod_A(\CCC)\to \LMod_{A'}(\CCC)$.

According to \cite[4.8.5.5, 4.8.5.16]{HA},
\begin{itemize}
\item $\theta_{\CCC}$ is fully faithful,

\item $\theta_{\CCC}$ is promoted to a symmetric monoidal functor.

\end{itemize}
Note that $n\CCC$ belongs to $\CAlg(\vCat)$.
We consider the sequence of the symmetric monoidal functors
\[
\Alg_1(n\CCC)\stackrel{\theta_{n\CCC}}{\longrightarrow}  \Mod_{n\CCC}(\vCat)_{n\CCC/}\longrightarrow \Mod_{n\CCC}(\vCat)
\]
such that the right one is the forgetful functor.
By abuse of notation, we usually write $\theta_{n\CCC}$ for the composite.
This composite carries $\AAA\in \Alg_1(n\CCC)$ to $\LMod_{\AAA}(n\CCC)$.
$\LMod_{\AAA}(n\CCC)$ is dualizable, and its dual is given by $\RMod_{\AAA}(n\CCC)=\LMod_{\AAA^{opp}}(n\CCC)$.
Here $\AAA^{opp}$ means the opposite algebra of $\AAA$.
This implies that $\LMod_{\AAA}(n\CCC)$ is $\kappa_0$-compact (cf. \cite[Proposition 5.1.10]{Ste}).
Thus, we have symmetric monoidal functors
\[
\theta_{n\CCC}:\Alg_1(n\CCC) \longrightarrow \big{(}(n+1)\CCC\big{)}_{n\CCC/}\longrightarrow (n+1)\CCC.
\]
such that the left one is fully faithful, and the right one is the forgetful functor.
By Dunn additivity theorem there is a canonical equivalence $\Alg_{m+1}(\CCC)\simeq \Alg_{m}(\Alg_1(\CCC))$.
Using this equivalence and $\theta_{n\CCC}$ we have
\[
\theta_{n\CCC,m}:\Alg_{m+1}(n\CCC)\simeq \Alg_m(\Alg_1(n\CCC)) \longrightarrow \Alg_{m}((n+1)\CCC).
\]
\end{Construction}

Our primary interest lies in the case $n\PRST, n\PRTT\in \CAlg(\vCat)$:

\begin{Definition}
Take $\CCC^\otimes$ to be $n\PRST$.
We have symmetric monoidal functors
\[
\theta_{n\PRST}:\Alg_1(n\PRST) \longrightarrow \big{(}(n+1)\PRST\big{)}_{n\PRST/}\longrightarrow (n+1)\PRST
\]
such that the left one is fully faithful, and the right one is the forgetful functor.
This induces
\[
\theta_{n\PRST,m}:\Alg_{m+1}(n\PRST) \longrightarrow \Alg_{m}((n+1)\PRST).
\]
\end{Definition}

\begin{Definition}
\label{thetadeftt}
Take $\CCC^\otimes$ to be $n\PRTT$.
We have symmetric monoidal functors
\[
\theta_{n\PRTT}:\Alg_1(n\PRTT) \longrightarrow \big{(}(n+1)\PRTT\big{)}_{n\PRTT/}\longrightarrow (n+1)\PRTT
\]
such that the left one is fully faithful, and the right one is the forgetful functor.
This induces
\[
\theta_{n\PRTT,m}:\Alg_{m+1}(n\PRTT) \longrightarrow \Alg_{m}((n+1)\PRTT).
\]
\end{Definition}

\begin{Remark}
\label{undermodule}
We consider $\theta_{n\PRST,m}:\Alg_{m+1}(n\PRST) \longrightarrow \Alg_{m}((n+1)\PRST)$.
The symmetric monoidal functors $\lambda_n:n\PRTT\to n\PRST$ and $\lambda_{n+1}:(n+1)\PRTT\to (n+1)\PRST$ in Definition~\ref{underlyingdef}
induce the diagram
\[
\xymatrix{
\Alg_{m+1}(n\PRTT) \ar[r]^{\theta_{n\PRTT,m}} \ar[d]& \Alg_{m}((n+1)\PRTT) \ar[d] \\
\Alg_{m+1}(n\PRST) \ar[r]^{\theta_{n\PRST,m}} & \Alg_{m}((n+1)\PRST)
}
\]
where the vertical functors are induced by $\lambda_n$ and $\lambda_{n+1}$.
There exists a canonical equivalence
\[
\LMod_{\mathcal{A}}(n\PRTT)\otimes_{n\PRTT}n\PRST\simeq \LMod_{\lambda_{n}(\mathcal{A})}(n\PRST)
\]
(cf. \cite[4.8.4.6]{HA})
so that 
the above diagram commutes up to canonical homotopy.
\end{Remark}

\begin{Definition}
If we take $\CCC^\otimes$ to be $n\PRTTT$, we then obtain a symmetric monoidal functor
\[
\theta_{n\PRTTT}:\Alg_1(n\PRTTT) \longrightarrow \big{(}(n+1)\PRTTT\big{)}_{n\PRTTT/}\longrightarrow (n+1)\PRTTT
\]
such that the left one is fully faithful, and the right one is the forgetful functor.
We also have
\[
\theta_{n\PRTTT,m}:\Alg_{m+1}(n\PRTTT) \longrightarrow \Alg_{m}((n+1)\PRTTT).
\]
\end{Definition}

\begin{Remark}
\label{comptheta}
Remember that 
\[
\LMod_{\mathcal{A}}(\PRTT)\otimes_{\PRTT}\PRTTT\simeq \LMod_{\mathcal{A}}(\PRTTT)
\]
(cf. \cite[4.8.4.6]{HA}).
There exists a diagram
\[
\xymatrix{
\Alg_{1}(n\PRTT) \ar[r]^(0.4){\theta_{n\PRTT}} \ar[d] & \big{(}(n+1)\PRTT\big{)}_{n\PRTT/} \ar[d] \\
\Alg_1(n\PRTTT) \ar[r]^(0.4){\theta_{n\PRTTT}} & \big{(}(n+1)\PRTTT\big{)}_{n\PRTTT/}
}
\]
which commutes up to canonical homotopy, where the vertical functors are determined by
completion functors $n\PRTT\to n\PRTTT$ and $(n+1)\PRTT\to (n+1)\PRTTT$
induced by the base change along $\PRTT\to \PRTTT$ (cf. \cite[4.8.5.1]{HA}).
\end{Remark}

\begin{Construction}
\label{tmodule}
We will construct
\[
\theta_{\SP}^{+}:\Alg_{1}(\SP)\to  (\PRTT)_{(\SP,\SPC)/},
\]
which carries $A$ to $(\LMod_A, \LMod_A^{\le0})$,
where $\LMod_A^{\le0}$ is the smallest full subcategory which contains $A$ and is closed under small colimits and extensions.
Recall that
$\theta_{\SP}:\Alg_1(\SP)\to  (\PRST)_{\SP/}$
is fully faithful, we denote the essential image by $T$. Let $S$ be the full subcategory of $(\PRTT)_{(\SP,\SPC)/}$
spanned by objects of the form $(\LMod_A, \LMod_A^{\le0})$ with $A\in\Alg_{1}(\SP)$.
The forgetful functor $\PRTT\to \PRST$ induces $f:S\to T$.
We claim that $f$ is an equivalence. 
We define $\theta_{\SP}^{+}$ to be the composition $g\circ \theta_{\SP}$
where $g$ is an inverse functor of $f$.
Clearly, $f$ is essentially surjective.
Note that every functor $F:\LMod_A\to \LMod_B$ in $(\PRST)_{\SP/}$
sends $A$ to $B$. Since $\LMod_A^{\le0}$ (resp. $\LMod_B^{\le0}$)
is the smallest full subcategory of $\LMod_A$ (resp. $\LMod_B$)
which contains $A$ (resp. $B$) and is closed under small colimits and extensions,
it follows that $F$ is right $t$-exact. Thus, $f$ is fully faithful.
Since $\theta_{\SP}$ and $\PRTT\to \PRST$ are promoted to symmetric monoidal functors,
thus $\theta_{\SP}^+$ is also promoted to a symmetric monoidal functor.

Suppose that $k$ is a commutative ring spectrum.
If we replace $\SP$ with $\Mod_k$, then $\theta_{\SP}^{+}$
is generalized to $\theta_k^+:\Alg_{1}(\Mod_k)\to  (\PRTT)_{(\Mod_k,\Mod_k^{\le0})/}$
in the same way.

If $A$ is connective, then $(\LMod_A,\LMod_A^{\le0})$ is a both left and right complete $t$-structure
(Example~\ref{connectivetstr}) so that $(\LMod_A,\LMod_A^\le0)$ belongs to $(\PRTTT)_{(\SP,\SPC)/}$.
We define
\[
\theta_{\SP}^{\pm}:\Alg_{1}(\SPC)\to  (\PRTTT)_{(\SP,\SPC)/}
\]
to be the restriction of $\theta_{\SP}^+$.
For a connective commutative ring spectrum $k$, we also have $\theta_k^\pm:\Alg_{1}(\Mod_k^{\le0})\to  (\PRTTT)_{(\Mod_k,\Mod_k^{\le0})/}$
in the same way.
\end{Construction}

\begin{Remark}
\label{tmodulerem}
From the construction,
$\Alg_1(\SP)\to (\PRST)_{\SP/}$
given informally by $A\mapsto (\LMod_A,A)$
factors as
\[
\Alg_{1}(\SP)\to  (\PRTT)_{(\SP,\SPC)/} \to (\PRT)_{(\SP,\SPC)/}\to  (\PRST)_{\SP/}
\]
where the second functor is the inclusion, and the third functor is the forgetful functor.
For three functors $\Alg_{1}(\SP)\to  (\PRTT)_{(\SP,\SPC)/} $ $\Alg_{1}(\SP) \to (\PRT)_{(\SP,\SPC)/}$ and $\Alg_{1}(\SP)\to  (\PRST)_{\SP/}$,
their right adjoint functors are given by the formula
$(\CCC,C)\mapsto \End_{\CCC}(C)^{op}$.
The composite $(\PRTT)_{(\SP,\SPC)/} \to  (\PRST)_{\SP/}$ has a right adjoint given by $(\CCC,C)\mapsto ((\CCC,\CCC^{\le0}=\CCC),C)$.
Here $\CCC^{\le0}=\CCC$ means that the connective part is the whole category.
In particular, right completion does not change the image in $\Alg_1(\SP)$ under the right adjoint functor  $ (\PRT)_{(\SP,\SPC)/}\to \Alg_1(\Mod_k)$.
\end{Remark}

\begin{Construction}
\label{thetaconst1}
For $m_1+m_2=m+1,\ 1\le m_2\le m+1$, using $\theta_{n\PRTT,m}$ we consider the sequence of symmetric monoidal fully faithful functors
\begin{eqnarray*}
\Alg_{m+1}(\SP)  &\longrightarrow&  \Alg_{m}((\PRTT)_{(\SP,\SPC)/})\simeq \Alg_{m}(\PRTT) \\
 &\longrightarrow&  \Alg_{m-1}(2\PRTT) \\
&\longrightarrow&  \Alg_{m-2}(3\PRTT) \\
&\cdots& \\
&\longrightarrow&   \Alg_{m_1-1}((m_2+1)\PRTT)
\end{eqnarray*}
where the first arrow is induced by $\theta_{\SP}^+$, and the subsequent arrows are induced by $\theta_{r\PRTT}$.
We denote the composite functor by 
\[
\theta_{\SP,m+1\to m_2}^+:\Alg_{m+1}(\SP)\to \Alg_{m_1-1}((m_2+1)\PRTT).
\]

Using $\theta_{n\PRST,m}$ instead of $\theta_{n\PRTT,m}$,
we also define 
\[
\theta_{\SP,m+1\to m_2}:\Alg_{m+1}(\SP)\to \Alg_{m_1-1}((m_2+1)\PRST).
\]
Note that $\theta_{k,m+1\to m_2}^+$ and $\theta_{k,m+1\to m_2}$
 are fully faithful.

For  $A\in \Alg_{m+1}(\SP)$ we write $m_2\PRTT_A$ for $\theta_{\SP,m+1\to m_2}^+(A)$.
For  $A\in \Alg_{m+1}(\SP)$ we write $m_2\PR_A$ for $\theta_{\SP,m+1\to m_2}(A)$.
For ease of notation,
we usually omit the monoidal structure from the notation.
\end{Construction}

\begin{Remark}
Let $\Alg_{m_1-1}((m_2+1)\PRTT) \to \Alg_{m_1-1}((m_2+1)\PRST)$ be the functor induced by $\lambda_{m_2+1}$ (cf. Definition~\ref{underlyingdef}).
From Remark~\ref{undermodule}, Construction~\ref{tmodule} and Remark~\ref{tmodulerem}
we see that there exists a canonical diagram
\[
\xymatrix{
\Alg_{m+1}(\SP)\ar[r]^{\theta^+_{\SP,m+1\to m_2}} \ar[rd]_{\theta_{\SP,m+1\to m_2}} &  \Alg_{m_1-1}((m_2+1)\PRTT) \ar[d] \\
  &  \Alg_{m_1-1}((m_2+1)\PRST)
}
\]
which commutes up to canonical homotopy. Namely, $m_2\PRTT_A\otimes_{m_2\PRTT}m_2\PRST\simeq m_2\PR_A$.
\end{Remark}

\begin{Construction}
\label{thetaconst2}
We consider the sequence of symmetric monoidal fully faithful functors
\begin{eqnarray*}
\Alg_{m+1}(\SPC)  &\longrightarrow&  \Alg_{m}((\PRTTT)_{(\SP,\SPC)/})\simeq \Alg_{m}(\PRTTT) \\
 &\longrightarrow&  \Alg_{m-1}(2\PRTTT) \\
&\longrightarrow&  \Alg_{m-2}(3\PRTTT) \\
&\cdots& \\
&\longrightarrow&   \Alg_{m_1-1}((m_2+1)\PRTTT)
\end{eqnarray*}
where the first arrow is induced by $\theta_{\SP}^\pm$, and the subsequent arrows are induced by $\theta_{r\PRTTT}$.
We denote the composite fully faithful functor by 
\[
\theta^\pm_{\SP,m+1\to m_2}:\Alg_{m+1}(\SPC)\to \Alg_{m_1-1}((m_2+1)\PRTTT).
\]

For  $A\in \Alg_{m+1}(\SPC)$ we write $m_2\PRTTT$ for $\theta_{\SP,m+1\to m_2}^{\pm}(A)$.
\end{Construction}

\begin{Remark}
For $A\in \Alg_{m+1}(\SPC)$, there is a canonical equivalence $(m_2+1)\widehat{L}(m_2\PRTT_A)=m_2\PRTT_A\otimes_{m_2\PRTT}m_2\PRTTT\simeq m_2\PRTTT_A$.
\end{Remark}

\subsubsection{$k$-linear contexts}
Let $k$ be a connective commutative ring spectrum.
We note that $(\Mod_k, \Mod_k^{\le0})$ belongs to 
$\CAlg(\PRTTT)$.
We set 
\[
\PRTT_k:=\Mod_{(\Mod_k,\Mod_k^{\le0})}(\PRTT)\ \ \ \  \textup{and}\ \ \ \  \PRTTT_k:=\Mod_{(\Mod_k,\Mod_k^{\le0})}(\PRTTT).
\]
We take $\CCC$ to be $\PRTT_k$ in Construction~\ref{cconstruction} instead of $\PRTT$,
then we get $n\PRTT_k\in\CAlg(\vCat)$. Similarly, if $\CCC$ is $\PRTTT_k$, we have $n\PRTTT_k\in\CAlg(\vCat)$.

We have symmetric monoidal functors
\[
\theta_{n\PRTT_k}:\Alg_1(n\PRTT_k) \longrightarrow \big{(}(n+1)\PRTT_k\big{)}_{n\PRTT_k/}\longrightarrow (n+1)\PRTT_k
\]
which induces
\[
\theta_{n\PRTT_k,m}:\Alg_{m+1}(n\PRTT_k) \longrightarrow \Alg_{m}((n+1)\PRTT_k)
\]
for $m\ge0$.

We define $\theta_{n\PRTTT_k}:\Alg_1(n\PRTTT_k) \rightarrow \big{(}(n+1)\PRTTT_k\big{)}_{n\PRTTT_k/}$
and $\theta_{n\PRTTT,m}:\Alg_{m+1}(n\PRTTT) \rightarrow \Alg_{m}((n+1)\PRTTT)$
in the same way.
Using them, we also define
\[
\theta_{k,m+1\to m_2}^+:\Alg_{m+1}(\Mod_k)\to \Alg_{m_1-1}((m_2+1)\PRTT_k).
\]
and 
\[
\theta_{k,m+1\to m_2}^\pm:\Alg_{m+1}(\Mod_k^{\le0})\to \Alg_{m_1-1}((m_2+1)\PRTTT_k).
\]
in a similar way to $\theta^+_{\SP,m+1\to m_2}$
and $\theta^\pm_{\SP,m+1\to m_2}$.
Similarly, in the setting of $\PR_k$ we have
\[
\theta_{k,m+1\to m_2}:\Alg_{m+1}(\Mod_k)\to \Alg_{m_1-1}((m_2+1)\PR_k).
\]
The functors $\theta_{k,m+1\to m_2}^+$, $\theta_{k,m+1\to m_2}^\pm$ and $\theta_{k,m+1\to m_2}$
are (promoted to) symmetric monoidal functors and are fully faithful.

In the $k$-linear setting, we continue to use the notation $m_2\PRTT_A$, $m_2\PRTTT_A$ and $m_2\PR_A$ 
for $\theta_{k,m+1\to m_2}(A)$, $\theta_{k,m+1\to m_2}^\pm(A)$ and $\theta_{k,m+1\to m_2}(A)$, respectively.
For instance, $2\PRTT_A$ is defined to be 
\[
\LMod_{\LMod_{\LMod_A}(\PRTT_k)}(2\PRTT_k).
\]
This notation overlaps with those in Construction~\ref{thetaconst1} and Construction~\ref{thetaconst2}.
However, this is harmless. Indeed, using definitions and \cite[7.1.3.1]{HA},
one sees that $m_2\PRTT_A$, $m_2\PRTTT_A$ and $m_2\PR_A$ in this $k$-linear setting
are equivalent to those in Construction~\ref{thetaconst1} and Construction~\ref{thetaconst2}.
By default, in the subsequent sections, we adopt $m_2\PRTT_A$, $m_2\PRTTT_A$ and $m_2\PR_A$ in the $k$-linear setting.


\section{Higher categorified Koszul duality}

\label{HCKDsec}

Let $k$ be a base field.

\subsection{Algebraization}

We introduce the notion of algebraization and discuss a relation with Koszul duals of $\eenu$-algebras.

\begin{Definition}
\label{plainalgebraizationdef}
Let $n\ge0$ be a nonnegative integer.
Let $\mathsf{P}$ be any of $\PR_k$ or $\PRTT_k$ or $\PRTTT_k$.
Let $(\LLL,\phi:n\mathsf{P}\to \LLL)$ be an object in $((n+1)\mathsf{P})_{n\mathsf{P}/}$.
Here $L$ is the image of $(n-1)\mathsf{P}\in n\mathsf{P}$ under $\phi$.
Conversely, this $L$ determines $\phi$ so that $(\LLL,\phi)$ amounts to the pair $(\LLL,L)$.
We consider the functor $j_{(\LLL,L)}:\Alg_{n+1}(\Mod_k)^{op}\to \wSSS$ given informally by $A\mapsto \Map_{((n+1)\mathsf{P})_{n\mathsf{P}/}}((n\mathsf{P}_A,(n-1)\mathsf{P}_A),(\LLL,L))$. We regard $(\LLL,L))$ as an object of $((n+1)\mathsf{P})_{n\mathsf{P}/}$.
Here $\mathsf{P}_A$ is $\PR_A$ or $\PRTT_A$ or $\PRTTT_A$, depending on the choice of $\mathsf{P}$.
If $j_{(\LLL,L)}$ is representable by some $C\in \Alg_{n+1}(\Mod_k)$, we say that $(\LLL,L)$ has an algebraization with respect to $\mathsf{P}$.
If the context is obvious, we say that $(\LLL,L)$ has an algebraization.
Let $(n\mathsf{P}_C,(n-1)\mathsf{P}_C) \to (\LLL,L)$ be the morphism
that corresponds to the identity $C\to C$ through the equivalence $\Map_{\Alg_{n+1}(\Mod_k)}(C,C)\simeq \Map_{((n+1)\mathsf{P})_{n\mathsf{P}/}}(n\mathsf{P}_C,(\LLL,L))$.
We refer to it as the algebraization functor (with respect to $\mathsf{P}$).
\end{Definition}

\begin{Notation}
Let $\mathsf{P}$ be any of $\PR_k$ or $\PRTT_k$ or $\PRTTT_k$.
Let $(\LLL,\phi:n\mathsf{P}\to \LLL)$ be an object in $((n+1)\mathsf{P})_{n\mathsf{P}/}$.
As explained in Definition~\ref{plainalgebraizationdef}, we usually use
a pair $(\LLL,L)$ to indicate an object of $(\LLL,\phi:n\mathsf{P}\to \LLL)$.
\end{Notation}

\begin{Definition}
\label{algebraizationdef}
Let $n\ge0$ be a nonnegative integer.
Let $\mathsf{P}$ be any of $\PR_k$ or $\PRTT_k$ or $\PRTTT_k$.
Let $\MMM^\otimes$ be an $\eone$-algebra object in $(n+1)\mathsf{P}$, that is, an object of $\Alg_1((n+1)\mathsf{P})$.
(For later use, more generally $\MMM^\otimes$ is allowed to be an $\eone$-algebra object in $\Mod_{n\mathsf{P}}(\vCat)$.)
We consider the functor $j_{\MMM^\otimes}:\Alg_{n+2}(\Mod_k)^{op}\to \wSSS$ given informally by $A\mapsto \Map_{\Alg_1((n+1)\mathsf{P})}(n\mathsf{P}_A,\MMM^\otimes)$.
Here $\mathsf{P}_A$ is $\PR_A$ or $\PRTT_A$ or $\PRTTT_A$.
If $j_{\MMM^\otimes}$ is representable by some $C\in \Alg_{n+2}(\Mod_k)$, we say that $\MMM^{\otimes}$ has an algebraization.
Let $n\mathsf{P}_C\to \MMM^\otimes$ be the morphism
that corresponds to the identity $C\to C$ through the equivalence $\Map_{\Alg_{n+2}(\Mod_k)}(C,C)\simeq \Map_{\Alg_1((n+1)\mathsf{P})}(n\mathsf{P}_C,\MMM^\otimes)$.
We refer to it as the algebraization functor  (with respect to $\mathsf{P}$).
\end{Definition}

\begin{Remark}
In the setting of Definition~\ref{algebraizationdef},
we suppose that  $\MMM^\otimes$ be an $\eone$-algebra object in $(n+1)\PRTTT_k$.
If we regard $\MMM^\otimes$ as an object of $\PRTT_k$, we have
\[
\Map_{\Alg_1((n+1)\PRTT_k)}(n\PRTT_A,\MMM^\otimes)\simeq \Map_{\Alg_1((n+1)\PRTTT_k)}(n\PRTTT_A,\MMM^\otimes).
\]
Note that $n\PRTTT_A=n\PRTT_A\otimes_{n\PRTT_k}n\PRTTT_k$.
If $\MMM^\otimes$ has an algebraization with respect to $\PRTT_k$, then $\MMM^\otimes$  has an algebraization with respect to $\PRTTT_k$.
\end{Remark}

We will relate the algebraizations and Koszul duals of algebras.
Let $A$ be an augmented $\mathbf{E}_{n+2}$-algebra in $\Mod_k$.
By the base change $n\PRTT_A\to n\PRTT_k$ induced by the augmentation $A\to k$,
we regard $n\PRTT_k$ as a left $n\PRTT_A$-module.
That is, $n\PRTT_k$ is an object of $\LMod_{n\PRTT_A}((n+1)\PRTT_k)\simeq (n+1)\PRTT_A$ (see Proposition~\ref{matomete} (1) for the equivalence).
We consider the endomorphism algebra $\End_{n\PRTT_A}^l(n\PRTT_k)$ of $n\PRTT_k\in (n+1)\PRTT_A$, defined in $\Alg_1( (n+1)\PRTT_k)$
(Lemma~\ref{endopre} for the existence).
Here the superscript $l$ and the subscript $n\PRTT_A$ in $\End_{n\PRTT_A}^l$ indicate ``left $n\PRTT_A$-linear functors".

\begin{Proposition}
\label{endoalgebraization}
The following holds.

\begin{enumerate}

\item 
Let $B$ be an augmented $\mathbf{E}_{n+1}$-algebra in $\Mod_k$.
By the base change $(n-1)\PRTT_B\to (n-1)\PRTT_k$ induced by the augmentation $B\to k$,
we regard $(n-1)\PRTT_k$ as a left $(n-1)\PRTT_B$-module.
Then  $(n\PRTT_B,(n-1)\PRTT_k)$ has an algebraization $C\in \Alg_{n+1}(\Mod_k)$.
Moreover, $C$ is equivalent to the $\mathbf{E}_{n+1}$-Koszul dual $\DD_{n+1}(B)$.

\item
Let $A$ be an augmented $\mathbf{E}_{n+2}$-algebra in $\Mod_k$.
Then there exists an endomorphism algebra $\End_{n\PRTT_A}^l(n\PRTT_k)$, and $\End_{n\PRTT_A}^l(n\PRTT_k)$ has an algebraization $C\in \Alg_{n+2}(\Mod_k)$.
Moreover, $C$ is equivalent to the $\mathbf{E}_{n+2}$-Koszul dual $\DD_{n+2}(A)$.

\end{enumerate}
\end{Proposition}

\begin{Remark}
The $\mathbf{E}_{n+2}$-algebra $C$ has a canonically defined augmentation $p_C:C\to k$,
and $(C,p_c)$ is equivalent to $\DD_{n+2}(A)$ as an augmented $\mathbf{E}_{n+2}$-algebra. See the proof of Proposition~\ref{endoalgebraization}.
\end{Remark}

The proof of Proposition~\ref{endoalgebraization} will require some preliminaries.

Recall that $n\PRTT_k$ belongs to $\CAlg(\vCat)$.
We have symmetric monoidal functors
\[
\theta_{n\PRTT_k}:\Alg_1(n\PRTT_k) \longrightarrow \big{(}(n+1)\PRTT_k\big{)}_{n\PRTT/}\longrightarrow (n+1)\PRTT_k
\]
such that the left one is fully faithful, and the right one is the forgetful functor (cf. Construction~\ref{cconstruction}).
This induces a fully faithful functor 
\[
\theta_{n\PRTT_k,m}:\Alg_{m+1}(n\PRTT_k) \longrightarrow \Alg_{m}((n+1)\PRTT_k).
\]
We consider whether or not $\theta_{n\PRTT_k,m}$ has a right adjoint.
From the construction $\theta_{n\PRTT_k,m}$ is extended to
\[
\Alg_{m+1}(n\PRTT_k) \to \Alg_{m}(\Mod_{n\PRTT_k}(\vCat))
\]
For $\DDD\in(n+1)\PRTT_k$ and for $D\in \DDD$, its Hom (or endomorphism algebra) object 
$\Hom_{\DDD}(D,D)$ exists in $\Mod_{(n-1)\PRTT_k}(\vCat)$.
We can ask if $\Hom_{\DDD}(D,D')$ exists in $n\PRTT_k$.
Motivated by this question, we consider the following condition:

\begin{Definition}
\label{goodendomorphism}
Let $n\ge0$ be a nonnegative integer.
Let $\mathsf{P}$ be any of $\PR_k$ or $\PRTT_k$ or $\PRTTT_k$.

\begin{enumerate}
\item
Let $(\LLL,L)$ be an object in $((n+1)\mathsf{P})_{n\mathsf{P}/}$.
If there exists an endomorphism algebra object of $L$ in $\Alg_1(n\mathsf{P})$, 
we say that $(\LLL,L)$ has an endomorphism algebra.
If exists,
we write $\mathsf{E}((\LLL,L))$ for an endomorphism algebra.

Let $\uni$ denote a unit object in $\mathsf{E}((\LLL,L))$.
Suppose that $(\mathsf{E}((\LLL,L)),\uni)$ also has an endomorphism algebra which we denote by $\mathsf{E}^2((\LLL,L))\in \Alg_2((n-1)\mathsf{P}),\ldots,
\mathsf{E}^{n-1}((\LLL,L))$ with a unit also has an endomorphism algebra which we denote by $\mathsf{E}^{n}((\LLL,L))\in \Alg_n(\mathsf{P})$.
In this case, we say that $(\LLL,L)$ is a good endomorphism algebra.

\item
Let $\MMM^\otimes$ be an $\eone$-algebra object in $(n+1)\mathsf{P}$, that is, an object of $\Alg_1((n+1)\mathsf{P})$.
Let $\uni_{\MMM}$ denote a unit object of $\MMM^\otimes$.
If $(\MMM,\uni_{\MMM})$ has an endomorphism algebra, 
we say that $\MMM^\otimes$ has an endomorphism algebra of a unit.
If exists,
we write $\mathsf{E}(\MMM^\otimes)$ for an endomorphism algebra object of the unit.
Suppose that $\mathsf{E}(\MMM^\otimes)$ also has an endomorphism algebra of a unit which we denote by $\mathsf{E}^2(\MMM^\otimes)\in \Alg_3((n-1)\mathsf{P}),\ldots,\mathsf{E}^{n-1}(\MMM^\otimes)$ also has an endomorphism algebra of a unit which we denote by $\mathsf{E}^{n}(\MMM^\otimes)\in \Alg_{n+1}(\mathsf{P})$.
In this case, we say that $\MMM^\otimes$ has a good endomorphism algebra of a unit.

\end{enumerate}
\end{Definition}

\begin{Lemma}
\label{goodalgebraization}
Let $n\ge0$ be a nonnegative integer.
Let $\mathsf{P}$ be either $\PRTT_k$ or $\PR_k$.
The following holds.
\begin{enumerate}
\item
Let $(\LLL,L)$ be an object in $((n+1)\mathsf{P})_{n\mathsf{P}/}$.
Suppose that $(\LLL,L)$ has a good endomorphism algebra.
Then $(\LLL,L)$ has an algebraization.

\item
Let $\MMM^\otimes$ be an $\eone$-algebra object in $(n+1)\mathsf{P}$.
Suppose that $\MMM^\otimes$ has a good endomorphism algebra of a unit.
Then $\MMM^\otimes$ has an algebraization.
\end{enumerate}
\end{Lemma}

\Proof
We treat the case of $\mathsf{P}=\PRTT_k$; the case of $\mathsf{P}=\PR_k$ is similar. 
We first prove (1).
We treat the case of $n=0$.
In this situation, $(\LLL,L)\in (\PRTT_k)_{\Mod_k/}$.
For any $L,L'\in \LLL$ there exists a Hom object $\Hom_{\LLL}(L,L')$ in $\Mod_k$
since $\Mod_k$ is presentable.
Recall that $\theta_{0\PR_k,0}:\Alg_1(\Mod_k)\to (\PR_k)_{\Mod_k/}$ is extended to
$\theta_{0\PRTT_k,0}:\Alg_1(\Mod_k)\to (\PRTT_k)_{\Mod_k/}$ (see Construction~\ref{tmodule}).
The functor $\theta_{0\PR_k,0}$ has a right adjoint which carries $(\mathcal{A},A)$
to the opposite of the endomorphism algebra $\End(A)^{op}$ of the object $A$.
By the definition of the connective part of the standard $t$-structure (cf. Remark~\ref{unifiedtstr}),
$\theta_{0\PRTT_k,0}$ also has a right adjoint $E^+:(\PRTT_k)_{\Mod_k/}\to \Alg_1(\Mod_k)$
given informally by $(\mathcal{A},A)\mapsto \End(A)^{op}$.
Note that $\theta_{0\PRTT_k,0}$ is fully faithful.
Thus for any $A\in \Alg_1(\Mod_k)$ we have
\[
 \Map_{(\PRTT_k)_{\Mod_k/}}((\LMod_A,A),(\LLL,L))\simeq \Map_{\Alg_1(\Mod_k)}(A,\End(L)^{op}).
\]
Thus, the claim holds in the case of $n=0$.

We will show the case of $n\ge1$.
Let $((n+1)\PRTTT_k)_{n\PRTT_k/}^\star$
be the full subcategory of $((n+1)\PRTT_k)_{n\PRTT_k/}$
spanned by $f:n\PRTT_k\to \mathcal{B}$ such that a Hom object $\Hom(B,B)$ exists in $n\PRTT_k$ 
where $B:=f((n-1)\PRTT)$.
Note that the essential image of the fully faithful functor
$\theta_{n\PRTT_k,0}:\Alg_1(n\PRTT_k)\to ((n+1)\PRTT_k)_{n\PRTT_k/}$
is contained in $((n+1)\PRTT_k)_{n\PRTT_k/}^\star$
since for any $C\in \Alg_1(n\PRTT_k)$ the free module $C$ in $\LMod_C(n\PRTT_k)$ has the Hom object $C=\Hom(C,C)$.
Observe that the restricted functor 
$\theta_{n\PRTT_k,0}^\star:\Alg_1(n\PRTT_k)\to ((n+1)\PRTT_k)^\star_{n\PRTT_k/}$
has a right adjoint $E_{n\PRTT_k}^\star$.
For this purpose, it will suffice to show that for any $[f:n\PRTT_k\to \mathcal{B}]\in ((n+1)\PRTT_k)^\star_{n\PRTT_k/}$, the object $B:=f((n-1)\PRTT)$ has an endomorphism algebra defined as an object of $\Alg_1(n\PRTT_k)$.
According to \cite[4.7.1.40, 4.7.1.41]{HA}, this is equivalent to the condition that there exists a Hom object $\Hom(B,B)$ in $n\PRTT_k$. However, this is our assumption.
Since $(\LLL,L)$ has a good endomorphism algebra,
the image $E_{n\PRTT_k}^\star(\LLL,L)$
lies in $\Alg_1(n\PRTT_k)^{\star}$.
Its image under the functor
$E_{(n-1)\PRTT_k}^\star:\Alg_1(n\PRTT_k)^\star\to \Alg_2((n-1)\PRTT_k)$
lies in $\Alg_2((n-1)\PRTT_k)^\star$.
Iterating this procedure, we get the image $D_1$ in $\Alg_n(\PRTT_k)$.
Namely, if we denote by $((n+1)\PRTT_k)_{n\PRTT_k/}^{\star \star}$ the full subcategory spanned by those objects
$(\LLL,L)$ such that $(\LLL,L)$ has a good endomorphism algebra, then the composition yields the adjoint pair  
\[
U:\Alg_n(\PRTT_k) \rightleftarrows ((n+1)\PRTT_k)_{n\PRTT_k/}^{\star \star}:V.
\]
Consider the adjoint pair
\[
\Alg_{n}(\theta_{0\PRTT_k,0}):\Alg_{n+1}(\Mod_k)\rightleftarrows \Alg_n(\PRTT_k):\Alg_{n}(E^+) 
\]
induced by $(\theta_{0\PRTT_k,0},E^+)$.
Applying
$\Alg_{n}(E^+):\Alg_n(\PRTT_k)\to \Alg_{n+1}(\Mod_k)$
to $D_1$, we have $D_0:=\Alg_{n}(E^+)(D_1) \in\Alg_{n+1}(\Mod_k)$.
Using the sequence of adjoint pairs we have the canonical equivalence
\[
\Map_{\Alg_{n+1}(\Mod_k)}(C,D_0)\simeq \Map_{\Alg_n(\PRTT_k)}(\LMod_C,D_1)\simeq \Map_{ ((n+1)\PRTT_k)_{n\PRTT_k/}}((n\PRTT_{C},(n-1)\PRTT_C),(\LLL,L))
\]
(notice that $U(\LMod_C)\simeq n\PRTT_C$).
Thus $(\LLL,L)$ has an algebraization.

We next prove (2).
The proof proceeds by passing to algebra objects of the constructions in (1).
In the case of $n=0$, applying $\Alg_1$ to the adjoint pair $(\theta_{0\PRTT_k,0}, E^+)$ yields the adjoint pair 
\[
\theta_{0\PRTT_k,1}:\Alg_2(\Mod_k)\rightleftarrows \Alg_1(\PRTT_k):E^+,
\]
where the right adjoint $E^+$ is informally given by $\mathcal{A}^\otimes \mapsto \End(\uni_{\mathcal{A}})^{op}$.
Since $\theta_{0\PRTT_k,1}$ is fully faithful, the claim holds as in (1).

For $n\ge1$, according to \cite[4.7.1.40, 4.7.1.41]{HA}, the existence of an endomorphism algebra as an object of $\Alg_1(n\PRTT_k)$ is equivalent to the existence of a Hom object  in $n\PRTT_k$, which is exactly the defining condition of $((n+1)\PRTT_k)_{n\PRTT_k/}^\star$.
Consequently, we have a canonical equivalence of full subcategories $\Alg_1((n+1)\PRTTT_k)^\star \simeq \Alg_1(((n+1)\PRTT_k)_{n\PRTT_k/}^\star)$, where the left hand side is spanned by those objects having a good endomorphism algebra of a unit.
Taking algebra objects of the adjoint pair $(\theta_{n\PRTT_k,0}^\star,E_{n\PRTT_k}^\star)$ from (1) gives the adjoint pair
\[
\Alg_1(\theta_{n\PRTT_k,0}^\star):\Alg_2(n\PRTT_k)\rightleftarrows \Alg_1((n+1)\PRTT_k)^\star:\Alg_1(E_{n\PRTT_k}^\star).
\]
Iterating this procedure 
\[
\xymatrix{
\Alg_{n+2}(\Mod_k)   \ar@<0.5ex>[r]^{\theta_{0\PRTT_k,1}}  & \Alg_{n+1}(\PRTT_k)  \ar@<0.5ex>[l]   \ar@<0.5ex>[r]^(0.4){U'}  &  \Alg_1((n+1)\PRTT_k)^{\star \star}    \ar@<0.5ex>[l]^(0.6){V'} 
}
\]
and evaluating the resulting sequence of right adjoints on $\MMM^\otimes$ 
to obtain $D_1 \in \Alg_{n+1}(\PRTT_k)$ and $D_0 \in \Alg_{n+2}(\Mod_k)$, we deduce the canonical equivalence
\[
\Map_{\Alg_{n+2}(\Mod_k)}(C,D_0)\simeq \Map_{\Alg_{n+1}(\PRTT_k)}(\LMod_C,D_1)\simeq \Map_{ \Alg_1((n+1)\PRTT_k)}(n\PRTT_{C},\MMM^\otimes)
\]
exactly as in (1). Thus $\MMM^\otimes$ has an algebraization.
\QED

\begin{Proposition}
\label{matomete}
\begin{enumerate}
  \item Let $A$ be an $\mathbf{E}_{n+2}$-algebra in $\Mod_k$.
Then there exists a canonical equivalence of monoidal $\infty$-categories
\[
n\PRTT_A\simeq\LMod_{(n-1)\PRTT_A}(n\PRTT_k).
\]
Similarly, $n\PRTTT_A\simeq \LMod_{(n-1)\PRTTT_A}(n\PRTTT_k)$.

\item Let $A\to B$ be a morphism of $\mathbf{E}_{n+2}$-algebras in $\Mod_k$.
The induced morphism $n\PRTT_A\to n\PRTT_B$ in $\Alg_1((n+1)\PRTT_k)$
defines the restriction functor $\LMod_{n\PRTT_B}(\vCat)\to \LMod_{n\PRTT_A}(\vCat)$.
The image of $n\PRTT_B$ lies in $(n+1)\PRTT_A$.
Namely, $n\PRTT_B$ is $\kappa_0$-compact in $\LMod_{n\PRTT_A}(\vCat)$.
Similarly, $n\PRTTT_B$ is $\kappa_0$-compact in $\LMod_{n\PRTTT_A}(\vCat)$.

\item In the situation of (2),
the induced morphism $(n-1)\PRTT_A\to (n-1)\PRTT_B$ in $\Alg_2(n\PRTT_k)$ exhibits $(n-1)\PRTT_B$ as an algebra object in $\LMod_{(n-1)\PRTT_A}(n\PRTT_k)$.
There exist canonical equivalences
\[
n\PRTT_B\simeq\LMod_{(n-1)\PRTT_B}(n\PRTT_k)\simeq \LMod_{(n-1)\PRTT_B}((\LMod_{(n-1)\PRTT_A}(n\PRTT_k))\simeq\LMod_{(n-1)\PRTT_B}(n\PRTT_A).
\]
Similarly, $n\PRTTT_B\simeq \LMod_{(n-1)\PRTTT_B}(n\PRTTT_A)$.

\end{enumerate}

\end{Proposition}

\Proof
We treat the case of $\PRTT$. The case of $\PRTTT$ is the same.
We prove (1). We first note that $(n-1)\PRTT_A$ belongs to $n\PRTT_k$
since 
$(n-1)\PRTT_A=\LMod_{(n-2)\PRTT_A}((n-1)\PRTT_k)$ is dualizable and $\kappa_0$-compact in $\Mod_{(n-1)\PRTT}(\vCat)$ (see \cite[Proposition 5.1.10]{Ste}).
We observe an equivalence 
\[
\Mod_{(n-1)\PRTT_A}(\vCat)^{\kappa_0}\simeq \Mod_{(n-1)\PRTT_A}(n\PRTT_k). 
\]
This right-hand side is $\Mod_{(n-1)\PRTT_A}(\Mod_{(n-1)\PRTT_k}(\vCat)^{\kappa_0})$, and
\begin{eqnarray*}
\Mod_{(n-1)\PRTT_A}(\Mod_{(n-1)\PRTT_k}(\vCat)^{\kappa_0}) &\subset& \Mod_{(n-1)\PRTT_A}(\Mod_{(n-1)\PRTT_k}(\vCat))  \\
&\simeq& \Mod_{(n-1)\PRTT_A}(\vCat) 
\end{eqnarray*}
where $\subset$ means the inclusion as a full subcategory (see \cite[7.1.3.1, 7.1.3.3]{HA} for the latter monoidal equivalence).
We note that
$\Mod_{(n-1)\PRTT_A}(\Mod_{(n-1)\PRTT_k}(\vCat)^{\kappa_0})$ is the smallest full subcategory of
$\Mod_{(n-1)\PRTT_A}(\vCat)$ which contains the collection of objects
$\{(n-1)\PRTT_A\otimes_{(n-1)\PRTT_k}(n-1)\PRTT_k\otimes \PPP\}_{\PPP\in \PR}$
and is closed under small colimits (see \cite[Proposition 5.1.7]{Ste}).
Again by \cite[Proposition 5.1.7]{Ste},
$\Mod_{(n-1)\PRTT_A}(\vCat)^{\kappa_0}$ is the smallest full subcategory of
$\Mod_{(n-1)\PRTT_A}(\vCat)$ which contains the collection of objects
$\{(n-1)\PRTT_A\otimes \PPP\}_{\PPP\in \PR}$
and is closed under small colimits. Thus, we see the equivalence.

We prove (2) and (3). By (1),
$\phi:\RMod_{(n-1)\PRTT_A}(n\PRTT_k)\simeq \LMod_{(n-1)\PRTT_A}(n\PRTT_k)\simeq n\PRTT_A$
(the first equivalence follows from the fact that $(n-1)\PRTT_A$ is $\etwo$-algebra).
Thus we have
\begin{eqnarray*}
\RMod_{(n-1)\PRTT_B}(n\PRTT_A) &\simeq& \LMod_{(n-1)\PRTT_B}(n\PRTT_A) \\
&\simeq& \RMod_{(n-1)\PRTT_B}(\RMod_{(n-1)\PRTT_A}(n\PRTT_k)) \\
&\simeq& \RMod_{(n-1)\PRTT_B}(n\PRTT_k)  \\
&\simeq& \LMod_{(n-1)\PRTT_B}(n\PRTT_k)= n\PRTT_B.
\end{eqnarray*}
where the first and 4-th equivalences follow from the fact that $(n-1)\PRTT_B$ is a $\etwo$-algebra,
$\phi$ induces the second equivalence, and the third equivalence is induced by the forgetful functor
(see \cite[7.1.3.1]{HA}).
By \cite[Proposition 5.1.10]{Ste}, $\RMod_{(n-1)\PRTT_B}(n\PRTT_A)$ is $\kappa_0$-compact in $\LMod_{n\PRTT_A}(\vCat)$.
Thus, $n\PRTT_B\in (n+1)\PRTT_A$.
\QED

\begin{Lemma}
\label{endopre}
Let $A$ be an augmented $\mathbf{E}_{n+2}$-algebra.
There exists an endomorphism algebra object
$\End_{n\PRTT_A}^l(n\PRTT_k)\in \Alg_1((n+1)\PRTT_k)$ of $((n+1)\PRTT_A,n\PRTT_k\in (n+1)\PRTT_A)$.
Similarly, there exists an endomorphism algebra object
$\End_{n\PRTTT_A}^l(n\PRTTT_k)\in \Alg_1((n+1)\PRTTT_k)$ of $((n+1)\PRTTT_A, n\PRTTT_k\in (n+1)\PRTTT_A)$.
\end{Lemma}

\Proof
The case of $\End_{n\PRTTT_A}^l(n\PRTTT_k)$ is the same with $\End_{n\PRTT_A}^l(n\PRTT_k)$.
We prove the case of $\End_{n\PRTT_A}^l(n\PRTT_k)$.
By \cite[4.7.1.40, 4.7.1.41]{HA} (or \cite[3.2.2.4]{HA}),
it is enough to prove that there exists a
Hom object $E=\Hom_{n\PRTT_A}^l(n\PRTT_k,n\PRTT_k)$ in $(n+1)\PRTT_k$.
Since $\Mod_{n\PRTT_k}(\vCat)$ is a $\mathbb{V}$-presentable $\infty$-category,
$E$ exists in $\Mod_{n\PRTT_k}(\vCat)$.
There exist equivalences in  $\Mod_{n\PRTT_k}(\vCat)$:
\begin{eqnarray*}
E   &\simeq& \Hom_{n\PRTT_k}^l(n\PRTT_k\otimes_{n\PRTT_A}n\PRTT_k,n\PRTT_k)  \\
&\simeq&  \Hom_{n\PRTT_k}^l(n\PRTT_{k\otimes_Ak},n\PRTT_k)   \\
&\simeq&  \Hom_{n\PRTT_k}^l(\LMod_{(n-1)\PRTT_{k\otimes_Ak}}(n\PRTT_{k}),n\PRTT_k)   \\
&\simeq&  \Hom_{n\PRTT_k}^l(\RMod_{(n-1)\PRTT_{k\otimes_Ak}}(n\PRTT_{k}),n\PRTT_k)   \\
&\simeq&  \LMod_{(n-1)\PRTT_{k\otimes_Ak}}(n\PRTT_{k})
\end{eqnarray*}
where the first one is induced by base change along $n\PRTT_A\to n\PRTT_k$, the second one comes from Lemma~\ref{categorytensor},
the third equivalence arises from the symmetry of bar construction,
and the final equivalence is induced by the $\PRTT_k$-linear functor
$\LMod_{(n-1)\PRTT_{k\otimes_Ak}}(n\PRTT_{k})\otimes_{\PRTT_k} \RMod_{(n-1)\PRTT_{k\otimes_Ak}}(n\PRTT_{k})\to \PRTT_k$
induced by the integral kernel of the diagonal bimodule $(n-1)\PRTT_{k\otimes_Ak}$ (cf.\ \cite[4.8.4.1, 4.8.4.8]{HA}).
By \cite[Proposition 5.1.10]{Ste}, $\LMod_{(n-1)\PRTT_{k\otimes_Ak}}(n\PRTT_{k})$ belongs to $(n+1)\PRTT_k$.
\QED

\begin{Lemma}
\label{goodendoendo}
Suppose  $n\ge1$. The following holds.
\begin{enumerate}
\item
Let $B$ be an augmented $\mathbf{E}_{n+1}$-algebra.
The pair $(n\PRTT_B,(n-1)\PRTT_k)$ has a good endomorphism algebra.

\item
Let $A$ be an augmented $\mathbf{E}_{n+2}$-algebra.
The endomorphism algebra $\End_{n\PRTT_A}^l(n\PRTT_k)\in \Alg_1((n+1)\PRTTT_k)$ has a good endomorphism algebra of a unit.
\end{enumerate}
\end{Lemma}

\Proof
We prove (2). The proof of (1) is similar to that of (2): the argument works by replacing $k\otimes_Ak$ with $B$.
The proof of (2) is a generalization of the proof of Lemma~\ref{endopre}.
From the proof of Lemma~\ref{endopre}, there exists an equivalence
$\Hom_{n\PRTT_A}^l(n\PRTT_k,n\PRTT_k)\simeq n\PRTT_{k\otimes_Ak}$
in $(n+1)\PRTT_k$. This equivalence
sends the identity functor of $n\PRTT_k$ to the left $(n-1)\PRTT_{k\otimes_Ak}$-module $(n-1)\PRTT_{k}$,
which is the integral kernel of the identity functor.
We write $\textup{Bar}(A)$ for $k\otimes_Ak$. We regard it as an $\mathbf{E}_{n+1}$-algebra in $\Mod_k$, that is defined by
the bar construction (more generally, it is an $\mathbf{E}_{n+1}$-algebra in $\textup{CoAlg}_1(\Mod_k)$).
Next, we prove that there exists an endomorphism algebra
$\End_{(n-1)\PRTT_{\textup{Bar}(A)}}^l((n-1)\PRTT_k)$ in $\Alg_1(n\PRTT_k)$.
As in the proof of Lemma~\ref{endopre},
it is enough to show that there exists a Hom object
$\Hom_{(n-1)\PRTT_\textup{Bar}(A)}^l((n-1)\PRTT_k,(n-1)\PRTT_k)$
in $n\PRTT_k$.
By the same argument, 
we have an equivalence
\[
\Hom_{(n-1)\PRTT_{\textup{Bar}(A)}}^l((n-1)\PRTT_k,(n-1)\PRTT_k)\simeq (n-1)\PRTT_{\textup{Bar}^2(A)}
\]
in $n\PRTT_k$, where $\textup{Bar}^2(A)=\textup{Bar}(\textup{Bar}(A))=k\otimes_{\textup{Bar}(A)}k$ is an $\mathbf{E}_n$-algebra obtained by the bar construction.
This equivalence
maps the identity functor of $(n-1)\PRTT_k$ to the left $(n-2)\PRTT_{\textup{Bar}^2}(A)$-module $(n-2)\PRTT_{k}$.
Thus, $\mathsf{E}(\End_{n\PRTT_A}^l(n\PRTT_k))=\End_{(n-1)\PRTT_{\textup{Bar}(A)}}^l((n-1)\PRTT_k)$ exists in $\Alg_1(n\PRTT_k)$.
Repeating this argument, we see that for $1\le r\le n-1$ there exists 
$\mathsf{E}^r(\End_{n\PRTT_A}^l(n\PRTT_k))$ in $\Alg_1((n+1-r)\PRTT_k)$, whose underlying $(n-r)\PRTT_k$-module is 
$(n-r)\PRTT_{\textup{Bar}^{r+1}(A)}$.
In particular, 
$\mathsf{E}^{n-1}(\End_{n\PRTT_A}^l(n\PRTT_k))\simeq \PRTT_{\textup{Bar}^{n}(A)}=\LMod_{\LMod_{\textup{Bar}^{n}(A)}}(\PRTT_k)$ in $2\PRTT_k$.
Here $\LMod_{\textup{Bar}^{n}(A)}$ is regarded as an object of $\Alg_1(\PRTT_k)$.
Thus, it is enough to show that there exists an endomorphism algebra of $\Mod_k\in \PRTT_{\textup{Bar}^{n}(A)}$ in $\Alg_1(\PRTT_k)$.
The existence follows from Lemma~\ref{specialdegree}. 
Consequently, our assertion follows. 

\QED

\begin{Lemma}
\label{specialdegree}
Let $B$ be an augmented $\mathbf{E}_{2}$-algebra,
$\PRTT_B= \LMod_{\LMod_{B}}(\PRTT_k)$ in $2\PRTT_k$.
There exists an endomorphism algebra of $\Mod_k\in \PRTT_{B}$ in $\Alg_1(\PRTT_k)$.
\end{Lemma}

\Proof
There exists a Hom object from $\Mod_k$ to itself in $\PRT_k$, defined by
\[
H=(\Fun_{\LMod_B}(\Mod_k,\Mod_k), \Fun_{\LMod_B^{\le0}}(\Mod_k^{\le0},\Mod_k^{\le0})).
\]
Note that it is generally not right complete (cf. \cite{CKDA}).
Passing to $\widehat{R}:\PRT_k\to \PRTT_k$ which is a right adjoint of the symmetric monoidal inclusion $\PRTT_k\subset \PRT_k$,
we obtain $\widehat{R}(H)$, that is a Hom object defined in $\PRTT_k$.
Finally, by \cite[4.7.1.40, 4.7.1.41]{HA}, the Hom object $\widehat{R}(H)$ is
(uniquely) promoted to an endomorphism algebra defined in $\Alg_1(\PRTT_k)$.
\QED

\begin{Lemma}
\label{compff}
Let $A$ and $B$ be $\mathbf{E}_{n+1}$-algebras over $k$. 
Then the functor $\theta_{k,n+1\to n}^+:\Alg_{n+1}(\Mod_k)\to ((n+1)\PRTT_k)_{n\PRTT_k/}$ given by $R\mapsto (n\PRTT_R,(n-1)\PRTT_R)$ induces an equivalence
\[
\Map_{\Alg_{n+1}(\Mod_k)}(A,B)\to \Map_{((n+1)\PRTT_k)_{n\PRTT_k/}}((n\PRTTT_A,(n-1)\PRTT_A),(n\PRTT_B,(n-1)\PRTT_B)).
\]

Suppose that $A$ and $B$ are promoted to $\mathbf{E}_{n+2}$-algebras. 
Then the functor $\Alg_1(\theta_{k,n+1\to n}^+):\Alg_{n+2}(\Mod_k)\to\Alg_1((n+1)\PRTT_k)$ given by $R\mapsto n\PRTT_R$ induces an equivalence
\[
\Map_{\Alg_{n+2}(\Mod_k)}(A,B)\to \Map_{\Alg_1((n+1)\PRTT_k)}(n\PRTT_A,n\PRTT_B).
\]
Suppose furthermore that $B$ is connective.
Then 
the functor $\Alg_1(\theta_{k,n+1\to n}^\pm):\Alg_{n+2}(\Mod_k)\to\Alg_1((n+1)\PRTTT_k)$ given by $R\mapsto n\PRTTT_R$ induces an equivalence
\[
\Map_{\Alg_{n+2}(\Mod_k)}(A,B)\to \Map_{\Alg_1((n+1)\PRTTT_k)}(n\PRTTT_A,n\PRTTT_B).
\]

\end{Lemma}

\Proof
Recall the sequence of adjoint pairs from the proof of Lemma~\ref{goodalgebraization}:
\[
\Alg_{n+1}(\Mod_k)\stackrel{\theta_{0\PRTT_k,0}}{\rightleftarrows}   \Alg_{n}(\PRTT_k) \stackrel{U}{\rightleftarrows} ((n+1)\PRTTT_k)_{n\PRTT_k/}^{\star \star}.
\]
Here by abuse of notation we write $\theta_{0\PRTT_k,0}$  for $\Alg_n(\theta_{0\PRTT_k,0})$.
The functor in the statement is the composite $U\circ  \theta_{0\PRTT_k,0}$ of left adjoint functors.
The first assertion follows from the fact that $\theta_{0\PRTT_k,0}$ and $U$ are fully faithful.
Applying $\Alg_1(-)$ to the sequence of the above adjoint pairs,
we also see the second assertion.

It remains to prove the final assertion: We assume that
$B$ is connective.
Consider a variant of the above sequence of adjoint pairs:
\[
\Alg_{n+2}(\Mod_k)\stackrel{\theta_{0\PRTT_k,1}}{\rightleftarrows}   \Alg_{n+1}(\PRTT_k) \stackrel{\widehat{L}}{\rightleftarrows}\Alg_{n+1}(\PRTTT_k) \stackrel{U^{\pm}}{\rightleftarrows} \Alg_1((n+1)\PRTTT_k)^{\star \star}.
\]
The left adjoint $U^{\pm}$ is a $\PRTTT$-version of $U'$, that is defined in a similar way. 
Here by abuse of notation we write $\theta_{0\PRTT_k,1}$ and $\widehat{L}$ for $\Alg_1(\theta_{0\PRTT_k,1})$ and $\Alg_1(\widehat{L})$, respectively.
The map of mapping spaces is induced by the composite of left adjoint functors.
Note that the functors $\theta_{0\PRTT_k,1}$ and $U^{\pm}$ are fully faithful.
Thus, it is enough to show that 
$\widehat{L}$ induces an equivalence
\[
\Map_{\Alg_{n+1}(\PRTT_k)}(\LMod_A,\LMod_B)\to \Map_{\Alg_{n+1}(\PRTTT_k)}(\widehat{L}(\LMod_A),\widehat{L}(\LMod_B)).
\]
Indeed, the induced map is an equivalence since the middle adjoint pair is a localization
(whose right adjoint is fully faithful), and by the connectivity of $B$, the category $\LMod_B$ lies in $\Alg_{n+1}(\PRTTT_k)$.
\QED

{\it Proof of Proposition~\ref{endoalgebraization}.}
The existence of algebraizations of (1) and (2) follows from Lemma~\ref{goodalgebraization} and Lemma~\ref{goodendoendo}.

In the case (1), there exist $C\in \Alg_{n+1}(\Mod_k)$ and an algebraization functor $\alpha:(n\PRTT_C,(n-1)\PRTT_C)\to (n\PRTT_B,(n-1)\PRTT_k)$.
We will prove that $C$ is equivalent to $\DD_{n+1}(B)$.
Consider the composition $P_C:n\PRTT_C\stackrel{\alpha}{\to} n\PRTT_B\to n\PRTT_k$ with the forgetful functor $P_B:n\PRTT_B\to n\PRTT_k$.
Note that $P_C$ carries $(n-1)\PRTT_C$ to $(n-1)\PRTT_k$.
From Lemma~\ref{compff} we deduce that there exists an essentially unique augmentation $p_C:C\to k$ that induces/corresponds to $P_C$. 
By Lemma~\ref{compff}, for any $D\in \Alg_{n+1}(\Mod_k)$ there exists an equivalence in $\wSSS$:
\begin{eqnarray*}
\Map_{\Alg_{n+1}(\Mod_k)}(B\otimes_k D,k) \simeq \Map_{((n+1)\PRTT_k))_{n\PRTT_k/}}(n\PRTT_B\otimes n\PRTT_D,n\PRTT_k)
\end{eqnarray*}
By Lemma~\ref{compff} again, we have equivalences
$\Map_{\Alg_{n+1}(\Mod_k)}(B,k)\simeq \Map_{((n+1)\PRTT_k)_{n\PRTT_k/}}(n\PRTT_B,n\PRTT_k)$ and $\Map_{\Alg_{n+1}(\Mod_k)}(D,k)\simeq \Map_{((n+1)\PRTT_k)_{n\PRTT_k/}}(n\PRTT_D,n\PRTT_k)$.
Let $p_D:D \to k$ be an augmentation. By these equivalences, 
\begin{eqnarray*}
\textup{Pair}(p_B,p_D)=\Map_{\Alg_{n+1}(\Mod_k)}(B\otimes_k D,k)\times_{\Map(B,k)\times \Map(D,k)}\{(p_A,p_D)\}
\end{eqnarray*}
is equivalent to
\[
\mathbb{M}=\Map_{\Alg_1(n\PRTT_k)}((n-1)\PRTT_B\otimes (n-1)\PRTT_D,(n-1)\PRTT_k)\times_{ \Map((n-1)\PRTT_B,(n-1)\PRTT_k)\times \Map((n-1)\PRTT_D,(n-1)\PRTT_k)}\{(P_B,P_D)\}
\]
where $P_D:(n-1)\PRTT_D\to (n-1)\PRTT_k$ is induced by $p_D$ 
(we omit some subscripts in $\Map$
and also omit specified objects from the notation).
Since  $\End_{(n-1)\PRTT_B}^l((n-1)\PRTT_k)\in \Alg_1(n\PRTT_k)$ is an $\eone$-Koszul dual of $(n-1)\PRTT_B$ in $\Alg_1(n\PRTT_k)$ for $n\ge2$,
the latter space $\mathbb{M}$
is equivalent to $\Map_{\Alg_1^+(n\PRTT_k)}((n-1)\PRTT_D, \End_{(n-1)\PRTT_B}^l((n-1)\PRTT_k))$.
Since $\End_{(n-1)\PRTT_B}^l((n-1)\PRTT_k)$ has an algebraization $C$, 
we have
\[
\Map_{\Alg_1(n\PRTT_k)}((n-1)\PRTT_D, \End_{(n-1)\PRTT_B}^l((n-1)\PRTT_k))\simeq \Map_{\Alg_{n+1}(\Mod_k)}(D,C).
\]
When $n=1$,
we have
\[
\Map_{\Alg_1(\PRT_k)}(\Mod_D, \End_{\Mod_B}^l(\Mod_k))\simeq \Map_{\Alg_1(\PRTT_k)}(\Mod_D, \widehat{R}(\End_{\Mod_B}^l(\Mod_k)))\simeq \Map_{\Alg_{2}(\Mod_k)}(D,C).
\]
Combined with  $\Map_{\Alg_{n+1}(\Mod_k)}(D,k)\simeq \Map_{\Alg_1(n\PRTT_k)}((n-1)\PRTT_D,(n-1)\PRTT_k)$,
we have
\[
\textup{Pair}(p_B,p_D)\simeq \Map_{\Alg_1^+(n\PRTT_k)}((n-1)\PRTT_D, \End_{(n-1)\PRTT_B}^l((n-1)\PRTT_k))\simeq \Map_{\Alg^+_{n+1}(\Mod_k)}((D,p_D),(C,p_C)).
\]
It follows that $(C,p_C)\in \Alg_{n+1}^+(\Mod_k)$ is an $\mathbf{E}_{n+1}$-Koszul dual $\DD_{n+1}(B)$ of $B$.

We prove (2). First, by Lemma~\ref{endopre}, there exists an endomorphism algebra $\End_{n\PRTT_A}^l(n\PRTT_k)$.
The following argument is similar to (1).
In the case (2), there exist $C\in \Alg_{n+2}(\Mod_k)$ and an algebraization functor $n\PRTT_C\to \End_{n\PRTT_A}^l(n\PRTT_k)$.
It remains to prove that $C$ is equivalent to $\DD_{n+2}(A)$.
Consider the composite $P_C:n\PRTT_C\to \End_{n\PRTT_A}^l(n\PRTT_k)\to\End_{n\PRTT_k}^l(n\PRTT_k)=n\PRTT_k$ induced by the restriction along the augmentation $P_A:n\PRTT_A\to n\PRTT_k$.
From Lemma~\ref{compff}, we note that there exists an essentially unique augmentation $p_C:C\to k$ that induces/corresponds to $P_C$. 
By Lemma~\ref{compff}, there exists an equivalence in $\wSSS$.
Let $D\in \Alg_{n+2}(\Mod_k)$.
\begin{eqnarray*}
\Map_{\Alg_{n+2}(\Mod_k)}(A\otimes_k D,k) \simeq \Map_{\Alg_1((n+1)\PRTT_k)}(n\PRTT_A\otimes_{n\PRTT_k}n\PRTT_D,n\PRTT_k)
\end{eqnarray*}
By Lemma~\ref{compff} again, we have equivalences
$\Map_{\Alg_{n+2}(\Mod_k)}(A,k)\simeq \Map_{\Alg_1((n+1)\PRTT_k)}(n\PRTT_A,n\PRTT_k)$ and $\Map_{\Alg_{n+2}(\Mod_k)}(D,k)\simeq \Map_{\Alg_1((n+1)\PRTT_k)}(n\PRTT_D,n\PRTT_k)$.
Let $p_D:D\to k$ be an augmentation.
By these equivalences, 
\begin{eqnarray*}
\textup{Pair}(p_A,p_D)=\Map_{\Alg_{n+2}(\Mod_k)}(A\otimes_k D,k)\times_{\Map(A,k)\times \Map(C,k)}\{(p_A,p_D)\}
\end{eqnarray*}
is equivalent to
\[
\mathbf{M}=\Map_{\Alg_1((n+1)\PRTT_k)}(n\PRTT_A\otimes n\PRTT_D,n\PRTT_k)\times_{ \Map(n\PRTT_A,n\PRTT_k)\times \Map(n\PRTT_D,n\PRTT_k)}\{(P_A,P_D)\}
\]
where $P_D:n\PRTT_D\to n\PRTT_k$ is induced by $p_D$ 
(we omit some subscripts in $\Map$). Since  $\End_{n\PRTT_A}^l(n\PRTT_k)\in \Alg_1((n+1)\PRTT_k)$ is an $\eone$-Koszul dual of $n\PRTT_A$ in $\Alg_1((n+1)\PRTT_k)$,
then $\mathbf{M}$
is equivalent to $\Map_{\Alg_1^+((n+1)\PRTT_k)}(n\PRTT_D, \End_{n\PRTT_A}^l(n\PRTT_k))$.
Since $\End_{n\PRTT_A}^l(n\PRTT_k)$ has an algebraization $C$, 
we have
\[
\Map_{\Alg_1((n+1)\PRTT_k)}(n\PRTT_D, \End_{n\PRTT_A}^l(n\PRTT_k))\simeq \Map_{\Alg_{n+2}(\Mod_k)}(D,C).
\]
Combined with  $\Map_{\Alg_{n+2}(\Mod_k)}(D,k)\simeq \Map_{\Alg_1((n+1)\PRTT_k)}(n\PRTT_D,n\PRTT_k)$,
we have
\[
\textup{Pair}(p_A,p_D)\simeq \Map_{\Alg_1^+((n+1)\PRTT_k)}(n\PRTT_D, \End_{n\PRTT_A}^l(n\PRTT_k))\simeq \Map_{\Alg^+_{n+2}(\Mod_k)}((D,p_D),(C,p_C)).
\]
It follows that $(C,p_C)\in \Alg_{n+2}^+(\Mod_k)$ is an $\mathbf{E}_{n+2}$-Koszul dual $\DD_{n+2}(A)$ of $A$.
\QED

\begin{Remark}
\label{pairingalgebraization}
From the final portion of the above proof, the algebraization functor $n\PRTT_{\DD_{n+2}(A)}\to  \End_{n\PRTT_A}^l(n\PRTT_k)$
corresponds to the pairing $n\PRTT_A\otimes_{n\PRTT_k}n\PRTT_{\DD_{n+2}(A)}\to n\PRTT_k$ induced by the universal $\mathbf{E}_{n+2}$-Koszul dual pairing 
$A\otimes_k\DD_{n+2}(A)\to k$.
\end{Remark}

As a byproduct of the proof we reach a categorical understanding of the following well-known fact:

\begin{Corollary}
\label{barkoszuldual}
Let $A$ be an augmented $\mathbf{E}_{n+2}$-algebra.
Let $\DD_{n+2}(A)$ be an $\mathbf{E}_{n+2}$-Koszul dual of $A$.
Then there exists an equivalence $\DD_{n+2}(A)\simeq \DD_{n+1}(\textup{Bar}(A))$
in $\Alg_{n+1}(\Mod_k)$.
\end{Corollary}

\Proof
Recall from Proposition~\ref{endoalgebraization} and its proof
that $\DD_{n+2}(A)$ is the image of $\End_{n\PRTT_A}^l(n\PRTT_k)\in \Alg_1((n+1)\PRTT_k)$ under the right adjoint functors in the sequence of adjoint pairs
\[
\Alg_{n+2}(\Mod_k)\stackrel{\theta_{0\PRTT_k,0}}{\rightleftarrows}   \Alg_{n+1}(\PRTT_k) \stackrel{U}{\rightleftarrows} \Alg_1(((n+1)\PRTTT_k))^{\star \star}.
\]
See the proof of Lemma~\ref{goodalgebraization}
for this sequence. 
The equivalence
$\End_{n\PRTT_A}^l(n\PRTT_k)\simeq n\PRTT_k\otimes_{n\PRTT_A}n\PRTT_k\simeq n\PRTT_{k\otimes_Ak}$
carries the identity functor to $(n-1)\PRTT_k\in n\PRTT_{k\otimes_Ak}$ (see Lemma~\ref{categorytensor} for $n\PRTT_k\otimes_{n\PRTT_A}n\PRTT_k\simeq n\PRTT_{k\otimes_Ak}$).
We identify $(n\PRTT_{k\otimes_Ak}, (n-1)\PRTT_k)$ with the image of $\End_{n\PRTT_A}^l(n\PRTT_k)$ in $((n+1)\PRTTT_k)_{n\PRTT_k/}^{\star \star}$.
The above sequence of adjoint pairs is obtained from the following sequence of adjoint functors  
\[
\Alg_{n+1}(\Mod_k)\rightleftarrows   \Alg_{1}(n\PRTT_k)^{\star \star}\rightleftarrows ((n+1)\PRTTT_k)_{n\PRTT_k/}^{\star \star}
\]
by applying $\Alg_1(-)$. In particular, the right adjoint functor carries $(n\PRTT_{k\otimes_Ak}, (n-1)\PRTT_k)$ to $\DD_{n+2}(A)$ in  $\Alg_{n+1}(\Mod_k)$.
The image of $(n\PRTT_{k\otimes_Ak}, (n-1)\PRTT_k)$ in $\Alg_{1}(n\PRTT_k)$
under the right adjoint functor is $\End_{(n-1)\PRTT_{k\otimes_Ak}}^l((n-1)\PRTT_k)$.
By Proposition~\ref{endoalgebraization} and its proof again, the image of $\End_{(n-1)\PRTT_{k\otimes_Ak}}^l((n-1)\PRTT_k)$ under the right adjoint functor is $\DD_{n+1}(k\otimes_Ak)$.
Thus, $\DD_{n+2}(A)\simeq \DD_{n+1}(k\otimes_Ak)$.
\QED

\begin{Lemma}
\label{categorytensor}
Let $A$ be an $\mathbf{E}_{n+2}$-algebra.  There exists a canonical equivalence $n\PRTT_k\otimes_{n\PRTT_A}n\PRTT_k\simeq n\PRTT_{k\otimes_Ak}$.
\end{Lemma}

\Proof
The proof proceeds by induction on $n$.
When $n=0$, the functor
$c_0:\Mod_k\otimes_{\LMod_A}\Mod_k\to \LMod_{k\otimes_Ak}$ induced by the augmented simplicial diagram $A^{\otimes\bullet}\to k\to k\otimes_Ak$ of the bar construction of $k\otimes_Ak$
is an equivalence (see \cite[4.8.5.13]{HA}).
Next, we assume that  for an $\mathbf{E}_{n+1}$-algebra $A$ the canonical map $c_{n-1}:(n-1)\PRTT_k\otimes_{(n-1)\PRTT_A}(n-1)\PRTT_k\to (n-1)\PRTT_{k\otimes_Ak}$
induced by the bar construction is an equivalence.
We will prove that $c_n:n\PRTT_k\otimes_{n\PRTT_A}n\PRTT_k\to n\PRTT_{k\otimes_Ak}$ is an equivalence.
We note that $\theta_{n\PRTT_k}:\Alg_1(n\PRTT_k)\to ((n+1)\PRTT_k)_{n\PRTT_k/}$ is fully faithful (cf. Construction~\ref{cconstruction}),
and the composition $\Alg_1(n\PRTT_k)\to  ((n+1)\PRTT_k)_{n\PRTT_k/} \to ((n+1)\PRTT_k)$ with the forgetful functor
preserves geometric realizations into the essential images of the composite.
Thus, under the assumption, we are reduced to showing that there exists an equivalence $n\PRTT_k\otimes_{n\PRTT_A}n\PRTT_k\simeq \LMod_{D}(n\PRTT_k)$ for some $D\in \Alg_1(n\PRTT_k)$.
By definition and \cite[7.1.3.1]{HA}, $n\PRTT_k\otimes_{n\PRTT_A}n\PRTT_k$ is equivalent to 
\[
\LMod_{(n-1)\PRTT_k}(\LMod_{(n-1)\PRTT_A}(n\PRTT_k))\otimes_{\LMod_{(n-1)\PRTT_A}(n\PRTT_k)} \RMod_{(n-1)\PRTT_k}(\LMod_{(n-1)\PRTT_A}(n\PRTT_k))
\]
where we regard $(n-1)\PRTT_k$ as an $\eone$-algebra in $\LMod_{(n-1)\PRTT_A}(n\PRTT_k)$.
This is equivalent to 
\[
\mathcal{B}_{n\PRTT_k}=\RMod_{(n-1)\PRTT_k}(\LMod_{(n-1)\PRTT_k}(\LMod_{(n-1)\PRTT_A}(n\PRTT_k)))
\]
(cf. \cite[4.8.4.6]{HA}).
For ease of notation, we write $\RMod_{B}(\LMod_C(\LMod_{A'}(P)))$ for $\mathcal{B}_{n\PRTT_k}$;
$P=n\PRTT_k$, $A'=(n-1)\PRTT_A$, etc.
It suffices to prove that $\RMod_{B}(\LMod_C(\LMod_{A'}(P)))$ is equivalent to $\LMod_D(P)(=\RMod_{D^{op}}(P))$ for some $D\in \Alg_1(P)$.
Consider the sequence of adjunctions
\[
P\stackrel{f}{\rightleftarrows} \LMod_{A'}(P) \stackrel{g}{\rightleftarrows} \RMod_B(\LMod_{A'}(P))  \stackrel{h}{\rightleftarrows} \RMod_{B}(\LMod_C(\LMod_{A'}(P)))
\]
such that $f$, $g$, and $h$ are free functors, and left arrows are forgetful functors.
The functor $f$ is given informally by $A'\otimes (-)$. The functor 
$g$ is given by $B\otimes_{A'}(-)$. 
The functor $h$ is given by $C\otimes_{A'}(-)$.
More precisely, we let $\alpha:\LMod_{A'}(P) \rightarrow\LMod_C(\LMod_{A'}(P))$
be the (unique) right $\LMod_{A'}(P)$-module map which sends $A'$ to $C$.
Applying $\RMod_B$ to $\alpha$ we get a $P$-module map $h$.
We consider the composite of adjunctions
\[
F_P:P\rightleftarrows \mathcal{B}_P=\RMod_{B}(\LMod_C(\LMod_{A'}(P))):G_P.
\]
We let $C\otimes_{A'}B$ be the image of the unit object of $P$ under $F_P$, so that
$F_P$ is given by $C\otimes_{A'}B\otimes (-)$.
The right adjoint $G_P$ is described by using 
Hom objects in $P$; $G_P(M)\simeq \Hom_{\mathcal{B}}(C\otimes_{A'}B,M)$. We write $X:=C\otimes_{A'}B$.
Let $E$ be the endomorphism algebra of $X$ defined as an object of $\Alg_1(P)$.
We write $Q=\Mod_{(n-1)\PRTT_k}(\vCat)$. This symmetric monoidal $\infty$-category is a $\mathbb{V}$-presentably symmetric monoidal $\infty$-category,
and $P\subset Q$ is closed under tensor products and small colimits.
Then $(F,G)$ is extended to
\[
F:Q\rightleftarrows \mathcal{B}=\RMod_{B}(\LMod_C(\LMod_{A'}(Q))):G.
\]
Since $Q$ is a $\mathbb{V}$-presentably symmetric monoidal $\infty$-category, and $\mathcal{B}$ is a $Q$-module ($Q$-tensored),
$\mathcal{B}$ and $Q$ are $Q$-enriched (cf. \cite[4.2.1]{HA}, \cite[7.4.13]{GH}). 
We define $F':\RMod_{E}(Q)\to \mathcal{B}$ by $L\mapsto L\otimes_EX$, which is determined by the universal left action on $E$ on $X$.
We define $G':\mathcal{B}\to \RMod_{E}(Q)$ as a right adjoint of $F'$, which is given by $M\mapsto\Hom_{\mathcal{B}}(X,M)$. 
The canonical evaluation map
$\Hom_{\mathcal{B}}(X,M)\otimes X\to M$
induces $F'(G'(M))=\Hom_{\mathcal{B}}(X,M)\otimes_EX\to M$.
Thus, we have $\epsilon:F'\circ G' \to \textup{id}$.
The canonical map $L\otimes X\to L\otimes_EX$
induces $L\to \Hom_{\mathcal{B}}(X,L\otimes_EX)=G'(F'(L))$.
Thus, we have $\eta:\textup{id}\to G'\circ F'$.
We show that $\eta$ and $\epsilon$
are equivalences.
Note that $G'$ preserves $\mathbb{V}$-small colimits since the composite
$\mathcal{B}\stackrel{G'}{\to} \RMod_{E}(Q)\stackrel{\textup{forget}}{\to} Q$
preserves $\mathbb{V}$-small colimits, and the forgetful functor is conservative and preserves $\mathbb{V}$-small colimits.
In particular, $G'\circ F'$ preserves geometric realizations. We may assume that $L$ is a free right $E$-module $T\otimes E$ for some $T\in Q$. Then
$T\otimes \Hom_{\mathcal{B}}(X,X)\to \Hom_{\mathcal{B}}(X,T\otimes E\otimes_EX)$ is equivalent to an identity $T\otimes X\to T\otimes X$.
This shows that $\eta$ is an equivalence. 
Next, we show that $\epsilon_M:\Hom_{\mathcal{B}}(X,M)\otimes_EX \to M$ is an equivalence.
As above, we may assume that $M=T\otimes X$.
Then we have canonical equivalences $\Hom_{\mathcal{B}}(X,T\otimes X)\otimes_EX\simeq T\otimes \Hom_{\mathcal{B}}(X,X)\otimes_EX\simeq T\otimes X$,
and $\epsilon_M$ can be identified with the identity $T\otimes X\to T\otimes X$. Thus, we see that $\epsilon$ is an equivalence. 
Thus, $F'$ and $G'$ induce mutually inverse equivalences $\RMod_E(Q)\simeq \RMod_{B}(\LMod_C(\LMod_{A'}(Q)))$.
(For the proof of equivalences of $F'$ and $G'$, one may also use Barr-Beck-Lurie theorem by combining the adjunction $(F,G)$ and the free-forget adjunction
$Q\rightleftarrows \RMod_E(Q)$, cf. \cite[4.7.3.16]{HA}.)
By the definition of $G'$ (or $F'$), this induces an equivalence $\RMod_E(P)\simeq \RMod_{B}(\LMod_C(\LMod_{A'}(P)))$.
\QED


\subsection{Categorified Koszul duality of $t$-structured $(\infty,n)$-categories}

\begin{Definition}
\label{defartin}
Let $R$ be an augmented $\eenu$-algebra in $\Mod_k$. Let $n$ be a natural number or $\infty$.
The augmented $\eenu$-algebra $R$ is Artin if the following properties hold:
\begin{itemize}
\item $R$ is connective,

\item $\oplus_{i\in \ZZ}H^{i}(R)$ is finite dimensional over $k$,

\item if $I$ is the radical of $H^0(R)$, the canonical map $k\to H^0(R)/I$ is an isomorphism.

\end{itemize}

If the truncation $\tau^{\ge -m}R$ is Artin for any nonnegative integer $m\ge 0$, we say that
$R$ is almost Artin.

By convention, an $\mathbf{E}_{\infty}$-algebra is a commutative algebra object in $\Mod_k$, and $\Alg_\infty(\Mod_k)=\CAlg(\Mod_k)$.
We let $\textup{Art}_n$ denote the full subcategory on $\Alg_n^+(\Mod_k)$ spanned by Artin $\eenu$-algebras.
We let $\textup{AArt}_n$ denote the full subcategory on $\Alg_n^+(\Mod_k)$ spanned by almost Artin $\eenu$-algebras.
\end{Definition}

\begin{Remark}
In \cite{DAG}, the terminology small is adopted instead of Artin.
\end{Remark}

Let $R$ be an augmented  $\mathbf{E}_{n+1}$-algebra. Note that $(n-1)\PRTT_R$ and $(n-1)\PRTT_{\DD_{n+1}(R)}$
are $\eone$-algebra objects in $n\PRTT_k$.
We regard $(n-1)\PRTT_k$ as a left $(n-1)\PRTT_R\otimes_{(n-1)\PRTT_k} (n-1)\PRTT_{\DD_{n+1}(R)}$-module determined by
$(n-1)\PRTT_R\otimes_{(n-1)\PRTT_k} (n-1)\PRTT_{\DD_{n+1}(R)}\to (n-1)\PRTT_k$
induced by the universal pairing
$R\otimes_k\DD_{n+1}(R)\to k$. We have
\[
\xymatrix{
(-)\otimes_{(n-1)\PRTT_R}(n-1)\PRTT_k :\RMod_{(n-1)\PRTT_R}(n\PRTT_k) \ar[r] &  \LMod_{\PRTT_{\DD_{n+1}(R)}}(n\PRTT_k). 
}
\]
We write $\Phi_{R,n}$ for $(-)\otimes_{(n-1)\PRTT_R}(n-1)\PRTT_k$.
Similarly, we set 
\[
\Psi_{R,n}=(n-1)\PRTT_k\otimes_{(n-1)\PRTT_{\DD_{n+1}(R)}}(-):\LMod_{(n-1)\PRTT_{\DD_{n+1}(R)}}(n\PRTT_k)  \to \RMod_{(n-1)\PRTT_{R}}(n\PRTT_k).
\]

\begin{Theorem}
\label{HCKD}
Let  $R$ be an almost Artin $\mathbf{E}_{n+1}$-algebra.
Suppose that $n\ge1$. Then $\Psi_{R,n}$ and $\Phi_{R,n}$ induce mutually inverse equivalences 
\[
\xymatrix{
\widehat{\Psi}_{R,n}:\LMod_{(n-1)\PRTTT_{\DD_{n+1}(R)}}(n\PRTTT_k) \ar[r]^(0.55){\simeq} &  \RMod_{(n-1)\PRTTT_{R}}(n\PRTTT_k) : \widehat{\Phi}_{R,n} \ar[l]
}
\]
in $(n+1)\PRTTT_k$. The functor $\widehat{\Psi}_{R,n}$ sends $\mathcal{E}$ to $(n-1)\PRTTT_k\widehat{\otimes}_{(n-1)\PRTTT_{\DD_{n+1}(R)}}\EE$ where $\widehat{\otimes}_{(n-1)\PRTTT_{\DD_{n+1}(R)}}$
means the relative tensor product in $n\PRTTT_k$.  
The functor $\widehat{\Phi}_{R,n}$ sends $\mathcal{D}$ to $\mathcal{D}\widehat{\otimes}_{(n-1)\PRTTT_{R}}(n-1)\PRTTT_k$ where $\widehat{\otimes}_{(n-1)\PRTTT_{R}}$
means the relative tensor product in $n\PRTTT_k$.  
\end{Theorem}

We deduce this theorem from Proposition~\ref{caseone} and Proposition~\ref{moritaleft}.

\begin{Proposition}
\label{caseone}
When $n=1$, the completions of $\Phi_{R,1}$ and $\Psi_{R,1}$ induce mutually inverse equivalences
\[
\widehat{\Phi}_{R,1} :\RMod_{\LMod_R}(\PRTTT_k) \stackrel{\simeq}{\to}  \LMod_{\LMod_{\DD_{2}(R)}}(\PRTTT_k),\ \ \ \widehat{\Psi}_{R,1}:\LMod_{\LMod_{\DD_{2}(R)}}(\PRTTT)  \stackrel{\simeq}{\to} \RMod_{\LMod_{R}}(\PRTTT_k).
\]
Here $\LMod_R$ and $\LMod_{\DD_{2}(R)}$ are equipped with the standard $t$-structures (see Example~\ref{connectivetstr} and Remark~\ref{unifiedtstr}). 
\end{Proposition}

When $R$ is Artin, this is one of the main results of \cite{IH}. Proposition~\ref{caseone} generalizes the result of \cite{IH} to the case of almost Artin algebras.
We will prove Proposition~\ref{caseone} in the next section.

\begin{Proposition}
\label{moritaleft}
We assume that
for an almost Artin $\mathbf{E}_{n+1}$-algebra $R$,
\[
\Phi_{R,n} :\RMod_{(n-1)\PRTT_R}(n\PRTT_k) \to  \LMod_{(n-1)\PRTT_{\DD_{n+1}(R)}}(n\PRTT_k)
\]
and 
\[
\Psi_{R,n}:\LMod_{(n-1)\PRTT_{\DD_{n+1}(R)}}(n\PRTT)  \to \RMod_{(n-1)\PRTT_{R}}(n\PRTT)
\]
are mutually inverse equivalences after the completion.
Then for an almost Artin $\mathbf{E}_{n+2}$-algebra $R$,
the functors $\Psi_{R,n+1}:\LMod_{n\PRTT_{\DD_{n+2}(R)}}((n+1)\PRTT)  \to \RMod_{n\PRTT_{R}}((n+1)\PRTT)$ and
$\Phi_{R,n+1}:\RMod_{n\PRTT_R}((n+1)\PRTT_k) \to  \LMod_{n\PRTT_{\DD_{n+2}(R)}}((n+1)\PRTT_k)$
induce mutually inverse equivalences after the completion.
\end{Proposition}

The rest of this section is  devoted to the proof of Proposition~\ref{moritaleft}.

\begin{Construction}
\label{xidefine}
Let $R$ be an $\mathbf{E}_{n+2}$-algebra.
We write $D=\DD_{n+2}(R)$.
We define a right  $n\PRTT_D$-module map
\[
\xi:n\PRTT_{D}\to n\PRTT_k\otimes_{(n\PRTT_R)^{\otimes opp}}n\PRTT_k.
\]
The superscript $\otimes opp$ in $(n\PRTT_R)^{\otimes opp}$ indicates the opposite algebra of $n\PRTT_R$ in $\Alg_1((n+1)\PRTT_k)$.
This $\xi$ is
induced by the algebraization functor $n\PRTT_{D}\to \End_{n\PRTT_R}^l(n\PRTT_k)$ (see Definition~\ref{algebraizationdef} and Proposition~\ref{endoalgebraization})
and 
\[
m: \End_{n\PRTT_R}^l(n\PRTT_k)^{\otimes opp}  \simeq \End_{n\PRTT_R}^r(n\PRTT_k)\to n\PRTT_k\otimes_{(n\PRTT_R)^{\otimes opp}}n\PRTT_k
\]
defined as follows: The tautological left module action of $\End_{n\PRTT_R}^r(n\PRTT_k)$ on $n\PRTT_k$ on the right-hand side in $n\PRTT_k\otimes_{(n\PRTT_R)^{\otimes opp}}n\PRTT_k$
induces
\[
\End_{n\PRTT_R}^r(n\PRTT_k)\times n\PRTT_k\otimes_{(n\PRTT_R)^{\otimes opp}}n\PRTT_k\to n\PRTT_k\otimes_{(n\PRTT_R)^{\otimes opp}}n\PRTT_k
\]
whose evaluation at $(n-1)\PRTT_k\otimes (n-1)\PRTT_k\in n\PRTT_k\otimes_{(n\PRTT_R)^{\otimes opp}}n\PRTT_k$ determines
\[
m:\End_{n\PRTT_R}^r(n\PRTT_k)  = \End_{n\PRTT_R}^r(n\PRTT_k)\times \{ (n-1)\PRTT_k\otimes (n-1)\PRTT_k  \}\to n\PRTT_k\otimes_{(n\PRTT_R)^{\otimes opp}}n\PRTT_k,
\]
where $(n-1)\PRTT_k\otimes (n-1)\PRTT_k$ denotes the image of $ (n-1)\PRTT_k\otimes (n-1)\PRTT_k\in n\PRTT_k\otimes_{n\PRTT_k}n\PRTT_k$
under the functor $ n\PRTT_k\otimes_{n\PRTT_k}n\PRTT_k\to  n\PRTT_k\otimes_{(n\PRTT_R)^{\otimes opp}}n\PRTT_k$ induced by $n\PRTT_k\to n\PRTT_R$.
The functor sends the unit object of $\End_{n\PRTT_R}^r(n\PRTT_k)$ (the identity functor)
to $(n-1)\PRTT_k\otimes (n-1)\PRTT_k \in n\PRTT_k\otimes_{(n\PRTT_R)^{\otimes opp}}n\PRTT_k$.
The algebraization functor and $m$ determine a right module action of $n\PRTT_D$ on $n\PRTT_k\otimes_{(n\PRTT_R)^{\otimes opp}}n\PRTT_k$.
Then
$\xi$ is defined to be a (unique) right $n\PRTT_D$-linear functor $n\PRTT_D\to n\PRTT_k\otimes_{(n\PRTT_R)^{\otimes opp}}n\PRTT_k$, which sends $(n-1)\PRTT_D$ to $(n-1)\PRTT_k\otimes(n-1)\PRTT_k$.
\end{Construction}

\begin{Construction}
\label{etadefine}
Let $R$ be an $\mathbf{E}_{n+2}$-algebra.
We write $D=\DD_{n+2}(R)$.
We define a left $n\PRTT_R$-module map
\[
\eta:n\PRTT_k\otimes_{(n\PRTT_{D})^{\otimes opp}}n\PRTT_k\simeq  n\PRTT_k\otimes_{(n\PRTT_{D})}n\PRTT_k   \to \RMod_{(n-1)\PRTT_R}(n\PRTT_k).
\]
This $\eta$ is defined as follows.
If we look at $n\PRTT_k$ on the right-hand side in 
$n\PRTT_k\otimes_{(n\PRTT_{D})^{\otimes opp}}n\PRTT_k$,
this has the $(n\PRTT_{D})^{\otimes opp}$-$n\PRTT_R$-bimodule structure arising from the universal pairing $D\otimes_k R\to k$.
In particular, $n\PRTT_k\otimes_{(n\PRTT_{D})^{\otimes opp}}n\PRTT_k$ inherits a right $n\PRTT_R$-module structure.
We have a (unique) right $n\PRTT_R$-linear functor 
\[
\eta':n\PRTT_R \to n\PRTT_k\otimes_{(n\PRTT_D)^{\otimes opp}}n\PRTT_k\simeq\RMod_{(n-1)\PRTT_{k\otimes_Dk}}(n\PRTT_k),
\]
which sends $(n-1)\PRTT_R$ to $(n-1)\PRTT_k\in\RMod_{(n-1)\PRTT_{k\otimes_Dk}}(n\PRTT_k)$ (the $(n-1)\PRTT_{k\otimes_Dk}$-module structure on $(n-1)\PRTT_k$ is given by
the augmentation $(n-1)\PRTT_{k\otimes_Dk}\to (n-1)\PRTT_{k}$).
Taking duals $\Fun_{n\PRTT_k}(-,n\PRTT_k)$ 
we have
\[
\eta:\LMod_{(n-1)\PRTT_{k\otimes_{D}k}}(n\PRTT) \simeq n\PRTT_k\otimes_{n\PRTT_{D}}n\PRTT_k \stackrel{p}{\to} \RMod_{(n-1)\PRTT_R}(n\PRTT) 
\]
which is left $n\PRTT_R$-linear. See Lemma~\ref{categorytensor} for the equivalence.

\end{Construction}

The relevant feature of $\xi$ and $\eta$ is that $\xi$ is a right  $n\PRTT_D$-module map,
and $\eta$ is a left $n\PRTT_R$-module map.

\begin{Lemma}
\label{indxi}
Let
$\xi:n\PRTT_{\DD_{n+2}(R)}\to n\PRTT_k\otimes_{(n\PRTT_R)^{\otimes opp}}n\PRTT_k \simeq\RMod_{(n-1)\PRTT_{k\otimes_Rk}}(n\PRTT_k)$
be the morphism in Construction~\ref{xidefine}. See Lemma~\ref{categorytensor} for the latter equivalence.
We consider 
\[
\Psi_{k\otimes_Rk,n}:n\PRTT_{\DD_{n+1}(k\otimes_Rk)}=\LMod_{(n-1)\PRTT_{\DD_{n+1}(k\otimes_Rk)}}(n\PRTT_k)  \to \RMod_{(n-1)\PRTT_{k\otimes_Rk}}(n\PRTT_k)
\]
given informally by $(n-1)\PRTT_k\otimes_{(n-1)\PRTT_{\DD_{n+1}(k\otimes_Rk)}}(-)$.
Here $(n-1)\PRTT_k$
is regarded as a $((n-1)\PRTT_{\DD_{n+1}(k\otimes_Rk)})^{\otimes opp}$-$(n-1)(\PRTT_{k\otimes_Rk})^{\otimes}$-bimodule determined by
 $(n-1)\PRTT_{k\otimes_Rk}\otimes_{(n-1)\PRTT_k}(n-1)\PRTT_{\DD_{n+1}(k\otimes_Rk)}\to (n-1)\PRTT_k$ arising from the universal pairing $(k\otimes_Rk)\otimes_k\DD_{n+1}(k\otimes_Rk)\to k$.
Then $\xi$ is equivalent to $\Psi_{k\otimes_Rk,n}$ as a morphism in $(n+1)\PRTT_k$.
\end{Lemma}

\Proof
Write $D=\DD_{n+2}(R)$ and $\Psi=\Psi_{k\otimes_Rk,n}$.
We first note that $\xi$ carries $(n-1)\PRTT_D$ to $(n-1)\PRTT_k\otimes(n-1)\PRTT_k $.
By definition $\Psi$ also carries $(n-1)\PRTT_D$ to $(n-1)\PRTT_k\otimes (n-1)\PRTT_k$.
To prove the assertion, we look at the endomorphism algebras.  
Recall the sequence of adjoint pairs from Lemma~\ref{goodalgebraization}:
\[
\xymatrix{
\Alg_{n+1}(\Mod_k)   \ar@<0.5ex>[r]^{\Alg_{n}(\theta_{0\PRTT_k,0})}  & \Alg_{n}(\PRTT_k)  \ar@<0.5ex>[l]   \ar@<0.5ex>[r]^{U}  &  ((n+1)\PRTT_k)_{n\PRTT_k/}^{\star \star}    \ar@<0.5ex>[l]^{V}  .
}
\]
See Definition~\ref{thetadeftt} for $\theta_{0\PRTT_k,0}$.
We regard
$\RMod_{(n-1)\PRTT_{k\otimes_Rk}}(n\PRTT_k)$ with the specified object $(n-1)\PRTT_k$ as an object of $((n+1)\PRTT_k)_{n\PRTT_k/}^{\star \star}$ (see Lemma~\ref{goodendoendo}).
Let $B$ be the image of $P=(\RMod_{(n-1)\PRTT_{k\otimes_Rk}}(n\PRTT_k),(n-1)\PRTT_k)$ in $\Alg_{n+1}(\Mod_k)$ under the composite of right adjoint functors.
By the adjoint pairs, we have a canonical equivalence of mapping spaces 
\[
\alpha:\Map_{\Alg_{n+1}(\Mod_k)}(C,B)\simeq \Map_{(n+1)\PRTT_k)_{n\PRTT_k/}}((n\PRTT_C,(n-1)\PRTT_C),P).
\]
By Corollary~\ref{barkoszuldual}, 
there exists an equivalence $B\simeq \DD_{n+1}(k\otimes_Rk)\simeq \DD_{n+2}(R)=D$ in $\Alg_{n+1}(\Mod_k)$.

The universal map
$(n\PRTT_D,(n-1)\PRTT_D)\to P$ is induced by the universal morphism
\[
u:((n-1)\PRTT_D)^{\otimes opp} \to \End_{P}((n-1)\PRTT_k)
\]
in $\Alg_{1}(n\PRTT_k)$,
which induces an equivalence
\[
\Map_{\Alg_{n+1}(\Mod_k)}(C,D)\simeq \Map_{\Alg_1(n\PRTT)}(((n-1)\PRTT_C)^{\otimes opp}, \End_{P}((n-1)\PRTT_k))
\]
given by $f\mapsto u\circ (n-1)\PRTT_f$.
As explained in Remark~\ref{pairingalgebraization}, $u$ corresponds to the pairing  $(n-1)\PRTT_{k\otimes_Rk}\otimes_{(n-1)\PRTT_k}(n-1)\PRTT_{D}\to (n-1)\PRTT_k$ induced by
the universal pairing $(k\otimes_Rk)\otimes_k\DD_{n+1}(k\otimes_Rk)\to k$.
Thus, $u\simeq \Psi$ since $u$ is the functor defined by the integral kernel
$(n-1)\PRTT_k$ whose $((n-1)\PRTT_D)^{\otimes opp}$-$((n-1)\PRTT_{k\otimes_Rk})^{\otimes}$-bimodule structure is given by
 $(n-1)\PRTT_{k\otimes_Rk}\otimes_{(n-1)\PRTT_k}(n-1)\PRTT_{D}\to (n-1)\PRTT_k$.

We note that $\xi$ is the composite of 
the algebraization functor $n\PRTT_D\to \End_{n\PRTT_R}^l(n\PRTT_k)$
and the equivalence $m:\End_{n\PRTT_R}^l(n\PRTT_k)\simeq n\PRTT_k\otimes_{n\PRTT_R}n\PRTT_k$
in $(n+1)\PRTT_k$  (cf. Construction~\ref{xidefine}).
The equivalence $m$ sends the identity functor to $(n-1)\PRTT_k$.
From the construction of the algebraization functor, if we regard $\xi$ as a morphism 
$(n\PRTT_D,(n-1)\PRTT_D)\to P$ in $((n+1)\PRTT_k)_{n\PRTT_k/}$,
then it corresponds to the identity map $D\to D$ in $\Alg_{n+1}(\Mod_k)$ through $\alpha$.
Consequently, $\xi\simeq \Psi$.
\QED

\begin{Lemma}
\label{indeta}
Let
$\eta:\LMod_{(n-1)\PRTT_{k\otimes_{D}k}}(n\PRTT) \to \RMod_{(n-1)\PRTT_R}(n\PRTT) $
be the morphism in Construction~\ref{etadefine}.
We consider 
\[
\Psi_{R,n}:\LMod_{(n-1)\PRTT_{(k\otimes_Dk)}}(n\PRTT_k)  \to \RMod_{(n-1)\PRTT_{R}}(n\PRTT_k)
\]
given informally by $(-)\otimes_{(n-1)\PRTT_{(k\otimes_Dk)}}(n-1)\PRTT_k$.
We note that $k\otimes_{D}k= k\otimes_{\DD_{n+2}(R)}k\simeq \DD_{n+1}(R)$ (see Proposition~\ref{dbar}).
Here $(n-1)\PRTT_k$
is regarded as a $(n-1)\PRTT_{(k\otimes_Dk)}$-$(n-1)(\PRTT_{R})^{\otimes opp}$-bimodule determined by
 $(n-1)\PRTT_{R}\otimes_{(n-1)\PRTT_k}(n-1)\PRTT_{\DD_{n+1}(R)}\to (n-1)\PRTT_k$ arising from the universal pairing $R\otimes_k\DD_{n+1}(R)\to k$.
Then $\eta$ is equivalent to $\Psi_{R,n}$ as a morphism in $(n+1)\PRTT_k$.
\end{Lemma}

\Proof
Let $\eta':n\PRTT_R\to \RMod_{(n-1)\PRTT_{k\otimes_Dk}}(n\PRTT_k)$
be a (unique) right $n\PRTT_R$-linear functor
which sends $(n-1)\PRTT_R$ to $(n-1)\PRTT_k$.
The morphism $\eta$ is defined to be its $n\PRTT_k$-linear dual.
By Lemma~\ref{indxi}, $\eta'$
is determined by $n\PRTT_k\otimes_{n\PRTT_R}(-)$,
where the bimodule structure on $n\PRTT_k$
is induced by the universal $\mathbf{E}_{n+1}$-Koszul dual pairing 
\[
 R\otimes_k\DD_{n+1}(R)\simeq \DD_{n+1}(\DD_{n+1}(R))\otimes \DD_{n+1}(R)\to k.
\]
($R$ is almost Artin and $R\simeq \DD_{n+1}(\DD_{n+1}(R))$. See Proposition~\ref{almostbiduality})
Now our assertion follows from a general observation of Lemma~\ref{dualkernel}.
\QED

\begin{Lemma}
\label{dualkernel}
Let $\CCC^\otimes$ be an object of $\CAlg(\vCat)$, that is,
a symmetric monoidal $\infty$-category such that $\CCC$ has small colimits and the 
tensor product preserves small colimits separately in each variable.
Let $A$ and $B$ be $\eone$-algebra objects in $\CCC^\otimes$.
Let $\LMod_A(\CCC)$ and $\LMod_B(\CCC)$ be the $\infty$-categories 
of left $A$-modules and left $B$-modules, respectively.
Let $F:\LMod_A(\CCC)\to \LMod_B(\CCC)$ be the functor determined by
a $B$-$A$-bimodule ${}_BM_A$. Namely, $F$ is given by $P\mapsto {}_BM_A\otimes_AP$.
Let $\Fun_{\CCC}(\LMod_A(\CCC),\CCC)$ and $\Fun_{\CCC}(\LMod_B(\CCC),\CCC)$
denote the $\CCC$-linear functor categories.
Recall the equivalences 
$\RMod_A(\CCC)\simeq \Fun_{\CCC}(\LMod_A(\CCC),\CCC)$ and $\RMod_B(\CCC)\simeq \Fun_{\CCC}(\LMod_B(\CCC),\CCC)$.
We consider 
\begin{eqnarray*}
\RMod_B(\CCC)   &\simeq& \Fun_{\CCC}(\LMod_B(\CCC),\CCC)  \\
&\to&  \Fun_{\CCC}(\LMod_A(\CCC),\CCC) \\
&\simeq& \RMod_A(\CCC).
\end{eqnarray*}
Then the composite $F^\vee:\RMod_B(\CCC)\to \RMod_A(\CCC)$ is
determined by the $B$-$A$-bimodule ${}_BM_A$. Namely, $F^\vee$ is given by $Q\mapsto Q\otimes_B{}_BM_A$.
\end{Lemma}

\Proof
For $N_B\in\RMod_B(\CCC)$, the corresponding functor $\LMod_B(\CCC)\to \CCC$ is given informally by $N_B\otimes_B(-)$.
Its composition with $F$ is $N_B\otimes_B{}_BM_A\otimes_A(-):\LMod_A(\CCC)\to \CCC$.
Thus, $F^\vee$ is given informally by $N_B\mapsto N_B\otimes_B{}_BM_A$.
Namely, $F^\vee$ is determined by  $(-)\otimes_B{}_BM_A$. This is our assertion.
\QED

{\it Proof of Proposition~\ref{moritaleft}.}
We write $D=\DD_{n+2}(R)$.
For $\EE\in \LMod_{n\PRTT_{\DD_{n+2}(R)}}((n+1)\PRTT)$ we define
\[
U_{\EE}:\EE\to (n\PRTT_k\otimes_{n\PRTT_{D}}\EE)\otimes_{n\PRTT_R}n\PRTT_k \simeq (n\PRTT_k\otimes_{(n\PRTT_R)^{\otimes opp}}n\PRTT_k)\otimes_{n\PRTT_{D}}\EE
\]
to be the morphism obtained from $\xi:n\PRTT_{D}\to n\PRTT_k\otimes_{(n\PRTT_R)^{\otimes opp}}n\PRTT_k$
by taking the tensor product over $n\PRTT_{D}$.
For $\DDD\in \RMod_{n\PRTT_{R}}((n+1)\PRTT)$ we define
\[
V_{\DDD}:n\PRTT_k\otimes_{n\PRTT_{D}}(\DDD\otimes_{n\PRTT_{R}}n\PRTT_k) \simeq  \DDD\otimes_{n\PRTT_{R}}(n\PRTT_k\otimes_{(n\PRTT_{D})^{\otimes opp}}n\PRTT_k)\to \DDD
\]
to be the morphism obtained from $\eta:n\PRTT_k\otimes_{(n\PRTT_{D})^{\otimes opp}}n\PRTT_k\simeq n\PRTT_k\otimes_{\RMod_{(n-1)\PRTT_{D}}(n\PRTT)}n\PRTT_k\to n\PRTT_R$
by taking the tensor product over $n\PRTT_{R}$.

We set $\Phi=\Phi_{R,n+1}$ and $\Psi=\Psi_{R,n+1}$.
They determine the natural transformations $U:\textup{id}\to \Phi\circ \Psi$ and $V:\Psi \circ \Phi\to \textup{id}$.
We will prove that both $U$ and $V$ are equivalences after the completion.
Suppose that $\EE\in \LMod_{n\PRTTT_{\DD_{n+2}(R)}}((n+1)\PRTTT)$.
The completion of $U_{\EE}$ determines
\[
\widehat{U}_{\EE}:\EE=n\PRTTT_D\widehat{\otimes}_{n\PRTTT_{D}} \EE\to  (n\PRTTT_k\widehat{\otimes}_{(n\PRTTT_R)^{\otimes opp}}n\PRTTT_k)\widehat{\otimes}_{n\PRTTT_{D}}\EE.
\]
According to Lemma~\ref{indxi} and the inductive assumption,
the completion of $\xi$ induces an equivalence $n\PRTTT_{D}\simeq n\PRTTT_k\widehat{\otimes}_{(n\PRTTT_R)^{opp}}n\PRTTT_k$.
Thus, $\widehat{U}_{\EE}$ is an equivalence.
Next, we
suppose that $\DDD\in \RMod_{n\PRTT_{R}}((n+1)\PRTTT)$.
The completion of $V_{\DDD}$ determines
\[
\widehat{V}_{\DDD}:\DDD\widehat{\otimes}_{n\PRTTT_{R}}(n\PRTTT_k\widehat{\otimes}_{(n\PRTTT_{D})^{\otimes opp}}n\PRTTT_k)\to \DDD\widehat{\otimes}_{n\PRTTT_{R}}\RMod_{(n-1)\PRTTT_R}(n\PRTTT_k)\simeq \DDD.
\]
For the right equivalence we note that $(n-1)\PRTTT_R$ is $\etwo$-monoidal so that
$\RMod_{(n-1)\PRTTT_R}(n\PRTTT_k)$ is the opposite algebra of $\LMod_{(n-1)\PRTTT_R}(n\PRTTT_k)=n\PRTTT_R$ in $\Alg_1((n+1)\PRTTT_k)$ (see \cite[Proposition 2.15]{IH}).
According to Lemma~\ref{indeta} and the inductive assumption,
the completion of $\eta$ induces an equivalence $n\PRTTT_k\widehat{\otimes}_{(n\PRTTT_{D})^{\otimes opp}}n\PRTTT_k\simeq n\PRTTT_k\widehat{\otimes}_{n\PRTTT_{D}}n\PRTTT_k\simeq n\PRTTT_R$ so that $\widehat{V}_{\DDD}$ is an equivalence.
\QED

\section{the case of $n=1$}
\label{n1sec}

We treat the case of $n=1$.
Proposition~\ref{caseone} says that there exist mutually inverse equivalences
\[
\xymatrix{
\LMod_{\LMod_{\DD_{2}(R)}}(\PRTTT_k) \ar@<0.5ex>[r]  &  \RMod_{\LMod_R}(\PRTTT_k)   \ar@<0.5ex>[l]  
}
\]
which are induced by $\Phi_{R,1}$ (the left arrow) and $\Psi_{R,1}$ (the right arrow) 
together with the left completion functor.
Precisely, $\LMod_{\LMod_{\DD_{2}(R)}}(\PRTTT_k) \to  \RMod_{\LMod_R}(\PRTTT_k)$
is defined by $\widehat{L}(\Mod_k\otimes_{\LMod_{\DD_2(R)}}(-))$,
and $\RMod_{\LMod_R}(\PRTTT_k)\to \LMod_{\LMod_{\DD_{2}(R)}}(\PRTTT_k)$
is given by $(-)\otimes_{\LMod_{R}}\Mod_k$.
When $R$ is Artin, this is one of the main results of \cite{IH}. 
We give a slight generalization to the case when $R$ is almost Artin.
The argument of \cite[Section 3]{IH} goes through without essential
change, with the following exceptions:

\begin{Proposition}
\label{dbar}
If $R$ is an almost Artin $\etwo$-algebra, there
exists an equivalence $k\otimes_{\DD_2(R)}k\simeq \DD_1(R)$ in $\Alg_1(\Mod_k)$.
\end{Proposition}

\begin{Proposition}
\label{leftcomp2}
Suppose that $R$ is an almost Artin $\etwo$-algebra.
There exists a left completion functor $\eta:\Mod_k\otimes_{\RMod_{\DD_2(R)}}\Mod_k=\RMod_{k\otimes_{\DD_2(R)}k}\to \LMod_R$.
Moreover, the target $\LMod_R$ is a free $\LMod_R$-module of rank one.
\end{Proposition}

Proposition~\ref{dbar}
follows from the proof in \cite[Lemma 3.13]{IH}
if we prove the following:

\begin{Proposition}
\label{almostbiduality}
Let $R$ be an almost Artin $\etwo$-algebra.
Then the biduality morphisms $R\to \DD_1\DD_1(R)$ and $R\to \DD_2\DD_2(R)$ are equivalences.
Moreover, $\DD_1(R)$ and $\DD_2(R)$ are locally finite. Here regarding $R$ as an $\eone$-algebra we define $\DD_1(R)$. 
\end{Proposition}

\begin{Lemma}
\label{onestepbarisomo}
Let $l$ be a natural number $l\ge1$.
Let $f:A'\to A$ be a morphism of Artin $\mathbf{E}_l$-algebras.
Let $n$ be a nonnegative integer.
Assume that $f$ exhibits $A$ as the truncation $\tau^{\ge -n}(A')$. 
Suppose that $H^{i}(A')=0$ for any $i<-n-1$.
Then the induced morphism $H^i(k\otimes_{A'}k)\to H^i(k\otimes_Ak)$ is an isomorphism of finite dimensional vector spaces for any $-n\le i\le 0$.
\end{Lemma}

\Proof
Under our assumption, $A'$ is a square zero extension of $A$ by $H^{-n-1}(A')[n+1]$ (cf. \cite[7.4.1.28]{HA})
so that $f$ is a pullback of the canonical morphism $A\to A\oplus H^{-n-1}(A')[n+2]$.
Thus, there exists a sequence
\[
A'=A(0)\to A(1)\to \cdots \to A(m)=A
\]
in $\textup{Art}_l$
such that,
for each $0\le i\le m-1$,
there exists a pullback square
\[
\xymatrix{
A(i) \ar[r] \ar[d] & A(i+1) \ar[d] \\
k \ar[r] &  k\oplus k[n+2]
}
\]
in $\textup{Art}_l$.
We may and will assume that 
$A'$ is $k\times_{k\oplus k[n+2]} A$.
Then we have the fiber sequence
\[
k[n+1]\to A'\to A
\]
in $\Mod_{A'}$.
We note that $k\otimes_{A'}k\simeq (k\otimes_AA)\otimes_{A'}k\simeq k\otimes_A(A\otimes_{A'}k)$.
There exists a (co)fiber sequence
\[
k\otimes_{A'}k[n+1] \to k=A'\otimes_{A'}k\to A\otimes_{A'}k.
\]
We note that
$H^i(k\otimes_{A'}k[n+1])=0$ for $i\ge -n$.
Taking into account the long exact sequence of the (co)fiber sequence,
we conclude that 
$H^i(k)\to H^i(A\otimes_{A'}k)$ is an isomorphism for any $-n\le i\le 0$.
This means that the amplitude of $H^i(A\otimes_{A'}k)$ is concentrated in cohomological degrees $i<-n-1$ except $H^0(A\otimes_{A'}k)=k$. 
The fiber of the canonical morphism $A\otimes_{A'}k\to k$ of $A$-modules is concentrated in the cohomological degree 
 $i<-n-1$. Thus, the fiber of $k\otimes_{A'}k\to k\otimes_Ak$
  is concentrated in cohomological degrees $i<-n-1$. 
The induced morphism $H^i(k\otimes_{A'}k)\to H^i(k\otimes_Ak)$ is an isomorphism of vector spaces for any $-n\le i\le 0$.
It remains to show that $H^i(k\otimes_{A'}k)$ and $H^i(k\otimes_Ak)$ are finite dimensional.
This is a consequence of the observation that for an augmented connective locally finite $\eone$-algebra $R\to k$,
$k\otimes_Rk$ is also connective and locally finite (a direct computation or the Tor spectral sequence).
Thus, our assertion follows.
\QED

\begin{Corollary}
\label{onestepcor}
Let $f:A'\to A$ be a morphism of Artin $\mathbf{E}_l$-algebras.
Let $n$ be a nonnegative integer.
Assume that $f$ exhibits $A$ as $\tau^{\ge -n}(A')$. 
Suppose that $H^{i}(A')=0$ for any $i<-n-1$.
Then the induced morphism $H^{j}(\DD_l(A))\to H^{j}(\DD_l(A'))$  is an isomorphism of finite dimensional vector spaces for any $0\le j\le n$.
\end{Corollary}

\Proof
For an augmented algebra $A\to k$, we let $\textup{Bar}(A)$ denote the bar construction $k\otimes_Ak$.
This construction induces a functor $\textup{Bar}:\Alg_l(\Mod_k)_{/k}\to \Alg_{l-1}(\Mod_k)_{/k}$. 
The $l$-times iteration yields $\textup{Bar}^l:\Alg_l(\Mod_k)_{/k}\to \Mod_k$.
Then $\DD_l(A)$ can be identified with the $k$-linear dual of $\textup{Bar}^l(A)$ in $\Mod_k$ (see e.g. Corollary~\ref{barkoszuldual}).
Thus, it is enough to show that $H^{j}(\textup{Bar}^l(A'))\to H^{j}(\textup{Bar}^l(A))$
is an isomorphism of finite dimensional vector spaces for any $-n\le j\le 0$.
This follows from Lemma~\ref{onestepbarisomo}.
\QED

{\it Proof of Proposition~\ref{almostbiduality}.}
Let $R$ be either an almost Artin $\eone$-algebra or an almost Artin $\etwo$-algebra.
When $R$ is an almost Artin $\etwo$-algebra,
the argument in the $\eone$-case works, so we focus on the case when $R$ is an almost Artin $\eone$-algebra. 
We write $R_n$ for $\tau^{\ge -n}R$ for $n\ge0$.
We regard 
$R$ as the sequential limit $\lim_{n\ge0 }R_n$ of Artin $\eone$-algebras.
Then the Koszul duality functor carries small colimits to limits so that
\[
\DD_1(\colim_{n\ge0}\DD_1(R_n))\simeq \lim_{n\ge0}\DD_1(\DD_1(R_n))\simeq \lim_{n\ge0}R_n\simeq R
\]
where the right equivalence is the biduality equivalence for Artin $\eone$-algebras (cf. \cite[X, Section 3]{DAG}).
From Corollary~\ref{onestepcor} and  the fact that $\DD_1(R_n)$ is coconnective and locally finite (cf. \cite[X Lemma 4.5.9]{DAG}),
we see that $\colim_{n\ge0}\DD_1(R_n)$ is coconnective and locally finite.
Also, $\DD_1(R)$ is coconnective and locally finite (as it is the linear dual of $k\otimes_Rk$).
Let $\Alg_1^+(\Mod_k)^{\textup{cclf}}$ be the full subcategory
of $\Alg_1^+(\Mod_k)$ spanned by coconnective and locally finite algebras.
According to \cite[X, Corollary 3.1.15]{DAG}, 
(the restriction of) the Koszul duality functor $\DD_1:\Alg_1^+(\Mod_k)^{\textup{cclf}}\to \Alg_1^+(\Mod_k)^{op}$
is fully faithful. 
Since $\DD_1(\colim_{n\ge0}\DD_1(R_n))\simeq R$,
it follows that the canonical morphism $\colim_{n\ge0}\DD_1(R_n)\to \DD_1(R)$ is an equivalence.
Consequently, we have $R\simeq \DD_1(\DD_1(R))$.
\QED

From the above proof, we see:

\begin{Corollary}
\label{limitkoszul}
Let $R$ be an almost Artin $\etwo$-algebra.
We write $R_n$ for $\tau^{\ge -n}R$ for $n\ge0$.
Then the canonical morphism $\colim_{n\ge0}\DD_2(R_n)\to \DD_2(R)$ is an equivalence.

Similarly, we let $R$ be an almost Artin $\eone$-algebra.
We write $R_n$ for $\tau^{\ge -n}R$ for $n\ge0$.
Then the canonical morphism $\colim_{n\ge0}\DD_1(R_n)\to \DD_1(R)$ is an equivalence.
\end{Corollary}

For a morphism $u:A\to B$ of almost Artin $\eone$-algebras, the pushforward functor along
$u$ is defined to be the $k$-linear (colimit-preserving) functor $\Ind(\textup{LCoh}(B))\to \Ind(\textup{LCoh}(A))$
that extends the restriction functor $\textup{LCoh}(B)\to \textup{LCoh}(A)$. See Example~\ref{Indtstr} for $\textup{LCoh}(A)$.
The $!$-pullback functor along $u$ is defined as the right adjoint $\Ind(\textup{LCoh}(A))\to  \Ind(\textup{LCoh}(B))$
of the pushforward functor (see \cite[Appendix]{CKDA} for details).
For an almost Artin $\eone$-algebra $A$, the dualizing object $\omega_A\in \Ind(\textup{LCoh}(A))$ of $A$ is defined to be the $!$-pullback of $k\in \Ind(\textup{Coh}(k))=\Mod_k$ along $k\to A$.

Suppose that $R$ is an almost Artin $\etwo$-algebra.
We write $R_n$ for $\tau^{\ge -n}R$ for $n\ge0$.
Recall $\Ind(\textup{LCoh}(R))$ and $\Ind(\textup{LCoh}(R_n))$ from Example~\ref{Indtstr}.
There is a fully faithful monoidal functor
\[
\Upsilon_{R_n}:\LMod_{R_n}\to \Ind(\textup{LCoh}(R_n)),
\]
see \cite[Corollary 9.6.3]{Ind}, \cite{CKDA}.
In particular, $\Upsilon_{R_n}$ sends $R_n$ to the dualizing object $\omega_{R_n}\in \textup{LCoh}(R_n)$.
The functors $\Upsilon_{R_n}$ and $\Upsilon_{R_{n+1}}$ are compatible with the (usual derived) pullback functor $\LMod_{R_{n+1}}\to \LMod_{R_n}$
and the $!$-pullback functor $\Ind(\textup{LCoh}(R_{n+1}))\to \Ind(\textup{LCoh}(R_n))$.

Passing to the limits,
we have a monoidal fully faithful functor
\[
\widehat{\Upsilon}_R=\lim_{n} \Upsilon_{R_n}: \lim_{n}\LMod_{R_n}\to \lim_{n}\Ind(\textup{LCoh}(R_n)).
\]

\begin{Lemma}
\label{contiind}
The $!$-pullback functors $\Ind(\textup{LCoh}(R))\to \Ind(\textup{LCoh}(R_n))$ induce an equivalence
\[
\Ind(\textup{LCoh}(R))\simeq \lim_{n}\Ind(\textup{LCoh}(R_n)).
\]
\end{Lemma}

\Proof
When $R$ arises from a commutative algebra object, the assertion follows from a more general result \cite[Vol.I, Chapter 4, Proposition 6.4.3]{Gai2}.
When $R$ is a connective $\etwo$-algebra, there exists an equivalence
$\Ind(\textup{LCoh}(R_n))\simeq \RMod_{\DD_1(R_n)}$ (cf. \cite[Lemma 2.6 (1)]{CKDA}, \cite[proof of Lemma 3.11]{IH}, \cite[X, Theorem 3.5.1]{DAG}).
This also gives an equivalence $\textup{LCoh}(R_n)\simeq \textup{RPerf}_{\DD_1(R_n)}$ of full subcategories
where $\textup{RPerf}_{\DD_1(R_n)}$ denotes the full subcategory of compact objects.
We regard 
$R$ as the sequential limit $\lim_{n\ge0 }R_n$ of Artin $\eone$-algebras.
By \cite[Lemma 2.6 (2)]{CKDA}, the diagram
\[
\xymatrix{
\Ind(\textup{LCoh}(R_{n+1}))\ar[r] \ar[d] &  \RMod_{\DD_1(R_{n+1})} \ar[d] \\
\Ind(\textup{LCoh}(R_n))\ar[r] &    \RMod_{\DD_1(R_n)}
}
\]
commutes,
where the vertical functor on the left-hand side is the $!$-pullback functor, and the vertical functor on the right-hand side
is the restriction functor along $\DD_1(R_n)\to \DD_1(R_{n+1})$.
To prove the assertion, it will suffice to show that 
the functor
$\RMod_{\DD_1(R)}\to \lim_{n}\RMod_{\DD_1(R_n)}$
induced by the restriction functors along $\DD_1(R_n)\to \DD_1(R)$ is an equivalence.
The functor $\Theta:\Alg_1(\Mod_k)\to \PR$ sending $A$ to $\RMod_{A}$ preserves filtered colimits
(note that $\theta_{\Mod_k}:\Alg_1(\Mod_k)\to (\PR_k)_{\Mod_k/}$ is a left adjoint (see Section~\ref{thetasection} for $\theta_{\Mod_k}$), and the forgetful functor $(\PR_k)_{\Mod_k/}\to \PR_k\to \PR$ preserves sifted colimits).
By Corollary~\ref{limitkoszul}, we have a canonical equivalence $\RMod_{\DD_1(R)}\simeq \colim_{n}\RMod_{\DD_1(R_n)}$ in $\PR$
such that $\colim_{n}\RMod_{\DD_1(R_n)}$ is a colimit in $\PR$ of the diagram
$\RMod_{\DD_1(R_0)}\to\RMod_{\DD_1(R_1)}\to \cdots$ induced by the base change functors.
Passing to right adjoint functors we identify $\colim_{n}\RMod_{\DD_1(R_n)}$
with $\lim_{n}\RMod_{\DD_1(R_n)}$ (i.e., the limit of $\infty$-categories).
Thus, we conclude that $\RMod_{\DD_1(R)}\to \lim_{n}\RMod_{\DD_1(R_n)}$ is an equivalence.
\QED

By Lemma~\ref{contiind}, we have 

\begin{Corollary}
\label{upsilonlimit}
The diagram
\[
\xymatrix{
\LMod_R \ar[r]^{\Upsilon_R} \ar[d] & \Ind(\textup{LCoh}(R)) \ar[d]^{\simeq} \\
\lim_n\LMod_{R_n} \ar[r]_{\textup{fully faithful}}^{\widehat{\Upsilon}_R}&  \lim_n\Ind(\textup{LCoh}(R_{n+1}))
}
\]
commutes up to canonical homotopy.
In particular, $\Upsilon_R:\LMod_R \to \Ind(\textup{LCoh}(R))$ can be identified with 
the canonical functor $\LMod_R\to\lim_n\LMod_{R_n}$ induced by base changes through the embedding $\lim_n\LMod_{R_n}\subset \lim_n\Ind(\textup{LCoh}(R_{n+1}))\simeq \Ind(\textup{LCoh}(R))$.
\end{Corollary}

\begin{Corollary}
\label{omegaend}
Let $\omega_R\in \Ind(\textup{LCoh}(R))$ be the dualizing object, which is a unit object of the monoidal $\infty$-category
$\Ind(\textup{LCoh}(R))$.
The (opposite of the) endomorphism $\etwo$-algebra of $\omega_R$ is $R$.
In particular, $\Upsilon_R$ is fully faithful if we restrict the domain to the single object $R$. 
\end{Corollary}

\Proof
By the equivalences in the diagram of Corollary~\ref{upsilonlimit}, the (opposite of the) endomorphism $\etwo$-algebra of $\omega_R$
is the limit of the (opposite of the) endomorphism $\etwo$-algebra of $\omega_{R_n}$,
 which is canonically equivalent to $\lim_nR_n$. Namely, there exist equivalences
 \[
\Hom_{\Ind(\textup{LCoh}(R))}(\omega_R,\omega_R)\simeq \lim_n\Hom_{\Ind(\textup{LCoh}(R_n))}(\omega_{R_n},\omega_{R_n})\simeq \lim_nR_n\simeq R
\]
as desired. Here $\Hom$ indicates the Hom object defined in $\Mod_k$ (note that $\Ind(\textup{LCoh}(R))$ is a $k$-linear presentable $\infty$-category).
\QED

\begin{Definition}
We define $R'$ in $\Ind(\textup{LCoh}(R))$ as the limit of $\cdots \to R_n\to \cdots \to R_1\to R_0$
taken in $\Ind(\textup{LCoh}(R))$. Here $R_n\in \textup{LCoh}(R)\subset \Ind(\textup{LCoh}(R))$.

More generally, for $M\in \LMod_R$ 
we define $M'\in \Ind(\textup{LCoh}(R))$ as the limit of $\cdots \to \tau^{\ge -n}M\to \cdots \to \tau^{\ge -1}M\to \tau^{\ge0}M$
taken in $\Ind(\textup{LCoh}(R))$. Here $\tau^{\ge -n}M$ is regarded as an object of $\Ind(\textup{LCoh}(R))^{\ge -n}\simeq \LMod_R^{\ge -n}$.
\end{Definition}

Let 
$\psi_R:\Ind(\textup{LCoh}(R))\to \LMod_R$ be the left Kan extension of the inclusion $\textup{LCoh}(R)\to \LMod_R$.
We remark that in \cite{Ind} and \cite{CKDA} the notation $\Psi_R$ is used instead of $\psi_R$. 
We note that $ \Ind(\textup{LCoh}(R))\to \LMod_R$ induces an equivalence
$\Ind(\textup{LCoh}(R))^{\ge0}\simeq \LMod_R^{\ge0}$
(see \cite[Proposition 1.2.4]{Ind} for the commutative case. The proof of the case of $\etwo$-algebras is similar, cf. \cite[the proof of Lemma 3.11]{IH}).
The functor $\LMod_R\to \Ind(\textup{LCoh}(R))$ given by $M\mapsto M'$ is nothing but the right adjoint of 
$\psi_R$.
The quasi-inverse $\LMod_R^{\ge0}\to \Ind(\textup{LCoh}(R))^{\ge0}$ of $ \Ind(\textup{LCoh}(R))^{\ge0}\simeq\LMod_R^{\ge0}$
is given by the restriction of the right adjoint.

\begin{Proposition}
\label{homind}
There exists a canonical equivalence
$\Hom_{\Ind(\textup{LCoh}(R))}(R',M')\simeq M$
in $\Mod_k$.  Here $\Hom$ indicates the Hom object defined in $\Mod_k$.
\end{Proposition}

\Proof
By the definition of $M'$, there exists a canonical equivalence
\[
\Hom_{\Ind(\textup{LCoh}(R))}(R',M')\simeq \lim_n\Hom_{\Ind(\textup{LCoh}(R))}(R',\tau^{\ge-n}M).
\]
If we denote by $\tau^{\le *}_{\Ind}$ and $\tau^{\ge *}_{\Ind}$ the truncations with respect to the $t$-structure on $\Ind(\textup{LCoh}(R))$,
then $\tau^{\ge -n}_{\Ind}R'\simeq \tau^{\ge-n}R$.
Indeed, there exists a cofiber sequence 
\[
\tau^{\le-n}_{\Ind}(R')\to R' \to \tau^{>-n}_{\Ind}(R')
\]
in $\Ind(\textup{LCoh}(R))$.
Since $\tau^{\le-n}_{\Ind}$ preserves small limits
 we have a cofiber sequence
\[
\lim_m(\tau^{\le-n}_{\Ind}\tau^{\ge-m}R)\to \lim_mR_m \to \tau^{>-n}_{\Ind}(\lim_mR_m).
\]
Since 
$\tau^{\le-n}_{\Ind}\tau^{\ge-m}R\simeq \tau^{\le-n}\tau^{\ge-m}R\simeq \tau^{\ge-m}\tau^{\le-n}R$,
it follows that $\lim_m(\tau^{\ge-m}\tau^{\le-n}R)\simeq \tau^{\le-n}_{\Ind}(R')$.
Hence, the above cofiber sequence is obtained from the cofiber sequence
\[
\tau^{\le-n}R\to R \to \tau^{>-n}R
\]
by taking $\tau^{\ge-m}$ and passing to the limit $m\to \infty$. It follows that $\tau^{>-n}_{\Ind}(\lim_mR_m)\simeq \tau^{>-n}R$.
Consequently, we have
\[
\Hom_{\Ind(\textup{LCoh}(R))}(R',\tau^{\ge-n}M)\simeq \Hom_{\Ind(\textup{LCoh}(R))}(\tau^{\ge-n}_{\Ind}R',\tau^{\ge-n}M)\simeq \Hom_{\Ind(\textup{LCoh}(R))}(\tau^{\ge-n}R,\tau^{\ge-n}M).
\]
Recall the equivalence $\Ind(\textup{LCoh}(R))^{\ge-n}\simeq \LMod_R^{\ge-n}$.
There exists an equivalence 
\[
\Hom_{\Ind(\textup{LCoh}(R))}(\tau^{\ge-n}R,\tau^{\ge-n}M)\simeq \Hom_{\LMod_R}(\tau^{\ge-n}R,\tau^{\ge-n}M)\simeq \Hom_{\LMod_R}(R,\tau^{\ge-n}M)\simeq \tau^{\ge-n}M
\]
in $\LMod_R$.
Hence we conclude $\lim_n\Hom_{\Ind(\textup{LCoh}(R))}(R',\tau^{\ge-n}M)\simeq \lim_n\tau^{\ge-n}M\simeq M$.
\QED

\begin{Corollary}
\label{nonecor}
The functor $\psi_R:\Ind(\textup{LCoh}(R))\to \LMod_R$
sends $R'$ to $R$.
\end{Corollary}

\begin{Corollary}
\label{nonecor2}
The left completion functor $\psi_R$ induces an equivalence
$\Hom_{\Ind(\textup{LCoh}(R))}(R',R')\to \Hom_{\LMod_R}(R,R)\simeq R$.
\end{Corollary}

\begin{Corollary}
\label{CKDA2}
Let $R$ be an almost Artin $\etwo$-algebra.
Then there exists a canonical equivalence
\[
\Ind(\textup{LCoh}(R))\simeq\End_{\LMod_{\DD_2(R)}}(\Mod_k)
\]
 in $\Alg_1(\PR_k)$.
\end{Corollary}

\Proof
We note $R\simeq \lim_{n}R_n$. 
We consider the Koszul dual $\End_{\Ind(\textup{LCoh}(R_n))}(\Mod_k)$ of $\Ind(\textup{LCoh}(R_n))\to \Mod_k$ (in $\Alg_1^+(\PR_k)$) and 
the universal pairing 
\[
u:\Ind(\textup{LCoh}(R_n))\otimes \End_{\Ind(\textup{LCoh}(R_n))}(\Mod_k)\to \Mod_k,
\]
which is a morphism in $\Alg_1(\PR_k)$.
This induces $\Ind(\textup{LCoh}(R_n))\to \End_{\End_{\Ind(\textup{LCoh}(R_n))}(\Mod_k)}(\Mod_k)$ in $\Alg_1(\PR_k)$.
Let $\LMod_{\DD_2(R_n)}\to \End_{\Ind(\textup{LCoh}(R_n))}(\Mod_k)$ be the algebraization functor (defined as a morphism of $\Alg_1(\PR_k)$)
of $\End_{\Ind(\textup{LCoh}(R))}(\Mod_k)\in \Alg_1(\PR_k)$ (see \cite[Lemma 3.3]{CKDA}).
This induces the restriction $\End_{\End_{\Ind(\textup{LCoh}(R_n))}(\Mod_k)}(\Mod_k)\to \End_{\LMod_{\DD_2(R_n)}}(\Mod_k)$.
We consider the composite functor 
\[
\Ind(\textup{LCoh}(R_n))\to \End_{\End_{\Ind(\textup{LCoh}(R_n))}(\Mod_k)}(\Mod_k)\to \End_{\LMod_{\DD_2(R_n)}}(\Mod_k)
\]
in $\Alg_1(\PR_k)$.
The sequence $\Ind(\textup{LCoh}(R_n))\to \End_{\End_{\Ind(\textup{LCoh}(R_n))}(\Mod_k)}(\Mod_k)\to \End_{\LMod_{\DD_2(R_n)}}(\Mod_k)$
is covariantly functorial with respect to $\cdots \to R_2\to R_1\to R_0$ (cf. \cite[Section 3]{CKDA}).
Thus, we have
\[
\lim_{n}\Ind(\textup{LCoh}(R_n)) \to \lim_{n}\End_{\LMod_{\DD_2(R_n)}}(\Mod_k)\simeq \lim_{n}\RMod_{k\otimes_{\DD_2(R_n)}k}.
\]
The arrow on the left-hand side is a morphism in $\Alg_1(\PR_k)$, and the equivalence on the right-hand side means an equivalence in $\PR_k$.
The limit on the left-hand side is the limit of the diagram induced by $!$-pullback functors.
The limit on the right-hand side is the limit of the diagram induced by the restriction functors along $k\otimes_{\DD_2(R_0)}k\to k\otimes_{\DD_2(R_1)}k\to k\otimes_{\DD_2(R_2)}k\to  \cdots$.
From the proof of Lemma~\ref{contiind}, we have
\[
 \lim_{n}\RMod_{k\otimes_{\DD_2(R)}k}\simeq \lim_n\RMod_{\DD_1(R_n)}\simeq \RMod_{\DD_1(R)}.
\]
Also, by Proposition~\ref{dbar} there exist equivalences $\End_{\LMod_{\DD_2(R)}}(\Mod_k)\simeq \RMod_{k\otimes_{\DD_2(R)}k}\simeq \RMod_{\DD_1(R)}$ in $\PR_k$.
Consequently, the canonical morphism $\End_{\LMod_{\DD_2(R)}}(\Mod_k) \to  \lim_{n}\End_{\LMod_{\DD_2(R_n)}}(\Mod_k)$ in $\Alg_1(\PR_k)$
is an equivalence. By Lemma~\ref{contiind}, the canonical functor $\Ind(\textup{LCoh}(R))\to \lim_{n}\Ind(\textup{LCoh}(R_n))$
is an equivalence. Thus, it is enough to prove that each $\Ind(\textup{LCoh}(R_n)) \to \End_{\LMod_{\DD_2(R_n)}}(\Mod_k)$ is an equivalence.
This is proved in \cite[Theorem 1.2]{CKDA}.
\QED

The proof of the first assertion of Proposition~\ref{leftcomp2} in \cite{IH} also works in the case when $R$ is almost Artin.
There exist equivalences $\RMod_{k\otimes_{\DD_2(R)}k}\simeq \RMod_{\DD_1(R)}\simeq \Ind(\textup{LCoh}(R))$.
For the first equivalence, we use Proposition~\ref{dbar}.
Through these equivalences, the functor
$\eta:\Mod_k\otimes_{\RMod_{\DD_2(R)}}\Mod_k\to \LMod_R$ can be identified with
$\psi_R:\Ind(\textup{LCoh}(R))\to \LMod_R$, which is a left completion.
We also note that the left completion $\psi_R$
can be identified with the composite
\[
\Ind(\textup{LCoh}(R))\to \lim_{n} \Ind(\textup{LCoh}(R))^{\ge -n}\simeq \lim_{n} \LMod_R^{\ge -n}\simeq \LMod_R.
\]
The first functor is induced by the truncation functors $\tau_{\Ind}^{\ge -n}:\Ind(\textup{LCoh}(R))\to \Ind(\textup{LCoh}(R))^{\ge -n}$.
The middle equivalence is induced by equivalences $\Ind(\textup{LCoh}(R))^{\ge -n}=\Ind(\textup{LCoh}(R)^{\ge -n})\simeq \LMod_R^{\ge -n}$.
The right equivalence comes from the truncation functors on $\LMod_R$ with respect to the left complete $t$-structure ($R$ is connective).
From the proof of Proposition~\ref{homind}, $R'$ maps to the $\{\tau^{\ge -n}R\}\in \lim_{n} \LMod_R^{\ge -n}$ where 
$\{\tau^{\ge -n}R\}$ informally indicates the evident object of $\lim_{n} \LMod_R^{\ge -n}$ such that $\tau^{\ge -n}R\in \LMod_R^{\ge -n}$.

\vspace{3mm}

{\it Proof of Proposition~\ref{leftcomp2}.}
We prove the latter assertion:
We will prove that the target $\LMod_R$ is a free $\LMod_R$-module of rank one.
Taking left completions, the action functor $\LMod_R\times \RMod_{k\otimes_{\DD_2(R)}k} \to \RMod_{k\otimes_{\DD_2(R)}k}$
induces the module action functor $\LMod_R\times \LMod_R\to \LMod_R$, where $\LMod_R$ is identified with a left completion of $\RMod_{k\otimes_{\DD_2(R)}k}$.
By the same argument as in \cite{IH}, we are reduced to showing that
there exists an object $E\in \RMod_{k\otimes_{\DD_2(R)}k}$ such that the composite 
\[
\rho:\LMod_R\times \{E\}\to \LMod_R\times \RMod_{k\otimes_{\DD_2(R)}k}\to  \RMod_{k\otimes_{\DD_2(R)}k}  \stackrel{\eta}{\to} \LMod_R
\]
is an equivalence.
We take $E$ to be $R'$ through the equivalences
$\RMod_{k\otimes_{\DD_2(R)}k}\simeq  \Ind(\textup{LCoh}(R))$.
By Corollary~\ref{nonecor}, $(R,E)\in \LMod_R\times \{E\}$ maps to $R\in \LMod_R$.
By Corollary~\ref{upsilonlimit}, Corollary~\ref{omegaend}, Corollary~\ref{nonecor2} and Corollary~\ref{CKDA2}, 
we see that $\rho:\LMod_R\to \LMod_R$
is fully faithful after the restriction to the full subcategory consisting of the single object $R$.
Since $R$ is a single compact generator of $\LMod_R$, and $\rho$ is a colimit-preserving functor,
we conclude that $\rho:\LMod_R\to \LMod_R$ is an equivalence of stable $\infty$-categories.
\QED

\section{Hochschild cohomology}

\label{hochsec}

We study Hochschild cohomology of $t$-structured stable presentable $\infty$-categories.

\subsection{}

Let $n$ be a nonnegative integer. Let $k$ be a connective commutative ring spectrum.
Let $\DDD$ be any of an object of $(n+1)\PR_k$ or an object of $(n+1)\PRTT_k$ or an object of $(n+1)\PRTTT_k$.
We will define Hochschild cohomology of $\DDD$.

\begin{Definition}
\label{Hochschilddef}
Let $\mathsf{P}$ be any of $\PR_k$ or $\PRTT_k$ or $\PRTTT_k$.
Let $\DDD$ be an object of $(n+1)\mathsf{P}$.
For $n\ge1$, we consider an endomorphism algebra object $\End(\DDD)\in \Alg_1(\Mod_{n\mathsf{P}}(\vCat))$
(note that an internal Hom object $\Hom(\DDD,\DDD)$ does not generally exist in $(n+1)\mathsf{P}$ 
but the internal Hom object $\Hom(\DDD,\DDD)$ exists in $\Mod_{n\mathsf{P}}(\vCat)$ since $\Mod_{n\mathsf{P}}(\vCat)$ is $\mathbb{V}$-presentable).
We assume that $\End(\DDD)$ has an algebraization functor $n\mathsf{P}_C\to \End(\DDD)$ (see Definition~\ref{algebraizationdef}).
In this situation, we say that $\DDD$ has Hochschild cohomology, and
we refer to $C\in \Alg_{n+2}(\Mod_k)$ as Hochschild cohomology of $\DDD$.

Suppose that $n=0$. We consider the internal endomorphism algebra object $\End(\DDD)\in \Alg_1(\PR_k)$ when $\mathsf{P}=\PR_k$.
When $\mathsf{P}$ is either $\PRTT_k$ or $\PRTTT_k$, we consider the internal endomorphism algebra object $\End(\DDD)\in \Alg_1(\PRT_k)$ 
by regarding $\DDD$ as an object of $\PRT_k$.
Let $\LMod_C\to \End(\DDD)$ be an algebraization functor, which always exists in this case.
We refer to $C\in \Alg_{n+2}(\Mod_k)$ as Hochschild cohomology of $\DDD$.

We write $\HH^\bullet_{(n+1)\mathsf{P}}(\DDD/k)$ for $C$.
It is common to refer to $\HH^\bullet_{(n+1)\mathsf{P}}(\DDD/k)$ as the Hochschild cochains etc., but unless stated otherwise
we refer to it as Hochschild cohomology.
\end{Definition}

\begin{Warning}
If $(n+1)\mathsf{P}$, in which $\DDD$ is defined, is clear, we often write $\HH^\bullet(\DDD/k)$ for $\HH^\bullet_{(n+1)\mathsf{P}}(\DDD/k)$.
\end{Warning}

\begin{Example}
We assume that $\End(\DDD)$ has a good endomorphism algebra of a unit (see Definition~\ref{goodendomorphism}).
We consider $\mathsf{E}^n(\End(\DDD))$, which is an $\mathbf{E}_{n+1}$-algebra in $\mathsf{P}$ (see Definition~\ref{goodendomorphism} for $\mathsf{E}^n$).
In this situation, $\DDD$ has
Hochschild cohomology $\HH^\bullet_{(n+1)\mathsf{P}}(\DDD/k)\in \Alg_{n+2}(\Mod_k)$, defined as the endomorphism algebra of a unit object of $\mathsf{E}^n(\End(\DDD))$
(see Lemma~\ref{goodalgebraization}).
\end{Example}

\begin{Remark}
\label{plaineasy}
We consider the case of $\mathsf{P}=\PR_k$.
Let us attempt to define Hochschild cohomology.
Let $\CCC$ be an object of $\Mod_{i\PR_k}(\vCat)$. Here we assume $i\ge1$.
For any two objects $M,N$ of $\CCC$
there exists a Hom object $\Hom(M,N)$ defined as an object of $\Mod_{(i-1)\PR_k}(\vCat)$
(generally not in $i\PR_k$)
since $\Mod_{\mathsf{(i-1)\PR_k}}(\vCat)$ is $\kappa_0$-compactly generated such that $i\PR_k$ consists of $\kappa_0$-compact objects (\cite[5.1.7]{Ste}). 
Thus, one obtains $\mathsf{E}^n(\End(\DDD))\in \Alg_{n+1}(\Mod_{\Mod_k}(\vCat))$.
However, for two objects $M,N$ of $\Mod_{\Mod_k}(\vCat)$,
a Hom object $\Hom(M,N)$ does not generally exist in $\Mod_k$, while
$\Hom(M,N)$ exists as a $k$-module spectrum in the enlarged universe $\mathbb{V}$.
This motivates us to study when Hochschild cohomology $\HH^\bullet(\DDD/k)$ exists in $\Mod_k$ (or $\Alg_{n+2}(\Mod_k)$).
\end{Remark}

As observed in \cite[Theorem C]{Ao}, for $n\ge3$ stable monoidal presentable $(\infty,n)$-categories do not generally have a presentable endomorphism algebra of a unit.
Thus, we need to consider when $\DDD$ has Hochschild cohomology.
Moreover, if $\DDD \in (n+1)\PRTT_k$ has Hochschild cohomology $\HH^\bullet(\DDD/k)$,
one can ask whether $\HH^\bullet(\DDD/k)$ coincides with the Hochschild cohomology $\HH^\bullet(\underline{\DDD}/k)$ of
the underlying stable $(\infty,n+1)$-category $\underline{\DDD}\in (n+1)\PR_k$, when the latter exists. We remark that when $n=0$
Hochschild cohomology does not depend on $t$-structures, so that Hochschild cohomology is determined by the underlying $k$-linear presentable $\infty$-category (see \cite{IH}).
In the subsequent subsections, we study these issues.

\subsection{}

\label{enrichsec}

We consider the case of $t$-structured (stable) $(\infty,2)$-categories.
The typical example is $\LMod_{\MMM^\otimes}(\PRTT_k)$ such that $\MMM^\otimes$ is an object of $\Alg_1(\PRTT_k)$, that is,
a $k$-linear $\eone$-monoidal $\infty$-category
with an accessible right complete $t$-structure.

Let $\mathsf{P}$ be a commutative algebra object in $\vCat$.
We consider a small category enriched over $\mathsf{P}$.
There are several approaches to enriched $\infty$-categories, see e.g. \cite{GH}, \cite{Hinen}, \cite{RZ}. We adopt the formalism developed in Hinich \cite{Hinen}.
We briefly review the formulation of enriched categories.

For a small set $X$, we consider the functor category $\Fun(X,\mathsf{P})$, which we regard as a $\mathsf{P}$-module.
Let $\End(\Fun(X,\mathsf{P}))$
denote the internal endomorphism algebra of the $\mathsf{P}$-module, which is defined as an object of $\Alg_1(\Mod_{\mathsf{P}}(\vCat))$.
The underlying $\mathsf{P}$-module of $\End(\Fun(X,\mathsf{P}))$ is equivalent to
$\Fun(X\times X,\mathsf{P})$.
Indeed, $\Fun(X,\mathsf{P})\in \Mod_{\mathsf{P}}(\vCat)$ is a dualizable $\mathsf{P}$-module (in particular, $\Fun(X,\mathsf{P})\in 2\mathsf{P}$)
whose dual $\Fun(X,\mathsf{P})^\vee$ is equivalent to $\Fun(X,\mathsf{P})$. Thus, the Hom object
$\Hom_{\mathsf{P}}(\Fun(X,\mathsf{P}),\Fun(X,\mathsf{P}))$ exists in $2\mathsf{P}$,
and it is equivalent to $\Fun(X,\mathsf{P})^\vee\otimes_{\mathsf{P}}\Fun(X,\mathsf{P})\simeq \Fun(X\times X, \mathsf{P})$.
Here the subscript $\mathsf{P}$ in $\Hom_{\mathsf{P}}(\Fun(X,\mathsf{P}),\Fun(X,\mathsf{P}))$ is not the $\infty$-category
on which $\Fun(X,\mathsf{P})$ lives (but $\Fun(X,\mathsf{P})$ lies in $2\mathsf{P}$): The subscript $\mathsf{P}$ means "$\mathsf{P}$-linear".
For $A,B \in\Fun(X\times X, \mathsf{P})$, the tensor product $A\otimes B$ is given informally by
$A\otimes B(x,z)=\colim_{y\in X} A(x,y)\otimes B(y,z)$.
In this setting, a $\mathsf{P}$-enriched category $\mathcal{A}$ having the prescribed set $X$ of objects is an algebra object $\mathcal{A}\in \Alg_1(\End(\Fun(X,\mathsf{P})))$. 
We think of the underlying functor $\mathcal{A}:X\times X\to \mathsf{P}$ as the assignment 
$(x,y)\mapsto \mathcal{H}om(x,y)\in \mathsf{P}$, and the associative algebra structure on $\mathcal{A}$ encodes the data of composition and units.

We note that
$\Fun(X,\mathsf{P})$ has the tautological $\End(\Fun(X,\mathsf{P}))$-module structure in the natural way.
We define the $\infty$-category of $\mathcal{A}$-module objects in $\Fun(X,\mathsf{P})$:
\[
P_{\mathsf{P}}(\mathcal{A})=\LMod_{\mathcal{A}}(\Fun(X,\mathsf{P})).
\]
We think of an object $F$ of $P_{\mathsf{P}}(\mathcal{A})$ as a presheaf on the enriched category $\mathcal{A}$
taking values in $\mathsf{P}$: The underlying functor $F:X\to \mathsf{P}$ assigns each object $x\in X$ to $F(x)\in \mathsf{P}$.
Note that  $P_{\mathsf{P}}(\mathcal{A})$ is a $\mathsf{P}$-module, that is, an object of $\Mod_{\mathsf{P}}(\vCat)$.

\begin{Example}
\label{exendoalg}
Suppose that $X$ is a one-point category having a single object $\ast$ such that $\Map_X(\ast,\ast)$ is a contractible space.
Then $\mathcal{A}$ amounts to an object of $\Alg_1(\mathsf{P})$, and $P_{\mathsf{P}}(\mathcal{A})$ is $\LMod_{\mathcal{A}}(\mathsf{P})$.

\end{Example}

\begin{Proposition}
\label{enrichdual}
$P_{\mathsf{P}}(\mathcal{A})$ belongs to $2\mathsf{P}$. Namely, it is $\kappa_0$-compact in $\Mod_{\mathsf{P}}(\vCat)$. Moreover, $P_{\mathsf{P}}(\mathcal{A})$ is dualizable in $2\mathsf{P}$.
\end{Proposition}

\begin{Proposition}
\label{enrichhoch}
When $\mathsf{P}$ is any of $\PR_k$ or $\PRTT_k$ or $\PRTTT_k$, $P_{\mathsf{P}}(\mathcal{A})$ has Hochschild cohomology $\HH^\bullet_{2\mathsf{P}}(P_{\mathsf{P}}(\mathcal{A})/k)\in \Alg_3(\Mod_k)$.
\end{Proposition}

For applications in any context depending on $\mathsf{P}$, it is important to understand the construction of Hochschild cohomology. We proceed step by step.

{\it Proof of Proposition~\ref{enrichdual}.}
The tensor product in 
$\Mod_{\mathsf{P}}(\vCat)$
preserves colimits separately in each variable.
It is enough to show that $P_{\mathsf{P}}(\mathcal{A})$ is dualizable in $\Mod_{\mathsf{P}}(\vCat)$ whose dual is
given by $\RMod_{\mathcal{A}}(\Fun(X,\mathsf{P}))$.
For this purpose, we define an evaluation map and a coevaluation map.
For ease of notation, we write $\mathsf{P}^X$ for $\Fun(X,\mathsf{P})$.
According to \cite[4.8.4.6]{HA}, there exists a canonical equivalence
$\RMod_{\mathcal{A}}(\End(\mathsf{P}^X))\otimes_{\End(\mathsf{P}^X)}\mathsf{P}^X\simeq \RMod_{\mathcal{A}}(\mathsf{P}^X)$.
Moreover, $\mathsf{P}^X$ is dualizable and self-dual in $\Mod_{\mathsf{P}}(\vCat)$ so that $\mathsf{P}^X\in 2\mathsf{P}$.
Thus we have equivalences
\begin{eqnarray*}
\RMod_{\mathcal{A}}(\mathsf{P}^X)\otimes_{\mathsf{P}} \LMod_{\mathcal{A}}(\mathsf{P}^X) &\simeq& \RMod_{\mathcal{A}}(\End(\mathsf{P}^X))\otimes_{\End(\mathsf{P}^X)}\mathsf{P}^X\otimes_{\mathsf{P}}\mathsf{P}^X\otimes_{\End(\mathsf{P}^X)} \LMod_{\mathcal{A}}(\End(\mathsf{P}^X)) \\
&\simeq& \RMod_{\mathcal{A}}(\End(\mathsf{P}^X))\otimes_{\End(\mathsf{P}^X)}\Hom(\mathsf{P}^X,\mathsf{P}^X)\otimes_{\End(\mathsf{P}^X)} \LMod_{\mathcal{A}}(\End(\mathsf{P}^X)) \\
&\simeq& \RMod_{\mathcal{A}}(\End(\mathsf{P}^X))\otimes_{\End(\mathsf{P}^X)} \LMod_{\mathcal{A}}(\End(\mathsf{P}^X)) \\
&\simeq& \LMod_{\mathcal{A}^{op}\otimes\mathcal{A}}(\End(\mathsf{P}^X)).
\end{eqnarray*}
Using these equivalences, we define an evaluation map
$u:\RMod_{\mathcal{A}}(\mathsf{P}^X)\otimes_{\mathsf{P}} \LMod_{\mathcal{A}}(\mathsf{P}^X)\to \mathsf{P}$
as the composite of
\[
\LMod_{\mathcal{A}^{op}\otimes\mathcal{A}}(\End(\mathsf{P}^X))\to \End(\mathsf{P}^X)
\]
given by $\mathcal{E}\mapsto \mathcal{A}_{\Delta}\otimes_{\mathcal{A}^{op}\otimes\mathcal{A}}\mathcal{E}$ 
and the canonical evaluation map
\[
\End(\mathsf{P}^X)=(\mathsf{P}^X)^\vee\otimes_\mathsf{P} \mathsf{P}^X\to \mathsf{P}.
\]
Here $\mathcal{A}_{\Delta}$ denotes $\mathcal{A}$ endowed with the diagonal right $\mathcal{A}^{op}\otimes\mathcal{A}$-module structure.
The coevaluation map $v:\mathsf{P}\to \LMod_{\mathcal{A}}(\mathsf{P}^X)\otimes_{\mathsf{P}} \RMod_{\mathcal{A}}(\mathsf{P}^X)$
is defined to be
a unique $\mathsf{P}$-linear functor which carries the unit of $\mathsf{P}$ to
the diagonal module $\mathcal{A}_{\Delta}\in \LMod_{\mathcal{A}^{op}\otimes\mathcal{A}}(\End(\mathsf{P}^X))\simeq \LMod_{\mathcal{A}}(\mathsf{P}^X)\otimes_{\mathsf{P}} \RMod_{\mathcal{A}}(\mathsf{P}^X)$.
It remains to prove that $u$ and $v$ satisfy the triangle identities.
It suffices to show that 
\begin{eqnarray*}
\LMod_{\mathcal{A}}(\mathsf{P}^X) &\stackrel{v\otimes \textup{id}}\longrightarrow & \LMod_{\mathcal{A}}(\mathsf{P}^X)\otimes_{\mathsf{P}} \RMod_{\mathcal{A}}(\mathsf{P}^X)\otimes_{\mathsf{P}}\LMod_{\mathcal{A}}(\mathsf{P}^X) \\
&\stackrel{\textup{id}\otimes u}{\longrightarrow}& \LMod_{\mathcal{A}}(\mathsf{P}^X)
\end{eqnarray*}
is the identity functor (another triangle identity is similar).
Indeed, unfolding the definitions, the composite is 
given by $\mathcal{P}\mapsto (\mathcal{A}\otimes \mathcal{P})\otimes_{(\mathcal{A}^{op}\otimes\mathcal{A})}\mathcal{A}\simeq \mathcal{P}$.
Namely, the composite is equivalent to the identity.
\QED

For ease of notation, we write $P(\mathcal{A})=P_{\mathsf{P}}(\mathcal{A})$.  
By Proposition~\ref{enrichdual}, $P(\mathcal{A})\in 2\PRTT_k$, and it is dualizable in $2\mathsf{P}$.
Next we will consider $\End(P(\mathcal{A}))$.
Note that
$P(\mathcal{A})=\LMod_{\mathcal{A}}(\Fun(X,\mathsf{P}))$, and its dual is $\RMod_{\mathcal{A}}(\Fun(X,\mathsf{P}))$ (see the proof of Proposition~\ref{enrichdual}). 
Thus, the internal Hom object
$\Hom_{\mathsf{P}}(P(\mathcal{A}),P(\mathcal{A}))$ is equivalent to $\RMod_{\mathcal{A}}(\Fun(X,\mathsf{P}))\otimes_{\mathsf{P}}\LMod_{\mathcal{A}}(\Fun(X,\mathsf{P}))$.
Thus, $\End(P(\mathcal{A}))$ is defined as an object of $\Alg_1(2\mathsf{P})$. 
Its underlying $\infty$-category is ${}_{\mathcal{A}}\BMod_{\mathcal{A}}(\End(\Fun(X,\mathsf{P})))$ where the diagonal bimodule $\mathcal{A}$ is a unit object.

\begin{Lemma}
\label{bimodulehom}
Suppose that $\mathsf{P}$ has internal Hom objects and small limits.
For $\mathcal{E},\mathcal{F}\in {}_{\mathcal{A}}\BMod_{\mathcal{A}}(\End(\Fun(X,\mathsf{P})))$, there exists a Hom object
\[
\Hom_{{}_{\mathcal{A}}\BMod_{\mathcal{A}}(\End(\Fun(X,\mathsf{P})))}(\mathcal{E},\mathcal{F})
\]
in $\mathsf{P}$.
\end{Lemma}

\Proof
To observe this, using the comonad of the free-forget adjunction of bimodules
\[
\xymatrix{
\End(\Fun(X,\mathsf{P})) \ar@<0.5ex>[r]^{\mathcal{A}\otimes(-)\otimes\mathcal{A}}  &  {}_{\mathcal{A}}\BMod_{\mathcal{A}}(\End(\Fun(X,\mathsf{P})))   \ar@<0.5ex>[l]^{\textup{forget}}  
}
\]
we take a simplicial object $\mathcal{R}_\bullet$ in ${}_{\mathcal{A}}\BMod_{\mathcal{A}}(\End(\Fun(X,\mathsf{P})))$ with an augmentation $\epsilon:\mathcal{R}_\bullet\to \mathcal{E}$ 
such that each $n$-th term is a free bimodule $\mathcal{A}\otimes \mathcal{E}(n)\otimes \mathcal{A}$,
and $\epsilon$ exhibits $\mathcal{E}$ as a colimit (see \cite[Proposition 4.7.3.14, 4.7.3.13]{HA}).
(This simplicial object is usually called the bar resolution.) 
Then we have
\[
\Hom_{{}_{\mathcal{A}}\BMod_{\mathcal{A}}(\End(\Fun(X,\mathsf{P})))}(\mathcal{E},\mathcal{F})\simeq \lim_{[n]\in \Delta}\Hom_{\Fun(X\times X,\mathsf{P})}( \mathcal{E}(n),\mathcal{F}).
\]
By the assumption, $\mathsf{P}$ admits internal Hom objects so that
$\Fun(X\times X,\mathsf{P})$ admits Hom objects in $\mathsf{P}$.
Thus each $\Hom_{\Fun(X\times X,\mathsf{P})}( \mathcal{E}(n),\mathcal{F})$ exists in $\mathsf{P}$.
Since
$\mathsf{P}$ has small limits by our assumption, $\lim_{[n]\in \Delta}\Hom_{\Fun(X\times X,\mathsf{P})}( \mathcal{E}(n),\mathcal{F})$ exists in $\mathsf{P}$.
\QED

We now focus on the case when both $\mathcal{E}$ and $\mathcal{F}$ are the diagonal bimodule $\mathcal{A}$.

\begin{Corollary}
\label{bimoduleexist}
The following holds.
\begin{enumerate}

\item
Suppose that $\mathsf{P}$ is either $\PR_k$ or $\PRTT_k$.
There exists 
\[
\Hom_{{}_{\mathcal{A}}\BMod_{\mathcal{A}}(\End(\Fun(X,\mathsf{P})))}(\mathcal{A},\mathcal{A})
\]
in $\Alg_2(\mathsf{P})$.

\item
Suppose that $\mathsf{P}=\PRTTT_k$.
There exists a Hom object  
\[
\Hom_{{}_{\mathcal{A}}\BMod_{\mathcal{A}}(\End(\Fun(X,\PRTTT_k)))}(\mathcal{A},\mathcal{A})
\]
in $\PRTT_k$, where ${}_{\mathcal{A}}\BMod_{\mathcal{A}}(\End(\Fun(X,\PRTTT_k)))$ is regarded as a $\PRTT_k$-module.
\end{enumerate}
\end{Corollary}

\Proof
Note first that $\PRT_k$ and $\PRTT_k$ admit internal Hom objects (cf. Remark~\ref{inthomtstr}):
For $\mathcal{P}, \mathcal{Q}\in \PRT_k$, an internal Hom object $\Hom_{\PRT_k}(\mathcal{P},\mathcal{Q})$ from $\mathcal{P}$ to $\mathcal{Q}$ in $\PRT_k$ is 
given  by the pair $(\Hom_{\PR_k}(\underline{\mathcal{P}},\underline{\mathcal{Q}}), \Hom_{\PR_k}(\underline{\mathcal{P}},\underline{\mathcal{Q}})^{\le0})$
such that $\Hom_{\PR_k}(\underline{\mathcal{P}},\underline{\mathcal{Q}})\in \PR_k$ 
is the internal Hom object from the underlying object $\underline{\mathcal{P}}$ to $\underline{\mathcal{Q}}$, and 
 $\Hom_{\PR_k}(\underline{\mathcal{P}},\underline{\mathcal{Q}})^{\le0}$ is the full subcategory, spanned by right $t$-exact functors.
Thus, $\PRTT_k$ also admits internal Hom objects (by taking right completions of internal Hom objects in $\PRT_k$).
Moreover, $\PRTT_k$ has small limits (see Proposition~\ref{PRTcolimitlimit}).
The $\infty$-category $\PR_k$ has small limits (see \cite[5.5.3.13]{HTT}, \cite[4.2.3.3]{HA})
and has internal Hom objects (see e.g., \cite[Lemma 5.1]{ID}).
Thus, (1) follows from Lemma~\ref{bimodulehom}.

Next, we prove (2).
From the proof of Lemma~\ref{bimodulehom}, 
we take $\epsilon:\mathcal{R}_\bullet\to \mathcal{A}$ 
such that each $n$-th term is a free bimodule $\mathcal{A}\otimes \mathcal{E}(n)\otimes \mathcal{A}$,
and $\epsilon$ exhibits $\mathcal{A}$ as a colimit.
In this case, $\mathcal{E}(n)\simeq \mathcal{A}^{\widehat{\otimes} 2n+1}$ where $\widehat{\otimes}$ indicates the tensor product
in $\PRTTT_k$.
We consider $\Hom_{\Fun(X\times X,\PRTTT_k)}(\mathcal{A}^{\widehat{\otimes} 2n+1},\mathcal{A})$.
This is the 
Hom object defined as an object of $\PRTT_k$ (not of $\PRTTT_k$), where we regard $\Fun(X\times X,\PRTTT_k)$ as a $\PRTT_k$-module.
We have the limit $\lim_{[n]\in \Delta}\Hom_{\Fun(X\times X,\PRTTT_k)}( \mathcal{A}^{\widehat{\otimes} 2n+1},\mathcal{A})$
taken in $\PRTT_k$, so that the Hom object $\Hom_{{}_{\mathcal{A}}\BMod_{\mathcal{A}}(\End(\Fun(X,\PRTTT_k)))}(\mathcal{A},\mathcal{A})$
defined in $\PRTT_k$ is given by $\lim_{[n]\in \Delta}\Hom_{\Fun(X\times X,\PRTTT_k)}( \mathcal{A}^{\widehat{\otimes} 2n+1},\mathcal{A})$.
\QED

We note that
$\Hom_{{}_{\mathcal{A}}\BMod_{\mathcal{A}}(\End(\Fun(X,\mathsf{P})))}(\mathcal{A},\mathcal{A})$
is a categorical analogue of the Hochschild complex of $\mathcal{A}$, that is, a Drinfeld center.
Indeed, we can compute it as $\lim_{[n]\in \Delta}\Hom_{\Fun(X\times X,\mathsf{P})}( \mathcal{A}^{\otimes 2n+3},\mathcal{A})$, using the above simplicial object $\mathcal{R}_\bullet$
with an augmentation $\mathcal{R}_\bullet\to \mathcal{A}$.

Suppose that $\mathsf{P}=\PRTT_k$.
We recall the adjoint pair $\theta_{0\PRTT_k,0}:\Alg_1(\Mod_k)\rightleftarrows (\PRTT_k)_{\Mod_k/}:E_{\PRTT_k}^+$
where the symmetric monoidal functor $\theta_{0\PRTT_k,0}$ sends $B$ to $(\LMod_B,B)$, while $E^+$ sends $(\mathcal{A},A)$ to $\End(A)^{op}$
(see Lemma~\ref{goodalgebraization} and its proof).
Applying $\Alg_2(-)$ to the right adjoint, we have
$\Alg_2(E_{\PRTT_k}^+):\Alg_2(\PRTT_k)\to \Alg_2(\Alg_1(\Mod_k))\simeq \Alg_3(\Mod_k)$ where the final equivalence is Dunn additivity theorem.
Then Hochschild cohomology $\HH^\bullet(P(\mathcal{A})/k)\in \Alg_{3}(\Mod_k)$ is defined to be the image of $\Hom_{{}_{\mathcal{A}}\BMod_{\mathcal{A}}(\End(\Fun(X,\PRTT_k)))}(\mathcal{A},\mathcal{A})$ under
$\Alg_2(E_{\PRTT_k}^+)$, that is
the endomorphism algebra of a unit object of $\Hom_{{}_{\mathcal{A}}\BMod_{\mathcal{A}}(\End(\Fun(X,\mathsf{P})))}(\mathcal{A},\mathcal{A})$ (cf. Example~\ref{exendoalg}).
By construction, it has the required property.

Suppose that $\mathsf{P}=\PRTTT_k$. We take into account the localization adjoint pair $\widehat{L}:\PRTT_k\rightleftarrows \PRTTT_k$.
The Hochschild cohomology $\HH^\bullet(P(\mathcal{A})/k)\in \Alg_{3}(\Mod_k)$ is defined to be the image of 
\[
\Hom_{{}_{\mathcal{A}}\BMod_{\mathcal{A}}(\End(\Fun(X,\PRTTT_k)))}(\mathcal{A},\mathcal{A})\in \Alg_2(\PRTT_k)
\]
under
$\Alg_2(E_{\PRTT_k}^+)$.

Suppose that $\mathsf{P}=\PR_k$.
Using the adjoint pair $\theta_{0\PR_k,0}:\Alg_1(\Mod_k)\rightleftarrows (\PR_k)_{\Mod_k/}:E_{\PR_k}^+$ instead of $(\theta_{0\PRTT_k,0},E_{\PRTT_k}^+)$,
one can construct Hochschild cohomology in the same way as above.

\vspace{3mm}

{\it Proof of Proposition~\ref{enrichhoch}.} The above construction shows our assertion.
\QED

We prove that in this setting Hochschild cohomology is independent of the choice of $\mathsf{P}$.

\begin{Proposition}
\label{comparison1}
Let $\AAA$ be a $\PRTT_k$-enriched category having a prescribed set $X$ of objects (i.e., an algebra object $\mathcal{A}\in \Alg_1(\End(\Fun(X,\PRTT_k)))$). 

\begin{enumerate}
\item
Consider $P_{\PRTT_k}(\AAA)\in 2\PRTT_k$. Let  $\underline{P_{\PRTT_k}(\AAA)}\in 2\PR_k$ denote
its underlying $k$-linear stable presentable $(\infty,2)$-category $\underline{P_{\PRTT_k}(\AAA)}\in 2\PR_k$, which is defined by $\underline{P_{\PRTT_k}(\AAA)}=P_{\PRTT_k}(\AAA)\otimes_{\PRTT_k}\PR_k$.
Then there exists a canonical equivalence $\HH^\bullet_{2\PRTT_k}(P_{\PRTT_k}(\AAA)/k)\simeq \HH^\bullet_{2\PR_k}(\underline{P_{\PRTT_k}(\AAA)}/k)$ in $\Alg_3(\Mod_k)$.

\item
Suppose that $\AAA=\widehat{\AAA}$ is the $\PRTT_k$-enriched category obtained from a $\PRTTT_k$-enriched category $\widehat{\AAA}$.
Consider $P_{\PRTTT_k}(\widehat{\AAA})\in 2\PRTTT_k$ and $P_{\PRTT_k}(\widehat{\AAA})\in 2\PRTT_k$.
Then there exists a canonical equivalence $\HH^\bullet_{2\PRTT_k}(P_{\PRTT_k}(\widehat{\AAA})/k)\simeq \HH^\bullet_{2\PRTTT_k}(P_{\PRTTT_k}(\widehat{\AAA})/k)$ in $\Alg_3(\Mod_k)$.

\end{enumerate}
\end{Proposition}

\Proof
We prove (1).
By Proposition~\ref{enrichdual}, 
$\End(\underline{P(\AAA)})\in \Alg_1(2\PR_k)$
is canonically equivalent to $\End(P(\AAA))\otimes_{\PRTT_k}\PR_k$, where 
$(-)\otimes_{\PRTT_k}\PR_k$ indicates the symmetric monoidal base change functor $2\PRTT_k\to 2\PR_k$.
It gives rise to
\[
\Hom_{{}_{\mathcal{A}}\BMod_{\mathcal{A}}(\End(\Fun(X,\PRTT_k)))}(\mathcal{A},\mathcal{A})\to \Hom_{{}_{\underline{\mathcal{A}}}\BMod_{\underline{\mathcal{A}}}(\End(\Fun(X,\PR_k)))}(\underline{\mathcal{A}},\underline{\mathcal{A}})
\]
in $\Alg_2(\PR_k)$. Here $\underline{\AAA}\in \Fun(X\times X,\PR_k)$ is the $\PR_k$-enriched category determined by $\Fun(X\times X,\PRTT_k)\to \Fun(X\times X,\PR_k)$.
It is enough to prove that it induces an equivalence between endomorphism algebra objects of the unit objects (the identity functors) in $\Alg_3(\Mod_k)$.
We identify the left-hand side with
$\lim_{[n]\in \Delta}\Hom_{\Fun(X\times X,\PRTT_k)}( \mathcal{A}^{\otimes 2n+1},\mathcal{A})$. This limit is taken in $\PRTT_k$.
As for $\HH^\bullet(\underline{P(\AAA)}/k)$
we identify the right-hand side with $\lim_{[n]\in \Delta}\Hom_{\Fun(X\times X,\PR_k)}( \underline{\mathcal{A}}^{\otimes 2n+1},\underline{\mathcal{A}})$,
where the limit is taken in $\PR_k$.
We consider an intermediate objects $\lim_{[n]\in \Delta}\Hom_{\Fun(X\times X,\PRT_k)}( (\mathcal{A}')^{\otimes 2n+1},\mathcal{A}')$,
where $\AAA'\in \Fun(X\times X,\PRT_k)$ is the $\PRT_k$-enriched category determined by $\Fun(X\times X,\PRTT_k)\to \Fun(X\times X,\PRT_k)$ (induced by $\PRTT_k\hookrightarrow \PRT_k$),
and the limit is taken in $\PRT_k$. (We note that it ``computes'' $\Hom_{{}_{\mathcal{A}'}\BMod_{\mathcal{A}'}(\End(\Fun(X,\PRT_k)))}(\mathcal{A}',\mathcal{A}')$.)
For $\mathcal{P}, \mathcal{Q}\in \PRT_k$, an internal Hom object $\Hom_{\PRT_k}(\mathcal{P},\mathcal{Q})$ from $\mathcal{P}$ to $\mathcal{Q}$ in $\PRT_k$ is 
given  by the pair $(\Hom_{\PR_k}(\underline{\mathcal{P}},\underline{\mathcal{Q}}), \Hom_{\PR_k}(\underline{\mathcal{P}},\underline{\mathcal{Q}})^{\le0})$
such that $\Hom_{\PR_k}(\underline{\mathcal{P}},\underline{\mathcal{Q}})\in \PR_k$ 
is the internal Hom object from the underlying object $\underline{\mathcal{P}}$ to $\underline{\mathcal{Q}}$, and 
 $\Hom_{\PR_k}(\underline{\mathcal{P}},\underline{\mathcal{Q}})^{\le0}$ is the full subcategory, spanned by right $t$-exact functors (cf. the proof of Corollary~\ref{bimoduleexist}).
Moreover, the symmetric monoidal forgetful functor $\PRT_k\to \PR_k$ preserves small limits (see Proposition~\ref{PRTcolimitlimit}). Thus, we see that the canonical functor  
\[
\lim_{[n]\in \Delta}\Hom_{\Fun(X\times X,\PRT_k)}((\mathcal{A}')^{\otimes 2n+1},\mathcal{A}') \to \lim_{[n]\in \Delta}\Hom_{\Fun(X\times X,\PR_k)}( \underline{\mathcal{A}}^{\otimes 2n+1},\underline{\mathcal{A}})
\]
in $\PR_k$ is an equivalence.
There also exists a canonical map
\[
f:\lim_{[n]\in \Delta}\Hom_{\Fun(X\times X,\PRTT_k)}( \mathcal{A}^{\otimes 2n+1},\mathcal{A}) \to \lim_{[n]\in \Delta}\Hom_{\Fun(X\times X,\PRT_k)}((\mathcal{A}')^{\otimes 2n+1},\mathcal{A}') 
\]
in $\PRT_k$, which exhibits the left-hand side as a right-completion of the right-hand side (note that the internal Hom objects in $\PRTT_k$ are constructed 
as right completions of those in $\PRT_k$). It is enough to prove that $f$ induces an equivalence of endomorphism algebra objects of unit objects (identity functors).
We consider the diagram  
\[
\xymatrix{
\Alg_1(\Mod_k)  \ar[rr]^{\theta_{0\PRT_k,0}} \ar[rd]^{\theta_{\PRTT_k,0}} &  &  (\PRT_k)_{\Mod_k/} \ar@<0.5ex>[ld]^{\widehat{R}} \\
  & (\PRTT_k)_{\Mod_k/} \ar@<0.5ex>[ru]^{I^+}&  
}
\]
Both $\theta_{0\PRT_k,0}$ and $\theta_{\PRTT_k,0}$ send $A$ to $((\LMod_A,\LMod_A^{\le0}), A\in \LMod_A)$.
The functor $I^+$ is the inclusion, while $\widehat{R}$ is the right completion functor.
Note that $I^+ \circ \theta_{\PRTT_k,0}=\theta_{\PRTT_k,0}$, and the right adjoint functors of $\theta_{0\PRT_k,0}$ and $\theta_{\PRTT_k,0}$ are given informally by
$(\mathcal{L},L)\mapsto \End(L)^{op}$. From the right completion $f$ and this commutative diagram, 
we conclude that $f$ induces an equivalence of endomorphism algebra objects of unit objects.

We prove (2).
As in (1),
$\End(P_{\PRTTT_k}(\widehat{\AAA}))\in \Alg_1(2\PRTTT_k)$
is canonically equivalent to $\End(P_{\PRTT_k}(\AAA))\otimes_{\PRTT_k}\PRTTT_k$, where 
$(-)\otimes_{\PRTT_k}\PRTTT_k$ indicates the symmetric monoidal base change functor $2\PRTT_k\to 2\PRTTT_k$.
It gives rise to 
\[
\Hom_{{}_{\widehat{\mathcal{A}}}\BMod_{\widehat{\mathcal{A}}}(\End(\Fun(X,\PRTT_k)))}(\widehat{\mathcal{A}},\widehat{\mathcal{A}})\to \Hom_{{}_{\widehat{\mathcal{A}}}\BMod_{\widehat{\mathcal{A}}}(\End(\Fun(X,\PRTTT_k)))}(\widehat{\mathcal{A}},\widehat{\mathcal{A}})
\]
in $\Alg_2(\PRTT_k)$.
From the proof of Corollary~\ref{bimoduleexist}, it can be identified with the canonical morphism
\[
f:\lim_{[n]\in \Delta}\Hom_{\Fun(X\times X,\PRTT_k)}( \widehat{\mathcal{A}}^{\otimes 2n+1},\widehat{\mathcal{A}})\to \lim_{[n]\in \Delta}\Hom_{\Fun(X\times X,\PRTTT_k)}( \widehat{\mathcal{A}}^{\widehat{\otimes} 2n+1},\widehat{\mathcal{A}})
\]
in $\PRTT_k$. Here $\widehat{\mathcal{A}}^{\otimes 2n+1}$ is the tensor product in $\Fun(X\times X,\PRTT_k)$, while
$\widehat{\mathcal{A}}^{\widehat{\otimes} 2n+1}$ is the tensor product in $\Fun(X\times X,\PRTTT_k)$.
It is enough to prove that $f$ is an equivalence.
Each morphism
\[
\Hom_{\Fun(X\times X,\PRTT_k)}( \widehat{\mathcal{A}}^{\otimes 2n+1},\widehat{\mathcal{A}})\to \Hom_{\Fun(X\times X,\PRTTT_k)}( \widehat{\mathcal{A}}^{\widehat{\otimes} 2n+1},\widehat{\mathcal{A}})
\]
is an equivalence since it is induced by the left completion, and 
$\widehat{\mathcal{A}}^{\widehat{\otimes} 2n+1}$ is a left completion of $\widehat{\mathcal{A}}^{\otimes 2n+1}$.
Thus, $f$ is an equivalence. 
\QED

\subsection{Categorified sheaves}

\label{catsheaves}

Let $\XX$ be a perfect stack over $k$ in the sense of \cite{BFN} (see {\it loc. cit.}).
Typical examples include quasi-compact and separated derived schemes over $k$.
We consider the``categorified sheaves" on $\XX$, that is, sheaves of categories.
We briefly describe them.

For simplicity, we suppose that $\mathcal{X}$ is a quasi-compact and separated derived scheme for a moment.
Let $\QC(\mathcal{X})$ be the $k$-linear symmetric monoidal $\infty$-category of quasi-coherent complexes on $\mathcal{X}$ (cf. \cite[VIII]{DAG}, \cite{BFN}).
If we write $\textup{AffOp}_{\XX}$ for the poset of affine open subschemes, then 
$\QC(\XX)=\lim_{\Spec R(\mathcal{U})\in \textup{AffOp}_{\XX}}\Mod_{R(\mathcal{U})}$ where $i^*:\Mod_{R(\mathcal{U})}\to \Mod_{R(\mathcal{V})}$
is base change functor along $R(\mathcal{U})\to R(\mathcal{V})$
for each $i:\mathcal{V}\to \mathcal{U}$. By definition, if $\XX=\Spec R$, then $\QC(\XX)=\Mod_R$.
We consider the functor 
\[
2\QC:\textup{AffOp}_{\XX}^{op} \to 2\PR_k
\]
which carries $\mathcal{U}=\Spec R(\mathcal{U})$ to $\Mod_{\Mod_{R(\mathcal{U})}}(\PR_k)$
and sends a morphism $i:\mathcal{V}\to \mathcal{U}$ to $\Mod_{\Mod_{R(\mathcal{U})}}(\PR_k)\to \Mod_{\Mod_{R(\mathcal{V})}}(\PR_k)$
given by $\mathcal{D}\mapsto \mathcal{D}\otimes_{\Mod_{R(\mathcal{U})}}\Mod_{R(\mathcal{V})}$ (i.e., 
the base change along $i^*:\Mod_{R(\mathcal{U})}\to \Mod_{R(\mathcal{V})}$). We define the $\infty$-category of categorified sheaves on $\XX$ as
\[
\textup{QC}_{\textup{2-cat}}(\XX):=\lim_{\mathcal{U}\in \textup{AffOp}_{\XX}}2\QC(\mathcal{U})=\lim_{\mathcal{U}\in \textup{AffOp}_{\XX}}\Mod_{\Mod_{R(\mathcal{U})}}(\PR_k).
\]
Informally, an object $\mathcal{F}\in \textup{QC}_{\textup{2-cat}}(\XX)$ is an assignment 
$\Gamma(\mathcal{U}, \mathcal{F})\in\Mod_{\Mod_{R(\mathcal{U})}}(\PR_k) $ to each $\mathcal{U}\in \textup{AffOp}_{\XX}$
and an equivalence
\[
\Gamma(\mathcal{U}, \mathcal{F})\otimes_{\Mod_{R(\mathcal{U})}}\Mod_{R(\mathcal{V})}\simeq \Gamma(\mathcal{V}, \mathcal{F})
\]
to each $i:\mathcal{V}\to\mathcal{U}$, endowed with homotopy coherent compatibility data.

Each $\Mod_{\Mod_{R(\mathcal{U})}}(\PR_k)$ admits a base change functor $\Mod_{\QC(\XX)}(\PR_k)\to \Mod_{\Mod_{R(\mathcal{U})}}(\PR_k)$
along the pullback functor $\QC(\XX)\to \Mod_{R(\mathcal{U})}$.
They determine a functor
\[
\Mod_{\QC(\XX)}(\PR_k)\to \lim_{\mathcal{U}\in \textup{AffOp}_{\XX}}\Mod_{\Mod_{R(\mathcal{U})}}(\PR_k)=\textup{QC}_{\textup{2-cat}}(\XX).
\]
By a $1$-affineness of $\XX$ (cf. \cite[Theorem 2.1.1]{Gai1aff}), this functor is an equivalence.
Thus, we may identify $\Mod_{\QC(\XX)}(\PR_k)$ with $\textup{QC}_{\textup{2-cat}}(\XX)$.

The stable $\infty$-category  $\QC(\XX)$  has the ``standard" accessible left and right complete $t$-structure (see \cite[VIII, Proposition 2.3.16, Proposition 2.3.18]{DAG}).
An object $M\in \QC(\XX)$ is connective on this $t$-structure if and only if for any $\mathcal{U}\in\textup{AffOp}_{\XX}$ the restriction $M|_{\mathcal{U}}\in \Mod_{R(\mathcal{U})}$
is connective. 
We regard $\QC(\mathcal{X})$ as an object of $\Alg_1(\PRTT_k)$ that comes from $\CAlg(\PRTT_k)$.
One can consider the $t$-structured categorified sheaves.
We define
\[
\textup{QC}_{\textup{2-cat,t}}(\XX):=\lim_{\mathcal{U}\in \textup{AffOp}_{\XX}}\Mod_{\Mod_{R(\mathcal{U})}}(\PRT_k).
\]
Then, by the above $1$-affineness of $\XX$, one sees that a canonical equivalence 
$\Mod_{\QC(\XX)}(\PRT_k)\simeq \textup{QC}_{\textup{2-cat,t}}(\XX)$ (with some elaboration we also have a version of $\PRTT_k$). 
Therefore, we will adopt $\Mod_{\QC(\XX)}(\PRT_k)$ as the $t$-structured enhancement of  $\Mod_{\QC(\XX)}(\PR_k)$.

By \cite[VIII, Proposition 2.3.18]{DAG}, the $t$-structure on $\QC(\XX)$ is both left and right complete. 
We think of the monoidal $\infty$-category $\QC(\mathcal{X})$ as a $\PRTTT_k$-enriched category having a single object,
which we usually denote by $B\QC(\mathcal{X})$.
In this case, $P(B\QC(\mathcal{X}))=\Mod_{\QC(\mathcal{X})}(\PRTTT_k)$. Note that $\Mod_{\QC(\mathcal{X})}(\PRTTT_k)$ belongs to $2\PRTTT_k$.
We adopt this as the complete $t$-structured enhancement of  $\Mod_{\QC(\XX)}(\PR_k)$.
Note that the underlying stable $(\infty,2)$-category
$\Mod_{\QC(\mathcal{X})}(\PRTTT_k)\otimes_{\PRTTT_k}\PR_k$ is equivalent to $\textup{QC}_{\textup{2-cat}}(\XX)$.

\begin{Proposition}
\label{catsheafhoch}
Let $\XX$ be a perfect stack over $k$.
Let $\HH^\bullet(\Mod_{\QC(\mathcal{X})}(\PRTTT_k)/k)$ be the Hochschild cohomology of $\Mod_{\QC(\mathcal{X})}(\PRTTT_k)$ (cf. Proposition~\ref{enrichhoch}).
Then the underlying complex of Hochschild cohomology $\HH^\bullet(\Mod_{\QC(\mathcal{X})}(\PRTTT_k)/k)$
is equivalent to the Hom complex $\mathcal{H}om_{\QC(L\XX)}(\OO_{\Delta},\OO_\Delta)\in \Mod_k$
where $L\XX$ is the loop space of $\XX$, defined as the fiber product $\XX\times_{\XX\times_k\XX}\XX$.
The object $\OO_{\Delta}\in \QC(L\XX)$ is the pushforward of the structure sheaf $\OO_{\mathcal{X}}\in \QC(X)$
along the diagonal morphism $\XX\to \XX\times_{\XX\times_k\XX}\XX=L\XX$.

Suppose that $k$ is of characteristic zero, and $\XX$ is an ordinary scheme separated and of finite type.
There exists an equivalence
\[
\HH^\bullet(\Mod_{\QC(\mathcal{X})}(\PRTTT_k)/k)\simeq \mathcal{H}om_{\QC(L\XX)}(\OO_{\Delta},\OO_\Delta)\simeq  \Gamma(\XX,\Sym_{\OO_{\XX}}^{\bullet}(\mathbb{T}_{\XX/k}[-2]))
\]
in $\Mod_k$
where $\mathbb{T}_{\XX/k}$ is the tangent complex of $\XX$, defined as the (derived) dual $\mathcal{H}om_{\QC(\XX)}(\mathbb{L}_{\XX/k},\OO_{\XX})$
of the cotangent complex $\mathbb{L}_{\XX/k}$.
Here $\Gamma(\XX,-)$ means the (derived) global sections, and $\Sym_{\OO_{\XX}}^{\bullet}(-)$ means a free commutative algebra in $\QC(\XX)$, that is, the symmetric algebra.
\end{Proposition}

\Proof
By Proposition~\ref{comparison1},
there exist equivalences 
\[
\HH^\bullet(\Mod_{\QC(\mathcal{X})}(\PRTTT_k)/k)\simeq\HH^\bullet(\Mod_{\QC(\mathcal{X})}(\PRTT_k)/k)\simeq \HH^\bullet(\Mod_{\QC(\mathcal{X})}(\PR_k)/k).
\]
We may and will replace $\HH^\bullet(\Mod_{\QC(\mathcal{X})}(\PRTTT_k)/k)$ with $\HH^\bullet(\Mod_{\QC(\mathcal{X})}(\PR_k)/k)$ in the statement.
As observed in Section~\ref{enrichsec},
there exist canonical equivalences
\begin{eqnarray*}
\Hom_{\PR_k}(\Mod_{\QC(\mathcal{X})}(\PR_k),\Mod_{\QC(\mathcal{X})}(\PR_k)) &\simeq& \Mod_{\QC(\mathcal{X})}(\PR_k)\otimes_{\PR_k}\Mod_{\QC(\mathcal{X})}(\PR_k)  \\
&\simeq& \Mod_{\QC(\mathcal{X})\otimes \QC(\mathcal{X})}(\PR_k).
\end{eqnarray*}
Since there exists an equivalence $p_1\boxtimes p_2:\QC(\mathcal{X})\otimes \QC(\mathcal{X})\to \QC(\mathcal{X}\times_k \mathcal{X})$ induced by the pullbacks of projections
$p_1,p_2:\mathcal{X}\times_k \mathcal{X}\rightrightarrows \mathcal{X}$ (see \cite[Theorem 4.7]{BFN}), we have the equivalence
\[
\Mod_{\QC(\mathcal{X})\otimes \QC(\mathcal{X})}(\PR_k)\simeq \Mod_{\QC(\mathcal{X}\times_k\mathcal{X})}(\PR_k).
\]
The identity functor $\Mod_{\QC(\mathcal{X})}(\PR_k)\to \Mod_{\QC(\mathcal{X})}(\PR_k)$ corresponds to
$\QC(\mathcal{X})$ endowed with the diagonal $\QC(\mathcal{X})$-$\QC(\mathcal{X})$-bimodule structure.
This bimodule is a $\QC(\mathcal{X}\times_k\mathcal{X})$-module $\QC(\mathcal{X})$ induced by the symmetric monoidal pullback functor
$\Delta^*:\QC(\mathcal{X}\times_k\mathcal{X})\to \QC(\mathcal{X})$
along the diagonal morphism $\Delta:\mathcal{X}\to \mathcal{X}\times_k\mathcal{X}$.
Moreover, $\QC(\mathcal{X})$ is dualizable in $\PR_k$ whose dual is given by $\QC(\mathcal{X})$ itself.
Indeed, according to \cite[Corollary 4.8]{BFN} $\QC(\XX)$ is self-dual in $\Mod_{\QC(\XX\times_k\XX)}(\PR_k)$. 
Therefore, we have
\[
\Hom_{\Mod_{\QC(\mathcal{X}\times_k \mathcal{X})}(\PR_k)}(\QC(\mathcal{X}),\QC(\mathcal{X}))\simeq \QC(\mathcal{X})\otimes_{\QC(\mathcal{X}\times_k \mathcal{X})}\QC(\mathcal{X})\simeq \QC(\mathcal{X}\times_{\mathcal{X}\times_k\mathcal{X}}\mathcal{X})
\]
in $\PR_k$, see \cite[Theorem 4.7]{BFN} for the latter equivalence.
Through the equivalence,
the identity functor $\QC(\XX)\to \QC(\XX)$ corresponds to the $\OO_{\Delta}\in \QC(L\XX)$, which is the pushforward of the structure sheaf $\OO_{\mathcal{X}}\in \QC(\XX)$
along the diagonal morphism $\Delta_L:\XX\to \XX\times_{\XX\times_k\XX}\XX=L\XX$.
Consequently, Hochschild cohomology is equivalent to $\mathcal{H}om_{\QC(L\XX)}(\OO_{\Delta},\OO_\Delta)\in \Mod_k$ (the equivalence can be promoted to the level of $\Alg_1(\Mod_k)$).

We show the second assertion.
Under our assumption,  Koszul duality \cite{MR} and \cite{H} imply that
there exists an equivalence
\[
\mathcal{H}om_{\QC(L\XX)}(\OO_{\Delta},\OO_\Delta)\simeq  \Gamma(\XX,\Sym_{\OO_{\XX}}^{\bullet}(\mathbb{T}_{\XX/k}[-2]))
\]
in $\Mod_k$.
\QED

\subsection{Higher module category}

Suppose that $\mathsf{P}$ is any of $\PR_k$ or $\PRTT_k$ or $\PRTTT_k$.
Let $\mathcal{A}$ be an object of $\CAlg(\mathsf{P})$.
We will associate to $\AAA$ a kind of higher module category.
We write $0\Mod_{\AAA}[\mathsf{P}]=\AAA$ and $1\Mod_{\AAA}[\mathsf{P}]=\Mod_{\AAA}(\mathsf{P})$.  We note that $0\Mod_{\AAA}\in \mathsf{P}$ and $1\Mod_{\AAA}\in 2\mathsf{P}$. 
We define $n\Mod_{\AAA}[\mathsf{P}]\in (n+1)\mathsf{P}$ as $\Mod_{(n-1)\Mod_{\AAA}[\mathsf{P}]}(n\mathsf{P})$.
For ease of notation we often write $n\Mod_{\AAA}$ for $n\Mod_{\AAA}[\mathsf{P}]$ if the setting is clear.

\begin{Proposition}
\label{highermodulehoch}
If $\mathsf{P}$ is any of $\PR_k$ or $\PRTT_k$ or $\PRTTT_k$,
then $n\Mod_{\AAA}[\mathsf{P}]\in (n+1)\mathsf{P}$ has Hochschild cohomology $\HH^\bullet(n\Mod_{\AAA}/k)\in \Alg_{n+2}(\Mod_k)$.
\end{Proposition}

We show this proposition step by step.

\begin{Lemma}
\label{algebramonoidal}
Let $i$ be a natural number.
Let $\BBB$ be an object of $\CAlg(i\mathsf{P})$.
Let $m$ be a nonnegative integer.
Then there exists an equivalence
$\Mod_{\BBB}(i\mathsf{P})\otimes S^m\to\Mod_{\BBB\otimes S^m}(i\mathsf{P})$ in $\CAlg(i\mathsf{P})$. 
Here $\Mod_{\BBB}(i\mathsf{P})\otimes S^m$ denotes the tensor with $S^m$ in $\CAlg((i+1)\mathsf{P})$,
and $\BBB\otimes S^m$ denotes the tensor with the sphere $S^m$ in $\CAlg(i\mathsf{P})$.
\end{Lemma}

\Proof
The proof proceeds by induction on $m$.
If $m=0$, there are canonical equivalences
\[
\Mod_{\BBB}(i\mathsf{P})\otimes S^0\simeq \Mod_{\BBB}(i\mathsf{P})\otimes\Mod_{\BBB}(i\mathsf{P})\simeq \Mod_{\BBB\otimes \BBB}(i\mathsf{P})
\]
where the middle is the coproduct of $\Mod_{\BBB}(i\mathsf{P})$ in $\CAlg((i+1)\mathsf{P})$ as well as the underlying object in $(i+1)\mathsf{P}$ is the tensor product.
The right equivalence follows from \cite[4.8.5.16]{HA}: The right map is defined by the base changes along  $\mathsf{1}\otimes \BBB\to \BBB\otimes \BBB\leftarrow \BBB\otimes\mathsf{1}$.
($\mathsf{1}$ is the unit algebra of $i\mathsf{P}$).  Thus, the case of $m=0$ holds.
Next we assume that $\Mod_{\BBB}(i\mathsf{P})\otimes S^m\simeq \Mod_{\BBB\otimes S^m}(i\mathsf{P})$.
More precisely, there is a canonically defined functor $\Mod_{\BBB}(i\mathsf{P})\otimes S^m\to \Mod_{\BBB\otimes S^m}(i\mathsf{P})$, which is assumed to be an equivalence: If we have a canonically defined equivalence
$\Mod_{\BBB}(i\mathsf{P})\otimes S^{m-1}\simeq \Mod_{\BBB\otimes S^{m-1}}(i\mathsf{P})$, then it gives rise to  
\begin{eqnarray*}
\Mod_{\BBB}(i\mathsf{P})\otimes S^m\simeq \Mod_{\BBB}(i\mathsf{P})\otimes_{\Mod_{\BBB}(i\mathsf{P})\otimes S^{m-1}}\Mod_{\BBB}(i\mathsf{P}) &\simeq& \Mod_{\BBB}(i\mathsf{P})\otimes_{\Mod_{\BBB\otimes S^{m-1}}(i\mathsf{P})}\Mod_{\BBB}(i\mathsf{P})   \\
&\to& \Mod_{\BBB\otimes S^m}(i\mathsf{P})
\end{eqnarray*}
where the last map is induced by the canonical morphism from
the pushout of $\BBB\leftarrow \BBB\otimes S^{m-1}\to \BBB$ to $\BBB\otimes S^m$.
We will prove that the canonically defined morphism 
\begin{eqnarray*}
\Mod_{\BBB}(i\mathsf{P})\otimes S^{m+1}\simeq \Mod_{\BBB}(i\mathsf{P})\otimes_{\Mod_{\BBB}(i\mathsf{P})\otimes S^{m}}\Mod_{\BBB}(i\mathsf{P}) &\simeq& \Mod_{\BBB}(i\mathsf{P})\otimes_{\Mod_{\BBB\otimes S^{m}}(i\mathsf{P})}\Mod_{\BBB}(i\mathsf{P})  \\ &\stackrel{f}{\to}& \Mod_{\BBB\otimes S^{m+1}}(i\mathsf{P})
\end{eqnarray*}
is an equivalence.
To this end, we will show that the arrow $f$ is an equivalence.
We regard $\BBB$ as a diagonal $\BBB\otimes S^m$-module through $\BBB\otimes S^m\to \BBB\otimes \ast\simeq \BBB$ induced by the contraction $S^m\to\ast$.
We then see that $f$ is decomposed as
\begin{eqnarray*}
\Mod_{\BBB}(i\mathsf{P})\otimes_{\Mod_{\BBB\otimes S^{m}}(i\mathsf{P})}\Mod_{\BBB}(i\mathsf{P}) &\simeq& \Mod_{\BBB}(\Mod_{\BBB\otimes S^{m}}(i\mathsf{P}))\otimes_{\Mod_{\BBB\otimes S^{m}}(i\mathsf{P})}\Mod_{\BBB}(\Mod_{\BBB\otimes S^{m}}(i\mathsf{P})) \\
&\simeq& \Mod_{\BBB\otimes_{\BBB\otimes S^{m}}\BBB}(\Mod_{\BBB\otimes S^{m}}(i\mathsf{P})) \\
&\simeq& \Mod_{\BBB\otimes_{\BBB\otimes S^{m}}\BBB}(i\mathsf{P}) \\
&\simeq& \Mod_{\BBB\otimes(\ast\sqcup_{S^m}\ast)}(i\mathsf{P}) \\
&\simeq& \Mod_{\BBB\otimes S^{{m+1}}}(i\mathsf{P})
\end{eqnarray*}
where the first map and the third map are induced by the forgetful functor $\Mod_{\BBB\otimes S^{m}}(i\mathsf{P})\to i\mathsf{P}$ (see \cite[7.1.3.1]{HA}),
and the second equivalence follows from the symmetric monoidal structure on $\Mod_{\BBB\otimes S^{m}}(i\mathsf{P})$ (see \cite[4.8.5.16]{HA}).
The fourth equivalence follows from  the fact that the pushout $\BBB\leftarrow \BBB\otimes S^m\to \BBB$ in $\CAlg(i\mathsf{P})$ is 
naturally equivalent to the bar construction $\BBB\otimes_{\BBB\otimes S^m}\BBB$ in $i\mathsf{P}$. 
Thus, $f$ is an equivalence.
This completes the proof.
\QED

There exists an endomorphism algebra
$\End(n\Mod_{\AAA}[\mathsf{P}])\in \Alg_1((n+1)\mathsf{P})$ whose underlying object $\Hom(n\Mod_{\AAA}[\mathsf{P}],n\Mod_{\AAA}[\mathsf{P}])$
is equivalent to $\Mod_{(n-1)\Mod_{\AAA}}(n\mathsf{P})\otimes \Mod_{(n-1)\Mod_{\AAA}}(n\mathsf{P})$
since $\AAA$ is a commutative algebra object, and $n\Mod_{\AAA}[\mathsf{P}]=\Mod_{(n-1)\Mod_{\AAA}[\mathsf{P}]}(n\mathsf{P})$ is self-dual in $(n+1)\mathsf{P}$.
Here $\otimes$ means the tensor product in $(n+1)\mathsf{P}$.
According to \cite[4.8.5.16]{HA}, there exists a canonical equivalence
\[
\Mod_{(n-1)\Mod_{\AAA}}(n\mathsf{P})\otimes \Mod_{(n-1)\Mod_{\AAA}}(n\mathsf{P})\simeq \Mod_{(n-1)\Mod_{\AAA}\otimes (n-1)\Mod_{\AAA}}(n\mathsf{P})
\]
in $(n+1)\mathsf{P}$.

\begin{Lemma}
Let
$(n-1)\Mod_{\AAA}[\mathsf{P}]\in \Mod_{(n-1)\Mod_{\AAA}\otimes (n-1)\Mod_{\AAA}}(n\mathsf{P})$
be the $(n-1)\Mod_{\AAA}\otimes (n-1)\Mod_{\AAA}$-module induced by the multiplication map $(n-1)\Mod_{\AAA}\otimes (n-1)\Mod_{\AAA}\to (n-1)\Mod_{\AAA}$.
This corresponds to the identity functor in $\End(n\Mod_{\AAA})$.
Let 
\[
K_n:=\Hom_{(n-1)\Mod_{\AAA}\otimes (n-1)\Mod_{\AAA}}((n-1)\Mod_{\AAA}, (n-1)\Mod_{\AAA}) \in n\mathsf{P}
\]
be the Hom object from $(n-1)\Mod_{\AAA}\in \Mod_{(n-1)\Mod_{\AAA}\otimes (n-1)\Mod_{\AAA}}(n\mathsf{P})$ to
itself, where the subscript $(n-1)\Mod_{\AAA}\otimes (n-1)\Mod_{\AAA}$ indicates the $(n-1)\Mod_{\AAA}\otimes (n-1)\Mod_{\AAA}$-linear.
Then there exists a canonical equivalence
\[
K_n\simeq \Mod_{(n-2)\Mod_{\AAA}\otimes S^1}((n-1)\mathsf{P})
\]
in $n\mathsf{P}$.
Here $(n-2)\Mod_{\AAA}\otimes S^1$ means the tensor with the sphere $S^1\in \SSS$ in $\CAlg((n-1)\mathsf{P})$.
\end{Lemma}

\Proof
We note that
\[
(n-1)\Mod_{\AAA}=\Mod_{(n-2)\Mod_{\AAA}}((n-1)\mathsf{P})\simeq \Mod_{(n-2)\Mod_{\AAA}}(\Mod_{(n-2)\Mod_{\AAA}\otimes (n-2)\Mod_{\AAA}}((n-1)\mathsf{P})).
\]
The right equivalence is obtained by regarding $(n-2)\Mod_{\AAA}$ as the diagonal $(n-2)\Mod_{\AAA}\otimes (n-2)\Mod_{\AAA}$-module (cf. \cite[7.1.3.1]{HA}).
We put 
\[
(n-1)\Mod_{\AAA}^{\otimes 2}:=(n-1)\Mod_{\AAA}\otimes (n-1)\Mod_{\AAA}\simeq \Mod_{(n-2)\Mod_{\AAA}\otimes (n-2)\Mod_{\AAA}}((n-1)\mathsf{P}).
\]
Here Lemma~\ref {algebramonoidal} implies the right equivalence. 
By the canonical equivalence
$\Mod_{(n-2)\Mod_{\AAA}}((n-1)\Mod_{\AAA}^{\otimes2})\simeq \Mod_{(n-2)\Mod_{\AAA}}((n-1)\mathsf{P})$ (again \cite[7.1.3.1]{HA})
and the self-duality of $\Mod_{(n-2)\Mod_{\AAA}}((n-1)\Mod_{\AAA}^{\otimes2})$, there exist equivalences
\begin{eqnarray*}
K_n &=&\Hom_{(n-1)\Mod_{\AAA}^{\otimes2}}(\Mod_{(n-2)\Mod_{\AAA}}((n-1)\mathsf{P}), \Mod_{(n-2)\Mod_{\AAA}}((n-1)\mathsf{P})) \\
&\simeq& \Hom_{(n-1)\Mod_{\AAA}^{\otimes2}}(\Mod_{(n-2)\Mod_{\AAA}}((n-1)\Mod_{\AAA}^{\otimes 2}), \Mod_{(n-2)\Mod_{\AAA}}((n-1)\Mod_{\AAA}^{\otimes 2})) \\
&\simeq& \Mod_{(n-2)\Mod_{\AAA}}((n-1)\Mod_{\AAA}^{\otimes2})\otimes_{(n-1)\Mod_{\AAA}^{\otimes2}}\Mod_{(n-2)\Mod_{\AAA}}((n-1)\Mod_{\AAA}^{\otimes2}) \\
&\simeq& \Mod_{(n-2)\Mod_{\AAA}\otimes(n-2)\Mod_{\AAA}}((n-1)\Mod_{\AAA}^{\otimes 2}) \\
&\simeq& \Mod_{(n-2)\Mod_{\AAA}\otimes_{(n-2)\Mod_{\AAA}\otimes (n-2)\Mod_{\AAA}}(n-2)\Mod_{\AAA}}((n-1)\mathsf{P}).
\end{eqnarray*}
We note that
there exists an equivalence 
\[
(n-2)\Mod_{\AAA}\otimes_{(n-2)\Mod_{\AAA}\otimes (n-2)\Mod_{\AAA}}(n-2)\Mod_{\AAA}\simeq (n-2)\Mod_{\AAA}\otimes S^1
\]
in 
$\CAlg((n-1)\mathsf{P})$,
where the right-hand side means the tensor with the circle $S^1\in \SSS$ in $\CAlg((n-1)\mathsf{P})$.
Thus, we have an equivalence
\[
K_n\simeq \Mod_{(n-2)\Mod_{\AAA}\otimes S^1}((n-1)\mathsf{P})
\]
in $n\mathsf{P}$.
\QED

We take the next step.
Let 
\[
K_{n-1}=\Hom_{(n-2)\Mod_{\AAA}\otimes S^1}((n-2)\Mod_{\AAA}, (n-2)\Mod_{\AAA}) \in (n-1)\mathsf{P}
\]
be the Hom object from $(n-2)\Mod_{\AAA}\in \Mod_{(n-2)\Mod_{\AAA}\otimes S^1}((n-1)\mathsf{P})$ to
itself, where the subscript $(n-2)\Mod_{\AAA}\otimes S^1$  indicates the $(n-2)\Mod_{\AAA}\otimes S^1$-linear.
Using Lemma~\ref{algebramonoidal}, the same argument shows that there exists an equivalence
\[
K_{n-1}\simeq \Mod_{(n-3)\Mod_{\AAA}\otimes_{(n-3)\Mod_{\AAA}\otimes S^1}(n-3)\Mod_{\AAA}}((n-2)\mathsf{P})
\]
in $(n-1)\mathsf{P}$. Here the contraction map $S^1\to \ast$ to the contractible space $\ast$
induces a morphism $(n-3)\Mod_{\AAA}\otimes S^1\to (n-3)\Mod_{\AAA}$ in $\CAlg((n-1)\mathsf{P})$,
which exhibits $(n-3)\Mod_{\AAA}$ as a $(n-3)\Mod_{\AAA}\otimes S^1$-module.
Since the pushout $\ast\sqcup_{S^1}\ast$ is $S^2$, we have $(n-3)\Mod_{\AAA}\otimes_{(n-3)\Mod_{\AAA}\otimes S^1}(n-3)\Mod_{\AAA}\simeq (n-3)\Mod_{\AAA}\otimes S^2$.
Consequently, we have the following.
\begin{Lemma}
There exists a canonical equivalence $K_{n-1}\simeq \Mod_{(n-3)\Mod_{\AAA}[\mathsf{P}]\otimes S^2}((n-2)\mathsf{P})$ in $(n-1)\mathsf{P}$.
\end{Lemma}

Repeating this procedure, we have

\begin{Lemma}
\label{kirep}
If we 
put $K_{i+1}=\Hom_{i\Mod_{\AAA}\otimes S^{n-i-1}}(i\Mod_{\AAA}, i\Mod_{\AAA})\in \Alg_{n-i+1}(i\mathsf{P})$, there exists a canonical equivalence
\[
K_{i+1}\simeq \Mod_{(i-1)\Mod_{\AAA}\otimes S^{n-i}}(i\mathsf{P}).
\]
\end{Lemma}

We arrive at 
\begin{eqnarray*}
K_2=\Hom_{\Mod_{\AAA}\otimes S^{n-2}}(\Mod_{\AAA}(\mathsf{P}), \Mod_{\AAA}(\mathsf{P})) &\simeq& \Mod_{\AAA}(\mathsf{P})\otimes_{\Mod_{\AAA}\otimes S^{n-2}}\Mod_{\AAA}(\mathsf{P})  \\
&\simeq& \Mod_{\AAA}\otimes S^{n-1}
\end{eqnarray*}
Using Lemma~\ref{algebramonoidal}, we see $K_2\simeq \Mod_{\AAA\otimes S^{n-1}}$.

Next, we let $\AAA_{\Delta}$ denote the $\AAA\otimes S^{n-1}$-module
$\AAA$ whose $\AAA\otimes S^{n-1}$-module structure is determined by
the symmetric monoidal functor $\AAA\otimes S^{n-1}\to \AAA$ induced by $S^{n-1}\to\ast$.
This diagonal module $\AAA_{\Delta}$ corresponds to the identity functor.
Next, we will show

\begin{Lemma}
\label{higherendo}
Suppose that $\mathsf{P}$ is either $\PRTT_k$ or $\PRTTT_k$.
Then there exists a Hom object $\Hom_{\AAA\otimes S^{n-1}}(\AAA_{\Delta},\AAA_{\Delta})$ in $\PRTT_k$ (the subscript again indicates $\AAA\otimes S^{n-1}$-linear).
In particular, there exists an endomorphism algebra $\End_{\AAA\otimes S^{n-1}}(\AAA_{\Delta})$ in $\Alg_{n+1}(\PRTT_k)$.

Suppose that $\mathsf{P}$ is $\PR_k$. 
Then there exists a Hom object $\Hom_{\AAA\otimes S^{n-1}}(\AAA_{\Delta},\AAA_{\Delta})$ in $\PR_k$,
and there exists an endomorphism algebra $\End_{\AAA\otimes S^{n-1}}(\AAA_{\Delta})$ in $\Alg_{n+1}(\PR_k)$.
\end{Lemma}

\Proof
As in the case of $t$-structured $(\infty,2)$-categories (see Section~\ref{enrichsec}),  from the free-forget adjunction of modules
\[
\xymatrix{
\PRTT_k \ar@<0.5ex>[r]^{(-)\otimes(\AAA\otimes S^{n-1})}  &  \Mod_{\AAA\otimes S^{n-1}}(\PRTT_k)   \ar@<0.5ex>[l]^{\textup{forget}}  
}
\]
we obtain a simplicial object $\mathcal{R}_\bullet$ in $\Mod_{\AAA\otimes S^{n-1}}(\PRTT_k)$ with an augmentation $\epsilon:\mathcal{R}_\bullet\to \mathcal{A}_{\Delta}$ 
such that each $m$-th term is a free module $\mathcal{E}(m)\otimes (\mathcal{A}\otimes S^{n-1})$,
and $\epsilon$ exhibits $\mathcal{A}_{\Delta}$ as a colimit.
We then have
\[
\Hom_{\AAA\otimes S^{n-1}}(\AAA_{\Delta},\AAA_{\Delta})\simeq \lim_{[m]\in \Delta}\Hom(\mathcal{E}(m),\AAA)
\]
where each $\Hom(\mathcal{E}(m),\AAA)$ is a Hom object, defined as an object of $\PRTT_k$,
from $\mathcal{E}(m)$ to $\AAA$. Since $\PRTT_k$ and $\PRTTT_k$ admit Hom objects defined in $\PRTT_k$ and have small limits,
it follows that $\Hom_{\AAA\otimes S^{n-1}}(\AAA_{\Delta},\AAA_{\Delta})$ exists in $\PRTT_k$.

This argument also works when $\mathsf{P}=\PR_k$.
\QED

\begin{Corollary}
\label{goodendohighermodule}
If $\mathsf{P}$ is either $\PRTT_k$ or $\PR_k$,
$\End(n\Mod_{\AAA}[\mathsf{P}])$ has a good endomorphism algebra of a unit.
\end{Corollary}

\Proof
Apply Lemma~\ref{kirep} and Lemma~\ref{higherendo}.
\QED

{\it Proof of Proposition~\ref{highermodulehoch}.}
We first suppose that $\mathsf{P}$ is either $\PRTT_k$ or $\PRTTT_k$.
As in Section~\ref{enrichsec},
applying $\Alg_{n+1}(-)$ to the right adjoint $E^+_{\PRTT_k}:(\PRTT_k)_{\Mod_k/}\to \Alg_1(\Mod_k)$, we have
$\Alg_{n+1}(E^+_{\PRTT_k}):\Alg_{n+1}(\PRTT_k)\to \Alg_{n+2}(\Mod_k)$.
The Hochschild cohomology $\HH^\bullet(n\Mod_{\AAA}[\mathsf{P}]/k)\in \Alg_{n+2}(\Mod_k)$ is constructed as the image of $\End_{\AAA\otimes S^{n-1}}(\AAA_{\Delta})$  under
$\Alg_{n+2}(E^+_{\PRTT_k})$. By construction, this has the required property (when $\mathsf{P}=\PRTTT_k$, note that $\PRTTT_k$ is a localization of $\PRTT_k$).

Next, suppose that $\mathsf{P}=\PR_k$.
Using the adjoint pair $\theta_{0\PR_k,0}:\Alg_1(\Mod_k)\rightleftarrows (\PR_k)_{\Mod_k/}:E_{\PR_k}^+$ instead of $(\theta_{0\PRTT_k,0},E_{\PRTT_k}^+)$,
one can construct Hochschild cohomology in the same way as above.
\QED

\begin{Proposition}
\label{comparison2}
The following holds.
\begin{enumerate}
\item
For $\mathcal{A}\in \CAlg(\PRTT_k)$, we
let $\underline{\AAA}$ denote the image of $\AAA$ under the forgetful functor $\CAlg(\PRTT_k)\to \CAlg(\PR_k)$.
We obtain $n\Mod_{\underline{\AAA}}[\PR_k]\in (n+1)\PR_k$.
We note that $n\Mod_{\underline{\AAA}}[\PR_k]$ is canonically equivalent to $n\Mod_{\AAA}[\PRTT_k]\otimes_{n\PRTT_k}n\PR_k$,
that is, the base change along $(n+1)\PRTT_k\to (n+1)\PR_k$.
Then there exists a canonical equivalence
$\HH^\bullet(n\Mod_{\AAA}[\PRTT_k]/k)\simeq \HH^\bullet(n\Mod_{\underline{\AAA}}[\PR_k]/k)$
in $\Alg_{n+2}(\Mod_k)$.

\item Let $\AAA$ be an object of $\CAlg(\PRTT_k)$. Suppose that $\AAA=\widehat{\AAA}$ belongs to $\PRTTT_k$ 
(when we regard $\AAA$ as an object of $\PRTTT_k$, we write $\widehat{\AAA}$ for $\AAA$).
Consider $n\Mod_{\AAA}[\PRTT_k]\in (n+1)\PRTT_k$ and $n\Mod_{\widehat{\AAA}}[\PRTTT_k]\in (n+1)\PRTTT_k$.
We then consider Hochschild cohomology $\HH^\bullet(n\Mod_{\AAA}[\PRTT_k]/k)$, obtained from $n\Mod_{\AAA}[\PRTT_k]\in (n+1)\PRTT_k$,
and Hochschild cohomology $\HH^\bullet(n\Mod_{\widehat{\AAA}}[\PRTTT_k]/k)$ obtained from $n\Mod_{\widehat{A}}[\PRTTT_k]\in (n+1)\PRTTT_k$.
Then there exists a canonical equivalence
$\HH^\bullet(n\Mod_{\AAA}[\PRTT_k]/k)\simeq \HH^\bullet(n\Mod_{\widehat{\AAA}}[\PRTTT_k]/k)$
in $\Alg_{n+2}(\Mod_k)$.
\end{enumerate}
\end{Proposition}

\Proof
The proof is similar to Proposition~\ref{comparison1}.
We prove (1).
Note first that $\End(n\Mod_{\AAA}[\PRTT_k])\otimes_{n\PRTT_k}n\PR_k\simeq \End(n\Mod_{\underline{\AAA}}[\PR_k])$.
The algebraization functor $n\PRTT_{\HH^\bullet(n\Mod_{\AAA}[\PRTT_k]/k)}\to \End(n\Mod_{\AAA}[\PRTT_k])$ (in $\Alg_1((n+1)\PRTT_k)$)
and the algebraization functor $n\PR_{\HH^\bullet(n\Mod_{\underline{\AAA}}[\PR_k]/k)}\to \End(n\Mod_{\underline{\AAA}}[\PR_k])$ (in $\Alg_1((n+1)\PR_k)$)
induce  
\begin{eqnarray*}
n\PR_{\HH^\bullet(n\Mod_{\AAA}[\PRTT_k]/k)}\simeq  n\PRTT_{\HH^\bullet(n\Mod_{\AAA}[\PRTT_k]/k)}\otimes_{n\PRTT_k}n\PR_k &\to&  n\PR_{\HH^\bullet(n\Mod_{\underline{\AAA}}[\PR_k]/k)}  \\
&\to& \End(n\Mod_{\underline{\AAA}}[\PR_k]).
\end{eqnarray*}
Then the composite  $n\PR_{\HH^\bullet(n\Mod_{\AAA}[\PRTT_k]/k)}\to n\PR_{\HH^\bullet(n\Mod_{\underline{\AAA}}[\PR_k]/k)}$
gives rise to a canonical morphism $\HH^\bullet(n\Mod_{\AAA}[\PRTT_k]/k)\to \HH^\bullet(n\Mod_{\underline{\AAA}}[\PR_k]/k)$ in $\Alg_{n+2}(\Mod_k)$.
As in Proposition~\ref{comparison1}, we are reduced to showing that
\[
\pi:\Hom_{\AAA\otimes S^{n-1}}(\AAA_{\Delta},\AAA_{\Delta}) \to \Hom_{\underline{\AAA}\otimes S^{n-1}}(\underline{\AAA_{\Delta}},\underline{\AAA_{\Delta}})  
\]
induces an equivalence of endomorphism algebras of the identity functors.
The right-hand side denotes the 
Hom object in $\PR_k$, which is defined in the same way as the left-hand side (the subscript $\underline{\AAA}\otimes S^{n-1}\in \CAlg(\PR_k)$ indicates the $\underline{\AAA}\otimes S^{n-1}$-linear).
By the proof of Lemma~\ref{higherendo}, $\pi$ can be identified with
\[
f:\lim_{[m]\in \Delta}\Hom(\mathcal{E}(m),\AAA)\to \lim_{[m]\in \Delta}\Hom(\underline{\mathcal{E}(m)},\underline{\AAA})
\]
where each $\underline{\mathcal{E}(m)}$ is the image of $\mathcal{E}(m)$ in $\PR_k$ (the left-hand side is the Hom object in $\PRTT_k$ discarding the $t$-structure,
and the right-hand is the Hom object in $\PR_k$). 
By the same argument as in Proposition~\ref{comparison1}, we see that $f$ 
induces an equivalence of endomorphism algebras of the identity functors.

The proof of (2) is similar to that of Proposition~\ref{comparison1} (2).
The same argument applies to this setting.
\QED

\subsection{Higher categorified sheaves}

\label{highercatsheaves}

Let $\XX$ be a perfect stack over $k$.
We will consider a higher version of Section~\ref{catsheaves}.
We take $\AAA$ to be $\QC(\mathcal{X})\in \CAlg(\PRTTT_k)$ (cf. Section~\ref{catsheaves}).
We consider $n\Mod_{\QC(\mathcal{X})}[\PRTTT_k]\in (n+1)\PRTTT_k$ and its Hochschild cohomology (see Proposition~\ref{highermodulehoch}).
The underlying stable $(\infty,n+1)$-category of 
$n\Mod_{\QC(\mathcal{X})}[\PRTTT_k]$ is $n\Mod_{\QC(\mathcal{X})}[\PR_k] \in (n+1)\PR_k$, which we take as a version of
higher categorified sheaves on $\XX$.
Here we will not discuss the affineness problem, which relates it to the sheaf $\textup{Aff}_{/\XX}^{op}\to (n+1)\PR_k$ given by $\mathcal{U}\mapsto n\PR_{R(\mathcal{U})}$.

The Hochschild cohomology $\HH^\bullet(n\Mod_{\QC(\mathcal{X})}[\PRTTT_k]/k)\in \Alg_{n+2}(\Mod_k)$ is defined to be the image of $\End_{\QC(\mathcal{X})\otimes S^{n-1}}(\QC(\mathcal{X}))$  under
$\Alg_{n+2}(E^+_{\PRTT_k})$ (cf. the proof of Proposition~\ref{highermodulehoch}).
According to Proposition~\ref{comparison2}, $\HH^\bullet(n\Mod_{\QC(\mathcal{X})}[\PRTTT_k]/k)$ is equivalent to 
Hochschild cohomology $\HH^\bullet(n\Mod_{\QC(\mathcal{X})}[\PR_k]/k)$
of $n\Mod_{\QC(\mathcal{X})}[\PR_k]\in (n+1)\PR_k$ when we regard $\QC(\XX)$ as an object of $\CAlg(\PR_k)$.

Let $m$ be a nonnegative integer.
Let $\Map(S^{m}, \mathcal{X})$ be the functor $\CAlg_k^{\textup{conn}}\to \SSS$, defined by $R\mapsto \Map(S^{m}, \mathcal{X}(R))$,
where $\mathcal{X}(R)$ is the space of $R$-valued points, and the sphere $S^{m}$ is regarded as an object of $\SSS$. Here $\CAlg_k^{\textup{conn}}$
is the full subcategory of $\CAlg_k$, spanned by connective spectra, and $\SSS$ is the $\infty$-category of small spaces.
This functor is represented by a perfect stack over $k$,
so that we think of $\Map(S^{m}, \mathcal{X})$ as a derived stack. Indeed, this stack is constructed as an iterated fiber product.
When $m=1$, $\Map(S^{m}, \mathcal{X})$ is the loop space $L\XX$.
This derived stack $\Map(S^{m}, \mathcal{X})$ admits a closed immersion $\iota_m:\XX\to \Map(S^{m}, \mathcal{X})$,
given by $\XX(R)\to \Map(S^m,\XX(R))$ determined by the assignment $x\mapsto c_x$,
where $c_x$ is the constant map taking the value $x$.

\begin{Proposition}
\label{highercatprop}
\begin{enumerate}
\item Let $\OO_\Delta$ in $\QC(\Map(S^n,\mathcal{X}))$ be the pushforward of the structure sheaf $\OO_{\XX}$ of $\XX$ along $\iota_n$.
Let $\mathcal{H}om_{\QC(\Map(S^n,\mathcal{X}))}(\OO_{\Delta}, \OO_{\Delta})$ in $\Mod_k$ denote the Hom complex.
Then there exists an equivalence 
\[
\HH^\bullet(n\Mod_{\QC(\mathcal{X})}[\PRTTT_k]/k)\simeq \HH^\bullet(n\Mod_{\QC(\mathcal{X})}[\PR_k]/k)\simeq \mathcal{H}om_{\QC(\Map(S^n,\mathcal{X}))}(\OO_{\Delta}, \OO_{\Delta})
\]
in
$\Mod_k$.

\item
Suppose that $k$ is of characteristic zero, and $\XX$ is an ordinary scheme separated and of finite type.
There exists an equivalence
\[
\HH^\bullet(n\Mod_{\QC(\mathcal{X})}[\PRTTT_k]/k)\simeq  \Gamma(\XX,\Sym_{\OO_{\XX}}^{\bullet}(\mathbb{T}_{\XX/k}[-n-1]))
\]
in $\Mod_k$, induced by Koszul duality (cf. the proof).
\end{enumerate}
\end{Proposition}

This proposition is a direct generalization of (the first part of) Proposition~\ref{catsheafhoch}.
The proof is also a rather straightforward generalization.

\Proof
We have already proved $\HH^\bullet(n\Mod_{\QC(\mathcal{X})}[\PRTTT_k]/k)\simeq \HH^\bullet(n\Mod_{\QC(\mathcal{X})}[\PR_k]/k)$.
According to \cite[Theorem 4.7, Corollary 4.12]{BFN}, there exists an equivalence
\[
\QC(\mathcal{X})\otimes S^{n-1}\simeq \QC(\Map(S^{n-1},\mathcal{X}))
\]
in $\CAlg(\PR_k)$.
Now we consider the underlying object of $\End_{\QC(\mathcal{X})\otimes S^{n-1}}(\QC(\mathcal{X}))$.
This is equivalent to $\Fun_{\QC(\Map(S^{n-1},\mathcal{X}))}(\QC(\mathcal{X}), \QC(\mathcal{X}))\in \PR_k$, which denotes the Hom object from $\QC(\mathcal{X})$ to itself.
Here $\QC(\mathcal{X})$ is the $\QC(\Map(S^{n-1},\mathcal{X}))$-module induced by the pullback functor $\iota_{n-1}^*:\QC(\Map(S^{n-1},\mathcal{X}))\to \QC(\mathcal{X})$
along the morphism $\iota_{n-1}:\mathcal{X}\to \Map(S^{n-1},\mathcal{X})$.
Since $\QC(\mathcal{X})$ is self-dual as an object of $\Mod_{\QC(\Map(S^{n-1},\mathcal{X}))}(\PR_k)$ (see \cite[Corollary 4.8]{BFN}), we have
\[
\Fun_{\QC(\Map(S^{n-1},\mathcal{X}))}(\QC(\mathcal{X}), \QC(\mathcal{X}))\simeq \QC(\mathcal{X})\otimes_{\QC(\Map(S^{n-1},\mathcal{X}))}\QC(\mathcal{X}).
\]
Moreover, by \cite[Theorem 4.7]{BFN} there exist equivalences
\[
\QC(\mathcal{X})\otimes_{\QC(\Map(S^{n-1},\mathcal{X}))}\QC(\mathcal{X})\simeq \QC(\mathcal{X}\times_{\Map(S^{n-1},\mathcal{X})}\mathcal{X})\simeq \QC(\Map(S^n,\mathcal{X}))
\]
in $\PR_k$. Through these equivalences, $\OO_\Delta$ corresponds to the identity functor  $\QC(\mathcal{X})\to \QC(\mathcal{X})$.
Then the underlying complex of the Hochschild cohomology is equivalent to
$\mathcal{H}om_{\QC(\Map(S^n,\mathcal{X}))}(\OO_{\Delta}, \OO_{\Delta})$.

Since $k$ is of characteristic zero, the affinization of $S^n$ is the trivial square zero extension $\Spec k\oplus k\epsilon_n$ of $\Spec k$ by $k\cdot \epsilon_n$
where $\epsilon_n$ is a ``formal" element of homological degree $n$.
Using this, we have an equivalence of derived schemes
$\Map(S^n,\mathcal{X})\simeq \Spec \Sym_{\OO_{\mathcal{X}}}^{\bullet}(\mathbb{L}_{\mathcal{X}}[n])$.
Then by Koszul duality \cite[Proof of Theorem 3.8, Theorem 4.3, Corollary 4.4]{H}, there exists an equivalence
$\mathcal{H}om_{\QC(\Spec\Sym_{\OO_{\mathcal{X}}}^{\bullet}(\mathbb{L}_{\mathcal{X}}[n]))}(\OO_{\Delta}, \OO_{\Delta})\simeq \Sym_{\OO_{\XX}}^{\bullet}(\mathbb{T}_{\XX/k}[-n-1])$
in $\QC(\mathcal{X})$. Taking global sections, our assertion follows.
\QED

\begin{Remark}
\label{cobor}
It is natural to ask if the equivalence of Proposition~\ref{highercatprop} (2) is promoted to an equivalence
of dg Lie algebras such that the underlying dg Lie algebra $\HH^\bullet(n\Mod_{\QC(\mathcal{X})}[\PRTTT_k]/k)[n+1]$ is obtained from
$\mathbf{E}_{n+2}$-algebra, which describes the deformation theory of $n\Mod_{\QC(\mathcal{X})}[\PRTTT_k]$ on commutative Artin algebras (see Section~\ref{deformsec}).
The right-hand side $T=\Gamma(\XX,\Sym_{\OO_{\XX}}^{\bullet}(\mathbb{T}_{\XX/k}[-n-1]))[n+1]$ also has a kind of Schouten bracket.
The work of Toën \cite[Section 5]{Toen} gives an equivalence between the $\mathbf{E}_{n+1}$-Hochschild cohomology $HH^{\mathbf{E}{n+1}}(\mathcal{X}/k)$ of $\QC(\mathcal{X})$ and $T$ as dg Lie algebras, where the bracket on $T$ is defined using the homology operad of the $\mathbf{E}_{n+1}$-operad. 
This construction also provides a link to the theory of shifted Poisson structures.
We note that the $\mathbf{E}_{n+2}$-algebra structure on $HH^{\mathbf{E}_{n+1}}(\mathcal{X}/k)$ in \cite{Toen} is defined by constructing
an $\mathbf{E}_{n+1}$-operad action on $S^n$ at the level of cobordisms (see Theorem 3.1 of {\it loc.\ cit.} for details).
Thus, the definition of the $\mathbf{E}_{n+2}$-algebra $HH^{\mathbf{E}_{n+1}}(\mathcal{X}/k)$ differs from the definition of Hochschild cohomology by means of universality (i.e., centers or algebraization/iterated endomorphisms),
though these will likely coincide:
the key step is a comparison of $\mathbf{E}_{n+1}$-monoidal structures on $\QC(\Map(S^n,\mathcal{X}))$ 
between that in \cite{Toen} and  ours (or one constructed as the center), which is beyond the scope of this paper.
We also remark that if one uses \cite{Toen}, then we have an equivalence in (2), at the level of $\Mod_k$ for a perfect stack $\XX$, which is a priori different from the above one. 
\end{Remark}


\section{Deformation theory}

\label{deformsec}

We will give applications to the deformation theory of $t$-structured stable presentable $(\infty,n)$-categories.

\subsection{Deformation theory}

\subsubsection{Deformations and Hochschild cohomology}
Let $R$ be an almost Artin $\mathbf{E}_{n+1}$-algebra (see Definition~\ref{defartin}).
Let $\DDD$ be an object of $n\PRTTT_k$.
Suppose that there exists the endomorphism algebra $\End(\DDD)\in \Alg_1(n\PRTTT_k)$. 

We define $\textup{Deform}_{\DDD}(R)^{\textup{cat}}\in \wCat$ to be the pullback of
\[
\xymatrix{
   & \{\DDD\}  \ar[d] \\
\RMod_{(n-1)\PRTTT_R}(n\PRTTT_k)  \ar[r]^(0.7){\beta^\pm}  & n\PRTTT_k
}
\]
where $\beta^\pm$ is the base change functor along 
$(n-1)\PRTTT_R \to (n-1)\PRTTT_k$.
We may regard an object of  $\textup{Deform}_{\DDD}(R)^{\textup{cat}}$ as 
\[
(\DDD_R, \phi:\beta^\pm(\DDD_R)=\DDD_R\widehat{\otimes}_{(n-1)\PRTTT_R}(n-1)\PRTTT_k\simeq \DDD)
\]
such that $\DDD_R$ is a $(n-1)\PRTTT_R$-module, and $\phi$ is an equivalence in $n\PRTTT_k$.
We define $\textup{Deform}_{\DDD}(R)=(\textup{Deform}_{\DDD}(R)^{\textup{cat}})^\simeq$ (cf. Section~\ref{notation}).
We shall refer to $\textup{Deform}_{\DDD}(R)$ as the space of deformations of $\DDD$ to $R$.
We will describe $\textup{Deform}_{\DDD}(R)$  in terms of Hochschild cohomology.

Let us recall
\[
\Phi_{R,n} :\RMod_{(n-1)\PRTT_R}(n\PRTT_k) \to  \LMod_{(n-1)\PRTT_{\DD_{n+1}(R)}}(n\PRTT_k) 
\]
Passing to the base change to $\PRTTT_k$,
we consider the commutative diagram 
\[
\xymatrix{
\RMod_{(n-1)\PRTTT_R}(n\PRTTT_k)  \ar[r]^{\widehat{\Phi}_{R,n}}_{\simeq} \ar[d]^{\beta^\pm} & \LMod_{(n-1)\PRTTT_{\DD_{n+1}(R)}}(n\PRTTT_k)   \ar[d]^{q^{\pm}} \\
n\PRTTT_k  \ar[r]^{\textup{id}}  & n\PRTTT_k
}
\]
where $q^\pm$ the forgetful functor and $\beta^\pm$ is the base change functor along 
$(n-1)\PRTTT_R \to (n-1)\PRTTT_k$.
By Theorem~\ref{HCKD}, $\widehat{\Phi}_{R,n}$ induces an equivalence.
It gives rise to an equivalence
\[
\widehat{\Phi}_{R,n}^{\DDD}:\RMod_{(n-1)\PRTTT_R}(n\PRTTT_k)\times_{\beta^{\pm},n\PRTTT_k}\{\DDD\}
\stackrel{\simeq}{\longrightarrow}  \LMod_{(n-1)\PRTTT_{\DD_{n+1}(R)}}(n\PRTTT_k) \times_{q^{\pm}, n\PRTTT_k}\{\DDD\}.
\]
We furthermore assume that $\DDD$ has Hochschild cohomology $\HH^\bullet(\DDD/k)\in \Alg_{n+1}(\Mod_k)$ (see Definition~\ref{Hochschilddef}).
By the universal property of $\End(\DDD)$ (\cite[6.1.2.40]{HA}) and the definition of $\HH^\bullet(\DDD/k)$ in terms of algebraizations,
we have equivalences of spaces
\begin{eqnarray*}
\big{(}\LMod_{(n-1)\PRTTT_{\DD_{n+1}(R)}}(n\PRTTT_k) \times_{q^\pm, n\PRTTT_k}\{\DDD\}\big{)}^{\simeq} &\simeq& \Map_{\Alg_1(\Mod_k)}((n-1)\PRTTT_{\DD_{n+1}(R)},\End(\DDD)) \\
&\simeq& \Map_{\Alg_{n+1}(\Mod_k)}(\DD_{n+1}(R),\HH^\bullet(\DDD/k))
\end{eqnarray*}
Combined with the equivalence $\widehat{\Phi}_{R,n}^{\DDD}$, we obtain an equivalence of spaces
\[
\textup{Deform}_{\DDD}(R)\simeq \Map_{\Alg_{n+1}(\Mod_k)}(\DD_{n+1}(R),\HH^\bullet(\DDD/k)).
\]
If $k\oplus \HH^\bullet(\DDD/k)\to k$ denotes the associated augmented $\mathbf{E}_{n+1}$-algebra,
then 
\[
\Map_{\Alg_{n+1}(\Mod_k)}(\DD_{n+1}(R),\HH^\bullet(\DDD/k))\simeq\Map_{\Alg_{n+1}^+(\Mod_k)}(\DD_{n+1}(R),k\oplus \HH^\bullet(\DDD/k)).
\]
Consequently, we see

\begin{Theorem}
\label{deformmain}
Suppose that $\DDD$ has Hochschild cohomology $\HH^\bullet(\DDD/k)$.
There exists a canonical equivalence 
\[
\textup{Deform}_{\DDD}(R)\simeq \Map_{\Alg_{n+1}^+(\Mod_k)}(\DD_{n+1}(R),k\oplus \HH^\bullet(\DDD/k)).
\]
\end{Theorem}

For instance, when $R=k\oplus k[j]$ is the trivial square extension by $k[j]$  ($j\ge0$), then
$\pi_0(\textup{Deform}_{\DDD}(R))\simeq \HH^{n}(\DDD/k)=H^{n+j+1}(\HH^\bullet(\DDD/k))$
since $\DD_{n+1}(R)$ is the free $\mathbf{E}_{n+1}$-algebra generated by $k[-n-j-1]$.

Next, we consider the deformation functor.

\begin{Definition}
\label{deffunctorcomplete}
Let 
$\theta_{k,n+1\to n}^\pm:\Alg_{n+1}(\Mod_k^{\le0})\to \Alg_1(n\PRTTT_k)$ be the functor given by $B\mapsto n\PRTTT_B$.
We note that $\theta_{k,n+1\to n}^\pm$ is fully faithful. We here write $\theta^\pm=\theta_{k,n+1\to n}^\pm$.
We set
\[
\textup{RLin}^{\wedge}(k):=\RMod(n\PRTTT_k)\times_{\Alg_1(n\PRTTT_k),\theta^\pm}\Alg_{n+1}(\Mod_k^{\le0}).
\]
where $\RMod(n\PRTTT_k)\to \Alg_1(n\PRTTT_k)$ is the forgetful functor (cf. Section~\ref{notation}).
Let
 $\textup{RLin}^\wedge(k)^\dagger$ be the subcategory spanned by coCartesian morphisms over $\Alg_{n+1}(\Mod_k^{\le0})$.
Let
$\textup{RLin}^\wedge(k)_{/(\DDD,k)}^\dagger$
be the overcategory where $(\DDD,k)$ means $\DDD\in n\PRTTT_k$ lying over $k\in \Alg_{n+1}(\Mod_k^{\le0})$.
We then have a left fibration $r:\textup{RLin}^\wedge(k)_{/(\DDD,k)}^\dagger\to \Alg_{n+1}^+(\Mod_k^{\le0})$.
Given $R\in \Alg_{n+1}^+(\Mod_k^{\le0})$ we refer to an object of the fiber $r^{-1}(R)$ over $R$ as a complete deformation of $\DDD$ to $R$.

Recall that $\SSS$ is the $\infty$-category of small spaces ($\infty$-groupoids),
 and $\wSSS$ is the $\infty$-category of $\mathbb{V}$-small spaces.

We define 
\[
\overline{\textup{Deform}}_{\DDD}:\Alg_{n+1}^+(\Mod_k^{\le0})\to \wSSS
\]
to be the functor which corresponds to the left fibration $r$. This is given informally by $A\mapsto \textup{Deform}_{\DDD}(A)$.
Let $\textup{Deform}_{\DDD}:\textup{Art}_{n+1}\to \wSSS$ be the restriction of $\overline{\textup{Deform}}_{\DDD}$ to $\textup{Art}_{n+1}$.
\end{Definition}

Next, we recall formal $\mathbf{E}_m$-moduli problems from \cite[Section 4]{DAG}.
Let $m$ be a natural number or $\infty$.
A functor $F:\textup{Art}_n\to \SSS$ from the $\infty$-category of Artin $\mathbf{E}_m$-algebras is a formal $\mathbf{E}_m$-moduli problem if
\begin{itemize}
  \item $F(k)$ is a contractible space,
\item for any pullback square \[
\xymatrix{
R \ar[r] \ar[d] &  R_0 \ar[d] \\
R_1 \ar[r] & R_{01}
}
\]
in $\textup{Art}_m$
such that $H^0(R_0)\to H^0(R_{01})$ and $H^0(R_1)\to H^0(R_{01})$ are surjective, the induced map
 \[
F(R) \to  F(R_0)  \times_{F(R_{01})}F(R_1)
\]
is an equivalence.
\end{itemize}
Let $\widehat{\mathsf{Stack}}_m^\ast$ be the full subcategory of $\Fun(\textup{Art}_m, \SSS)$, which consists of formal $\mathbf{E}_m$-moduli problems.
When $m$ is a natural number, thanks to a theorem of Lurie \cite[X, Theorem 4.0.8]{DAG} there exists an equivalence of $\infty$-categories
\[
\widehat{\mathsf{Stack}}_m^\ast \simeq \Alg_m^+(\Mod_k).
\]
Through the equivalence,
an augmented algebra $B\in \Alg_m^+(\Mod_k)$ corresponds to the formal $\mathbf{E}_m$-moduli problem $\mathcal{F}_{B}^{(m)}: \textup{Art}_m\to  \SSS$
defined informally by $R \mapsto \Map_{\Alg_m^+(\Mod_k)}(\DD_m(R),B)$.

\begin{Theorem}
\label{FMPmain}
The functor
$\textup{Deform}_{\DDD}$ is a formal $\mathbf{E}_{n+1}$-moduli problem $\mathcal{F}_{k\oplus \HH^\bullet(\DDD/k)}^{(n+1)}$, that is represented by
the augmented Hochschild cohomology $k\oplus\HH^\bullet(\DDD/k)\to k$ through $\widehat{\mathsf{Stack}}_{n+1}^\ast \simeq \Alg_{n+1}^+(\Mod_k)$.
\end{Theorem}

\Proof
It is enough to construct an equivalence
$\textup{Deform}_{\DDD}\simeq \mathcal{F}_{k\oplus \HH^\bullet(\DDD/k)}^{(n+1)}$.
It will suffice to construct an equivalence $\textup{Deform}_{\DDD}(R)\simeq \Map_{\Alg_{n+1}^+(\Mod_k)}(\DD_{n+1}(R),k\oplus \HH^\bullet(\DDD/k))$
in a functorial way with respect to $R\in \textup{Art}_{n+1}$.
This can be done by applying the procedure in \cite[X, Construction 5.3.18]{DAG}.
\QED

\subsubsection{Obstruction theory}
\label{obssec}
We consider obstructions for deformations.
Let $p:R'\to R$ be a morphism in $\textup{Art}_{n+1}$, which fits in the pullback square
\[
\xymatrix{
R' \ar[r]^p \ar[d] & R \ar[d]^u \\
k \ar[r] & k\oplus V[i] 
}
\]
in $\textup{Art}_{n+1}$, where $V$ is a finite dimensional vector space, $i$ is a positive integer,
and $k\oplus V[i]$ is the trivial square zero extension by $V[i]$.
We shall refer to $p$ as an elementary extension of $R$ by $V[i-1]$, associated to $u$.
Every map $R_1\to R_2$ in $\textup{Art}_{n+1}$ is the composition of a sequence of elementary extensions.

\begin{Theorem}
\label{obstructionmain}
Let $\DDD_R$ be a deformation of $\DDD$ to $R$, that is an object of $\textup{Deform}_{\DDD}(R)$. 
For the above square diagram, there exists a canonically defined element $\textup{Obs}(\DDD_R,p)\in \HH^{i+n+1}(\DDD/k)\otimes_kV$ 
such that $\textup{Obs}(\DDD_R,p)=0$ if and only if there exists $\DDD_{R'}\in \pi_0(\textup{Deform}_{\DDD}(R'))$ lying over $\DDD_R$ along  the map $p$.
Here $\HH^{i+n+1}(\DDD/k)$ is $H^{i+n+1}(\HH^\bullet(\DDD/k))$.
\end{Theorem}

\Proof
We note from Theorem~\ref{FMPmain} that there exists the canonical equivalence 
\[
\alpha:\textup{Deform}_{\DDD}(R')\simeq \textup{Deform}_{\DDD}(R)\times_{\textup{Deform}_{\DDD}(k\oplus V[i])}\textup{Deform}_{\DDD}(k).
\]
where $\textup{Deform}_{\DDD}(k)$ is a contractible space.
Since $\DD_{n+1}(k\oplus V[i])$ is a free $\mathbf{E}_{n+1}$-algebra generated by $V^{\vee}[-i-n-1]$ (\cite[X, Proposition 4.5.6]{DAG}),
there exists an equivalence $\phi:\pi_0(\textup{Deform}_{\DDD}(k\oplus V[i])) \simeq \HH^{i+n+1}(\DDD/k)\otimes_kV$.
Here, $V^\vee$ is the dual vector space.
Moreover, by construction,
$0\in \HH^{i+n+1}(\DDD/k)\otimes_kV$ corresponds to the trivial deformation, that is, the image of $\textup{Deform}_{\DDD}(k)\to \textup{Deform}_{\DDD}(k\oplus V[i])$.
We define $\textup{Obs}(\DDD_R,p)$ to be the image of $\DDD_R\in \pi_0(\textup{Deform}_{\DDD}(R))$
in $\HH^{i+n+1}(\DDD/k)\otimes_kV$ via $\phi$.
By the equivalence $\alpha$, we see that 
there exists some $\DDD_{R'}\in \textup{Deform}_{\DDD}(R')$ lying over $\DDD_R$ exactly when
the image of $\DDD_{R}$ in $\textup{Deform}_{\DDD}(k\oplus V[i])$ is equivalent to the trivial deformation.
This condition corresponds to the condition $\textup{Obs}(\DDD,p)=0$.
\QED

\subsection{Deformation of categories of Higher categorified sheaves}

We consider
$n\Mod_{\QC(\mathcal{X})}[\PRTTT_k]\in (n+1)\PRTTT_k$, which we thought of as the $(\infty,n+1)$-category of $t$-structured (higher) categorified sheaves.
Note that the underlying stable $(\infty,n+1)$-category is $n\Mod_{\QC(\mathcal{X})}[\PR_k]$.

\begin{Theorem}
\label{higherfinal}
Let $\XX$ be a perfect stack over $k$. 
Consider the deformation functor $\textup{Deform}_{n\Mod_{\QC(\mathcal{X})}[\PRTTT_k]}$.
Then the isomorphism classes of deformations of $n\Mod_{\QC(\mathcal{X})}[\PRTTT_k]$ to the trivial square zero extension $k\oplus k[i]$ ($i\ge0$) (i.e. $\pi_0(\textup{Deform}_{n\Mod_{\QC(\mathcal{X})}[\PRTTT_k]})(k\oplus k[i])$)
is isomorphic to
$H^{i+n+2}(\mathcal{H}om_{\QC(\Map(S^n,\mathcal{X}))}(\OO_{\Delta}, \OO_{\Delta}))$.

We further assume that $\XX$ is a scheme separated of finite type over $k$ of characteristic zero. 
Then there is an isomorphism
\[
\pi_0(\textup{Deform}_{n\Mod_{\QC(\mathcal{X})}[\PRTTT_k]})(k\oplus k[i])\simeq 
H^{i+n+2}(\Gamma(\XX,\Sym_{\OO_{\XX}}^{\bullet}(\mathbb{T}_{\XX/k}[-n-1]))).
\]
\end{Theorem}

\Proof
This assertion follows from Proposition~\ref{catsheafhoch}, Proposition~\ref{highercatprop} and Theorem~\ref{deformmain}
together with the fact that $\DD_{n+2}(k\oplus k[i])$ is a free $\mathbf{E}_{n+2}$-algebra generated by a shifted vector space $k[-i-n-2]$.
\QED


\begin{thebibliography}{99}



\bibitem{Ao}
K. Aoki,
Higher presentable categories and limits,
preprint available at arXiv:2510.13503

\bibitem{BSP}
C.~Barwick and C.~Schommer-Pries,
On the unicity of the theory of higher categories,
\emph{J.\ Amer.\ Math.\ Soc.} \textbf{34} (2021), no.~4, 1011--1058.

\bibitem{BBD}
A. Beilinson, J. Bernstein, and P. Deligne,
\emph{Faisceaux pervers},
Ast\'erisque, 100, Soci\'et\'e Math\'ematique de France, Paris, 1982.

\bibitem{BFN}
D. Ben-Zvi, J. Francis and D. Nadler,
Integral transforms in derived algebraic geometry,
\emph{J. Amer. Math. Soc.} (2010), 909--966.

\bibitem{BKP}
A. Blanc, L. Katzarkov, and P. Pandit,
Generators in formal deformations of categories,
\emph{Compos. Math.} \textbf{154.10} (2018).

\bibitem{BT}
L. Brantner and L. Taelman,
Deformations and lifts of Calabi-Yau varieties in characteristic $p$,
preprint available at arXiv:2407.09256

\bibitem{Ch}
Y. Chen,
Deformation theory of $\mathbb{E}_n$-monoidal categories,
preprint available at arXiv:2602.21431

\bibitem{FR}
J. Francis,
The tangent complex and Hochschild cohomology of $\mathbb{E}_n$-rings,
\emph{Compos. Math.} \textbf{149} (2013), no. 3, 430--480.

\bibitem{Ind}
D. Gaitsgory,
Ind-coherent sheaves,
\emph{Mosc. Math. J.}, (2013), Vol. 13, 399--528.

\bibitem{Gai1aff}
D. Gaitsgory,
Sheaves of categories and the notion of $1$-affineness,
\emph{Contemp. Math.} \textbf{643} (2015), 127--225.

\bibitem{Gai2}
D. Gaitsgory and N. Rozenblyum,
A study in Derived Algebraic Geometry II: Deformations, Lie Theory and Formal Geometry,
Mathematical Surveys and Monographs, 222, American Math. Soc., 2017.

\bibitem{Ger1}
M. Gerstenhaber,
On the deformation of rings and algebras,
\emph{Ann. of Math.} \textbf{79} (1964), 59--103.

\bibitem{Ger2}
M. Gerstenhaber,
The cohomology structure of an associative ring,
\emph{Ann. of Math.} \textbf{78} (1963), 267--288.

\bibitem{GH}
D. Gepner and R. Haugseng,
Enriched $\infty$-categories via non-symmetric operads,
\emph{Adv. Math.} \textbf{279} (2015), 575--716.

\bibitem{Hinen}
V. Hinich,
Yoneda lemma for enriched $\infty$-categories,
\emph{Adv. Math.} \textbf{400} (2022), 108268.

\bibitem{H}
I. Hisamichi,
Koszul duality for symmetric algebras and derived loop spaces,
preprint available at arXiv:2609.18115

\bibitem{ID}
I. Iwanari,
Differential calculus of Hochschild pairs for infinity-categories,
SIGMA $\mathbf{16}$ (2020), 97  57 pages,
Special Issue on Primitive Forms and Related Topics in honor of Kyoji Saito for his 77th birthday.

\bibitem{IM}
I. Iwanari,
Moduli theory associated to Hochschild pairs,
Trans. Amer. Math. Soc. to appear (61 pages); available at the author's webpage https://sites.google.com/site/isamuiwanarishomepage/

\bibitem{AL}
I. Iwanari,
Algebraic Analysis of linear categories,
available at the author's webpage

\bibitem{IH}
I. Iwanari,
Hochschild cohomology and deformation theory of stable infinity-categories with $t$-structures
available at the author's webpage or arXiv:2506.21867

\bibitem{CKDA}
I. Iwanari,
Categorified Koszul duality of algebras,
available at the author's webpage or arXiv:2512.12263

\bibitem{JFS}
T.~Johnson-Freyd and C.~Scheimbauer,
(Op)lax natural transformations, twisted quantum field theories, and ``even higher'' Morita categories,
\emph{Adv.\ Math.} \textbf{307} (2017), 147--223.

\bibitem{HTT}
J. Lurie,
Higher Topos Theory,
Annals of Math. Studies, 2009.

\bibitem{HA}
J. Lurie,
Higher Algebra,
Draft 2017 September, available at Lurie's webpage.

\bibitem{DAG}
J. Lurie,
Derived Algebraic Geometry Series,
available at Lurie's webpage.

\bibitem{Lu}
J.~Lurie,
On the classification of topological field theories,
\emph{Current Developments in Mathematics, 2008},
Int.\ Press, Somerville, MA, 2009, 129--280.

\bibitem{MR}
I. Mirkovic and S. Riche,
Linear Koszul duality,
\emph{Compositio Math.} \textbf{146} (2010), 233--258.

\bibitem{RZ}
D. Reutter and M. Zetto,
Enriched $\infty$-categories as marked module categories,
preprint available at arXiv:2501.07697

\bibitem{Si}
C.~Simpson,
\emph{Homotopy theory of higher categories},
New Mathematical Monographs, vol.~19,
Cambridge University Press, Cambridge, 2012.

\bibitem{Ste}
G. Stefanich,
Presentable $(\infty,n)$-categories,
preprint available at arXiv:2011.03035

\bibitem{TV}
B. To\"{e}n and G. Vezzosi,
Chern character, loop spaces and derived algebraic geometry,
In: Algebraic topology. Vol. 4.
Abel Symp. Berlin: Springer, 2009, pp. 331--354.

\bibitem{Toen}
B. To\"{e}n,
Operations on derived moduli spaces of branes,
preprint available at arXiv:1307.0405







\end{thebibliography}
\end{document}